\documentclass[oneside,12pt]{amsart}
\usepackage{amsmath,amsfonts,euscript,amscd,amsthm,amssymb,upref,graphics}

\theoremstyle{plain}
\newtheorem{problem}{Problem}

\swapnumbers

\theoremstyle{plain}

\newtheorem{theorem}[subsection]{Theorem}

\newtheorem{proposition}[subsection]{Proposition}
\newtheorem{lemma}[subsection]{Lemma}
\newtheorem{lemmadef}[subsection]{Lemma and Definition}
\newtheorem{corollary}[subsection]{Corollary}

\newtheorem{nothing}[subsection]{}

\newtheorem{subtheorem}{Theorem}[subsection]
\newtheorem{subproposition}[subtheorem]{Proposition}
\newtheorem{sublemma}[subtheorem]{Lemma}

\newtheorem{subnothing}[subtheorem]{}

\theoremstyle{definition}

\newtheorem{definition}[subsection]{Definition}

\newtheorem{nothing*}[subsection]{}
\newtheorem{example}[subsection]{Example}

\newtheorem{notation}[subsection]{Notation}

\newtheorem{subdefinition}[subtheorem]{Definition}

\newtheorem{subnothing*}[subtheorem]{}

\theoremstyle{remark}

\newtheorem*{remark}{Remark}

{\refstepcounter{subsection}%
{\par\smallskip\noindent\thesubsection\bf.~#1.}\it}{\par}

\newenvironment{gen*}[1]%
{\refstepcounter{subsection}
{\par\smallskip\noindent\thesubsection\bf.~#1.}\rm}{\par}

{{\par\smallskip\noindent\bf #1.}\it}{\par}

\newenvironment{gennonumber*}[1]%
{{\par\smallskip\noindent\bf #1.}\rm}{\par}

{\refstepcounter{subtheorem}%
{\par\smallskip\noindent\thesubtheorem\bf.~#1.}\it}{\par}

\newenvironment{subgen*}[1]%
{\refstepcounter{subtheorem}
{\par\smallskip\noindent\thesubtheorem\bf.~#1.}\rm}{\par}

{\begin{enumerate}\setlength{\itemsep}{#1}}
{\end{enumerate}}

{\begin{itemize}\setlength{\itemsep}{#1}}
{\end{itemize}}

\newcommand{\Aut}{	\operatorname{{\rm Aut}}}

\newcommand{\Vtx}{	\operatorname{{\rm Vtx}}}
\newcommand{\Sub}{	\operatorname{{\rm Sub}}}
\newcommand{\SUB}{	\operatorname{\overline{\Sub}}}
\newcommand{\hodge}[1]{\| {#1} \| }

\newcommand{\For}{	\operatorname{\mbox{\rm FO}^+}}
\newcommand{\RFSK}{	\operatorname{\mbox{\rm FO}}}

\newcommand{\transp}{\operatorname{\mbox{\rm Transp}}}
\newcommand{\SEC}{\mbox{$\Integ^*/\!\!\sim$}}

\newcommand{\setspec}[2]{\big\{\,#1\, \mid \,#2\, \big\}}

\newlength{\mylength}
\newcommand{\Integ}{\ensuremath{\mathbb{Z}}}
\newcommand{\Nat}{\ensuremath{\mathbb{N}}}

\newcommand{\Reals}{\ensuremath{\mathbb{R}}}

\newcommand{\proj}{\ensuremath{\mathbb{P}}}

\newcommand{\Xgoth}{{\ensuremath{\mathfrak{X}}}}

\newcommand{\Beul}{\EuScript{B}}
\newcommand{\Ceul}{\EuScript{C}}

\newcommand{\Geul}{\EuScript{G}}
\newcommand{\Heul}{\EuScript{H}}

\newcommand{\Leul}{\EuScript{L}}
\newcommand{\Meul}{\EuScript{M}}
\newcommand{\Neul}{\EuScript{N}}
\newcommand{\Oeul}{\EuScript{O}}
\newcommand{\Peul}{\EuScript{P}}

\newcommand{\Reul}{\EuScript{R}}

\newcommand{\Xeul}{\EuScript{X}}
\newcommand{\Yeul}{\EuScript{Y}}

\newcommand{\isom}{\cong}
\renewcommand{\epsilon}{\varepsilon}
\renewcommand{\phi}{\varphi}
\renewcommand{\emptyset}{\varnothing}
\newcommand{\emptyseq}{\varnothing}

\newcommand{\rien}[1]{}

\newcommand{\simref}[1]{\overset{\,\mbox{\tiny #1} \,}{\sim}}
\newcommand{\tequiv}[1]{\overset{\,\mbox{\tiny $#1$} \,}{\sim}}

\newcommand{\GG}{\Geul}
\newcommand{\LL}{\Leul}

\newcommand{\Gn}[1]{\Geul_{#1}}
\newcommand{\HH}{\Heul}

\newcommand{\Empty}{\mbox{\Large\bf\o}}
\newcommand{\skm}{\dashrightarrow}
\newcommand{\lskm}{\dashleftarrow}

\newcommand{\DR}{\cite{DaiRuss:Rul1}}

\begin{document}
\renewcommand{\baselinestretch}{1.07}
\setlength{\unitlength}{1mm}


\title[Classification of weighted graphs]%
{Classification of weighted graphs\\
up to blowing-up and blowing-down}

\author{Daniel Daigle}

\address{Department of Mathematics and Statistics\\
University of Ottawa\\
Ottawa, Canada\ \ K1N 6N5}

\email{ddaigle@uottawa.ca}

\begin{abstract}
We classify weighted forests up to the
blowing-up and blowing-down operations which are relevant for the 
study of algebraic surfaces.
\end{abstract}

\maketitle

The word ``graph'' in this text means a finite
undirected graph such that no edge connects a vertex to itself
and at most one edge joins any given pair of vertices.
A {\em weighted graph\/} is a graph in which each vertex
is assigned an integer (called its weight).
Two operations are performed on weighted graphs:
The blowing-up and its inverse, the blowing-down.
Two weighted graphs are said to be equivalent if one can be obtained
from the other by means of a finite sequence of 
blowings-up and blowings-down (see \ref{FromHere}--\ref{ToHere}).

These weighted graphs and operations are well known to geometers
who study algebraic surfaces.
Many problems in the geometry
of surfaces can be formulated in graph-theoretic terms
and solving these sometimes requires elaborate graph-theoretic
considerations. This gives rise to a variety of questions about
weighted graphs, all in connection with
the equivalence relation generated by blowing-up and blowing-down.

The present paper proposes a classification of weighted forests
up to equivalence.
In particular, Theorem~\ref{MainRFSK} defines an invariant $\bar Q(\GG)$
for any pseudo-minimal (\ref{LemDefRegForest}) weighted forest $\GG$,
and asserts that $\bar Q(\GG) = \bar Q(\GG')$ if and only if
$\GG$ is equivalent to $\GG'$.
Since $\bar Q(\GG)$ can actually be computed, this yields an algorithm
for deciding whether two weighted forests are equivalent 
(see Problem~\ref{DecisionProblem},
at the end of section~\ref{Section:RegForSkel}).
Apparently this decision problem was previously open,
even in the special case of ``linear chains'',
i.e., weighted graphs of the form:
$$
\setlength{\unitlength}{1mm}
\begin{picture}(32,4)(-1,-.5)
	\put(0,0){\line(1,0){14}}
	\put(0,0){\circle*{1}}
	\put(10,0){\circle*{1}}
	\put(20,0){\makebox(0,0)[c]{\dots}}
	\put(30,0){\circle*{1}}
	\put(30,0){\line(-1,0){4}}
	\put(0,1.5){\makebox(0,0)[b]{$\scriptstyle x_1$}}
	\put(10,1.5){\makebox(0,0)[b]{$\scriptstyle x_2$}}
	\put(30,1.5){\makebox(0,0)[b]{$\scriptstyle x_q$}}
\end{picture}\qquad (x_i\in\Integ).
$$
Note that, in the case of linear chains,
\ref{MainRFSK} simplifies to \ref{Isoms}.

We also contribute to the problem of listing all minimal elements in a
given equivalence class of weighted forests.
Section~\ref{Sec:MinimalReduc} reduces that problem to the
case of linear chains.
This special case is given a recursive solution in
Section~\ref{Sec:FurtherLin} and, in some simple cases, an explicit solution.
Incidentally, the cases that we are able to describe explicitely
are precisely those which arise from the study of algebraic surfaces.

Section \ref{Sec:StrictEqu} is concerned with
topological properties of graphs, but topology is never mentioned.
For instance, \ref{SameSkel} has the following consequence:
\begin{quote}
\it
Let $\GG, \GG'$ be weighted trees which are minimal and equivalent.
If $\GG$ is not a linear chain then the two trees are homeomorphic.
\end{quote}
Note, however, that if $n$ is a positive integer then the linear chains
$$
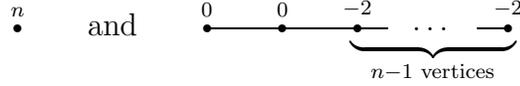
\setlength{\unitlength}{1mm}
\begin{picture}(2,9)(-1,-6)
	\put(0,0){\circle*{1}}
	\put(0,1.5){\makebox(0,0)[b]{$\scriptstyle n$}}
\end{picture}
\qquad \raisebox{5\unitlength}{\text{and}} \qquad
\begin{picture}(42,9)(-1,-6)
	\put(0,0){\line(1,0){24}}
	\put(40,0){\line(-1,0){4}}
	\put(0,0){\circle*{1}}
	\put(10,0){\circle*{1}}
	\put(20,0){\circle*{1}}
	\put(30,0){\makebox(0,0)[c]{\dots}}
	\put(40,0){\circle*{1}}
	\put(0,1.5){\makebox(0,0)[b]{$\scriptstyle 0$}}
	\put(10,1.5){\makebox(0,0)[b]{$\scriptstyle 0$}}
	\put(20,1.5){\makebox(0,0)[b]{$\scriptstyle -2$}}
	\put(40,1.5){\makebox(0,0)[b]{$\scriptstyle -2$}}
\put(30,-1.5){\makebox(0,0)[t]%
{$\displaystyle\underbrace{\rule{22\unitlength}{0mm}}_{n-1 \text{ vertices}}$}}
\end{picture}
$$
are minimal and equivalent, but not homeomorphic.
Because of this irregularity, and for other reasons as well,
the notions of skeleton and skeletal map (Section~\ref{SecSkel})
are better suited than topology for our purpose.
The ``topological'' results of section \ref{Sec:StrictEqu}
are of fundamental importance for
classifying weighted trees and forests (Section~\ref{Section:RegForSkel}),
but are not needed for the special case of linear chains.  
Readers only interested in that case may restrict themselves to sections
\ref{Sec:BasicDef},
\ref{SectionSequences},
\ref{Sec:ClassLinChains} and \ref{Sec:FurtherLin}.\footnote{%
In fact one needs \ref{AllLinChains} to prove \ref{BasicBup},
but \ref{AllLinChains} can be
proved independently from sections \ref{SecSkel} and~\ref{Sec:StrictEqu}.
Some ideas from Section~\ref{Section:TauEquiv}
are used in the proofs of Section~\ref{Sec:FurtherLin},
but not in a very crucial way.}

\medskip
\noindent\textbf{Acknowledgements.}
Papers \cite{Neumann:calculus} and \cite{Neumann:bil} classify
weighted forests up to an equivalence relation weaker than
the one considered here (the relation is generated
by blowing-up, blowing-down and other operations which are not
allowed here).
Result 3.2.1 of \cite{Rus:formal} classifies linear chains
but, again, this is relative to a weak equivalence relation.
Paper~\cite{Shastri:FinFundGp} uses the same equivalence relation as
we do, but only classifies a restricted class of weighted trees.

Proposition~3.2 of \cite{Rus:formal} almost%
\footnote{One also needs \ref{SubSeq} for the proof.}
implies the fact (\ref{EveryCan}) that each linear chain is 
equivalent to at least one canonical chain.
As we realized \textit{a posteriori,}
there is even some similarity between the cited result
and our method for proving \ref{EveryCan}.
We also noticed \textit{a posteriori\/} a certain resemblance between
our skeletons equipped with extra structure
(\ref{DefSkel}, \ref{DefSWw}, \ref{DefN})
and the ``W-graphs'' briefly described in \cite{Neumann:bil}.

\section{Basic definitions and facts}
\label{Sec:BasicDef}

\begin{nothing*}\label{NotaSequences}
If $E$ is a set then $E^* = \cup_{n=0}^\infty E^n$
denotes the set of finite sequences in $E$, including the empty
sequence $\emptyseq \in E^*$.
We write $A^-$ for the \textit{reversal\/} of $A\in E^*$,
i.e., if $A=(a_1,\dots,a_n)$ then $A^- =(a_n,\dots,a_1)$.
\end{nothing*}

\begin{nothing*}
If $\GG$ is a weighted graph, $\Vtx(\GG)$ is its vertex set.
If $v \in \Vtx(\GG)$ then
$w(v,\GG)$ denotes the weight of $v$ in $\GG$;
$\deg(v,\GG)$ denotes the degree of $v$ in $\GG$, that is,
the number of neibhors of $v$.
The empty graph is denoted $\Empty$.
\end{nothing*}

\begin{nothing*}\label{DefLinChain}\label{DefAdm}
A {\em linear chain\/} is a weighted tree
in which every vertex has degree at most two.
Given $(x_1, \dots, x_q) \in \Integ^*$, the linear chain
\begin{equation*}
\setlength{\unitlength}{1mm}
\begin{picture}(32,6)(-1,-2)
	\put(0,0){\line(1,0){14}}
	\put(0,0){\circle*{1}}
	\put(10,0){\circle*{1}}
	\put(20,0){\makebox(0,0)[c]{\dots}}
	\put(30,0){\circle*{1}}
	\put(30,0){\line(-1,0){4}}
	\put(0,1.5){\makebox(0,0)[b]{$\scriptstyle x_1$}}
	\put(10,1.5){\makebox(0,0)[b]{$\scriptstyle x_2$}}
	\put(30,1.5){\makebox(0,0)[b]{$\scriptstyle x_q$}}
\end{picture}
\end{equation*}
is denoted $[x_1, \dots, x_q]$.
So we distinguish between the graph $[x_1, \dots, x_q]$ 
and the sequence $(x_1, \dots, x_q)$ and we note that
$[x_1, \dots, x_q] = [x_q, \dots, x_1]$.

An {\em admissible chain} is a linear chain in which every weight
is strictly less than $-1$.  The empty graph $\Empty$ is an admissible chain.
\end{nothing*}

\begin{nothing*}\label{FromHere}
Let $\GG$ be a weighted graph. We define three types of ``blowing-up 
of $\GG$'':

\begin{enumerate}

\item If $v$ is a vertex of $\GG$ then the
{\em blowing-up of $\GG$ at $v$\/}
is the weighted graph $\GG'$ obtained from $\GG$ by adding
one vertex $e$ of weight $-1$, adding one edge joining $e$ to $v$,
and decreasing the weight of $v$ by $1$.
(This process is called a blowing-up ``at a vertex''.)

\item If $\epsilon = \{v_1,v_2\}$ is an edge of $\GG$
(so $v_1,v_2$ are distinct vertices of $\GG$),
then the {\em blowing-up of $\GG$ at $\epsilon$\/}
is the weighted graph $\GG'$ obtained from $\GG$ by adding
one vertex $e$ of weight $-1$,
deleting the edge
$\epsilon = \{v_1,v_2\}$, adding the two edges
$\{v_1,e\}$ and $\{e,v_2\}$, and 
decreasing the weights of $v_1$ and $v_2$ by $1$.
(This is called a blowing-up ``at an edge'', or a ``subdivisional''
blowing-up.)

\item The {\em free blowing-up of $\GG$\/} is the weighted graph
$\GG'$ obtained by taking the disjoint union of $\GG$ and of a
vertex $e$ of weight $-1$.

\end{enumerate}
In each of the above three cases, we call $e$ the 
vertex {\em created\/} by the blowing-up.
If $\GG'$ is a blowing-up of $\GG$
then there is a natural way to identify $\Vtx(\GG)$
with a subset of $\Vtx(\GG')$ (whose complement is $\{e\}$).
It is understood that, whenever a blowing-up is performed,
such an injective map $\Vtx(\GG)\hookrightarrow \Vtx(\GG')$ is chosen.
We stress that if $\GG'$ is a blowing-up of $\GG$ and 
$\GG''$ is a weighted graph isomorphic to $\GG'$, then
$\GG''$ is a blowing-up of $\GG$.
\end{nothing*}

\begin{nothing*}
A vertex $e$ of a weighted graph $\GG'$ is said to be
\textit{contractible\/} if the following three conditions hold:
(i) $e$ has weight $-1$;
(ii) $e$ has at most two neighbors;
(iii) if $v_1$ and $v_2$ are distinct neighbors of $e$
then $v_1,v_2$ are not neighbors of each other.

If $e$ is a contractible vertex of $\GG'$ then $\GG'$ is the
blowing-up of some weighted graph $\GG$ in such a way that $e$
is the vertex created by this process.  Up to isomorphism of weighted
graphs, $\GG$ is uniquely determined by $\GG'$ and $e$.
We say that $\GG$ is
obtained by  \textit{blowing-down\/} $\GG'$ at $e$.
The blowing-down is the inverse operation of the blowing-up.
\end{nothing*}

\begin{nothing*}\label{DefMinWG}
A  weighted graph is \textit{minimal\/} if it does not have
a contractible vertex.
\end{nothing*}

\begin{nothing*}\label{ToHere}
Two weighted graphs $\GG$ and $\HH$ are {\em equivalent\/}
(notation: $\GG \sim \HH$) if one can be obtained from
the other by a finite sequence of blowings-up and blowings-down.
\end{nothing*}

\begin{nothing*}\label{DefBilForm}
Given a weighted graph $\GG$,
consider the real vector space $V$ with basis $\Vtx(\GG)$
and define a symmetric bilinear form
$ B_\GG : V \times V \to \Reals$ by:
$$
B_\GG (u,v) = \begin{cases}
w(u, \GG),	& \text{if $u=v\in\Vtx(\GG)$}, \\
1,			& \text{if $u,v\in\Vtx(\GG)$ are distinct and joined by an edge}, \\
0,			& \text{if $u,v\in\Vtx(\GG)$ are distinct and
			not joined by an edge}.
\end{cases}
$$
One calls $B_\GG$ the \textit{intersection form\/} of $\GG$.
Then define the natural number
$\hodge{ \GG } = \displaystyle \max_{W} \dim W$,
where $W$ runs in the set of subspaces of $V$ satisfying
$$
\forall_{x\in W}\ \ B_\GG(x,x) \ge 0.
$$
Note that  $\hodge{ \GG } = 0$ iff $B_\GG$ is negative definite,
in which case we say that $\GG$ is negative definite.
\end{nothing*}

\begin{nothing}\label{HodgeInvar}
For weighted graphs $\GG$ and $\GG'$, 
$\GG \sim \GG' \implies \hodge{ \GG } = \hodge{ \GG' } $.
\end{nothing}

\begin{proof}
See for instance 1.14 of \cite{Rus:formal}.
\end{proof}

\begin{nothing*}\label{DefDet}
Consider a weighted graph $\GG$ and its intersection form
$B_\GG : V\times V\to \Reals$ (see \ref{DefBilForm}).
Let $v_1,\dots,v_n$ be the distinct vertices of $\GG$
(enumerated in any order) and let $M$ be the $n\times n$
matrix representing $B_\GG$ with respect to the basis
$(v_1,\dots,v_n)$ of $V$. That is,
$M_{ii}=w(v_i,\GG)$ and, if $i\neq j$,
$M_{ij}=1$ (resp.\ $0$) if $v_i,v_j$ are neighbors (resp.\ are not neighbors)
in $\GG$.
Note that $\det(-M)$ is independent of the choice of an ordering
for $\Vtx(\GG)$.  One defines the
{\em determinant\/} of the weighted graph $\GG$ by:
$$
\det(\GG)=\det(-M).
$$
Note that $\det(\GG)\in\Integ$.
By convention, $\det(\Empty) = 1$.
\end{nothing*}

The following is well-known, and easily verified:

\begin{nothing}\label{DetInvar}
For weighted graphs $\GG$ and $\GG'$, 
$\GG \sim \GG' \implies \det(\GG) = \det(\GG') $.
\end{nothing}

\begin{remark}
Without the minus sign in $\det(-M)$, \ref{DetInvar} would
only be true up to sign.
\end{remark}

\section{Skeletons and skeletal maps}
\label{SecSkel}

In this section,
all graphs are forests and (except in \ref{SkOfWeighted})
no graph is weighted.

\begin{definition}\label{DefSkel}
A \textit{skeleton\/} is a forest which contains no vertex of degree
zero or two.
\end{definition}

\begin{definition}\label{DefPT}
Given a forest $G$,
let $P(G)$ be the set of nonempty finite sequences
$\gamma=(v_0,\dots,v_n)$ of vertices of $G$ satisfying either
\begin{enumerate}

\item $n=0$ and  $\deg(v_0,G)=0$; or

\item $n>0$ and the following hold:
\begin{enumerate}

\item $v_0,\dots,v_n$ are distinct

\item $\{ v_0, v_1\}$,  $\dots$, $\{ v_{n-1}, v_n\}$ are
edges in $G$

\item $\setspec{ i }{ \deg(v_i,G) \neq 2 } = \{ 0, n \}$.

\end{enumerate}
\end{enumerate}
Note that $P(G)$ is the empty set if and only if $G = \Empty$.
If $\gamma = (v_0,\dots,v_n) \in P(G)$, then
$\gamma^- = (v_n,\dots,v_0)$ also belongs to $P(G)$.

Each element
$\gamma = (v_0,\dots,v_n)$ of $P(G)$ is of one of four types,
defined as follows:
\begin{itemize}

\item[] $\gamma$ is of type $(-,-)$ if $\deg(v_0,G)<2$ and $\deg(v_n,G)<2$;
\item[] $\gamma$ is of type $(+,-)$ if $\deg(v_0,G)>2$ and $\deg(v_n,G)<2$;
\item[] $\gamma$ is of type $(-,+)$ if $\deg(v_0,G)<2$ and $\deg(v_n,G)>2$;
\item[] $\gamma$ is of type $(+,+)$ if $\deg(v_0,G)>2$ and $\deg(v_n,G)>2$.
\end{itemize}
Observe that if $n=0$ then $\gamma$ is of type $(-,-)$.
\end{definition}

\begin{notation}\label{Nota:VtxCond}
Let $G$ be a graph and let $\Peul$ be a condition on the degree of a
vertex; then we define:
$$
\Vtx_{\Peul}(G) = \setspec{ x\in \Vtx(G) }
{ \deg( x, G ) \mbox{ satisfies } \Peul }.
$$
For instance,
$\Vtx_{\neq2}(G)=
\setspec{ x\in \Vtx(G) } { \deg( x, G ) \neq2 }$.
\end{notation}

\begin{definition}\label{DefSkelMap}
Let $G$ and $G'$ be forests.
\begin{enumerate}

\item A \textit{pre-skeletal map\/} from $G$ to $G'$ is a set map
$$
f : \Vtx_{\neq2}(G) \to \Vtx_{\neq2}(G')
$$
such that, given any $\gamma = (v_0,\dots, v_m) \in P(G)$,
there exists an element $\gamma' = (v_0',\dots, v_n')$ of  $P(G')$
satisfying $f(v_0)=v_0'$ and $f(v_m)=v_n'$.
Then $\gamma \mapsto\gamma'$ is
a well-defined map which we denote $\vec{f} : P(G)\to P(G')$.
This map satisfies $\vec{f}(\gamma^-) = \vec{f}(\gamma)^-$ for all
$\gamma \in P(G)$.

\item A \textit{skeletal map\/} from $G$ to $G'$ is a 
pre-skeletal map $f$ from $G$ to $G'$ which satisfies
the following additional conditions:
\begin{enumerate}

\item $\vec{f} : P(G)\to P(G')$ is surjective

\item if $u,v$ are distinct elements of the domain of $f$
satisfying $f(u)=f(v)$,
then there exists $(v_0,\dots, v_m) \in P(G)$ such that 
$v_0=u$ and $v_m=v$.

\end{enumerate}
\end{enumerate}
We write $f:G \skm G'$ to indicate that $f$ is a skeletal map from $G$ to $G'$.
The composition of two skeletal maps is a skeletal map.
If $G$ is a forest then the identity map on the set
$\Vtx_{\neq2}(G)$ is a skeletal map from $G$ to itself.
\end{definition}

%
%
%
%
%
%
We leave it to the reader to verify:

\begin{nothing}\label{PreservesType}
If $f: G \skm G'$ (where $G$ and $G'$ are forests)
then $\vec f : P(G) \to P(G')$ preserves type. That is,
if $\gamma \in P(G)$ then $\gamma$ and $\vec f(\gamma)$ are of the
same type (see \ref{DefPT}).
\end{nothing}

\begin{lemma}\label{Basic:Skel}
\begin{enumerate}

\item Suppose that
$G \overset{f}{\skm} G' \overset{\sigma}{\lskm} S$
are skeletal maps, where $G$, $G'$ are forests and $S$ is a skeleton.
Then $\sigma$ factors through $f$, i.e., there exists
$\theta : S \skm G$ such that $f\circ\theta = \sigma$.

\item If $S$ and $S'$ are skeletons
and $f: S \skm S'$ is a skeletal
map then $f$ is actually an isomorphism of graphs.

\end{enumerate}
\end{lemma}

\begin{proof}
Let $G, G', S, f$ and $\sigma$ be as in statement (1).
The definition of skeletal map implies that 
$f: \Vtx_{\neq2}(G) \to \Vtx_{\neq2}(G')$ is
surjective and is almost injective:
\begin{itemize}

\item[(i)] If $z\in \Vtx_{\neq2}(G')$ is such that 
$| f^{-1}(z) | > 1$, then $\deg(z,G')=0$ and $| f^{-1}(z) | = 2$.

\item[(ii)] If $z\in \Vtx_{=0}(G')$ is such that 
$ f^{-1}(z) = \{ v \}$, then $\deg(v,G)=0$.

\end{itemize}
Of course, $\sigma$ has similar properties. So, if we define
$Z=\Vtx_{=0}(G')$, $f$ and $\sigma$ restrict to bijections:
$$
f^{-1} \big( \Vtx_{\neq2}(G') \setminus Z \big)
\overset{f_1}{\longrightarrow} 
\Vtx_{\neq2}(G') \setminus Z
\overset{\sigma_1}{\longleftarrow} 
\sigma^{-1} \big( \Vtx_{\neq2}(G') \setminus Z \big),
$$
so we may define the bijection
$\theta_1 = f_1^{-1} \circ \sigma_1 : 
\sigma^{-1} \big( \Vtx_{\neq2}(G') \setminus Z \big) \to
f^{-1} \big( \Vtx_{\neq2}(G') \setminus Z \big)$.
Moreover, for each $z\in Z$ we have $| \sigma^{-1}(z) | = 2$ and 
$| f^{-1}(z) | \in \{ 1, 2 \}$, so we may define a surjection
$\theta_z : \sigma^{-1}(z) \to f^{-1}(z)$. Gluing $\theta_1$ with
the various $\theta_z$ gives a surjection
$\theta$ from  $\Vtx_{\neq2}(S) = \Vtx(S)$ to $\Vtx_{\neq2}(G)$ 
satisfying $f\circ\theta = \sigma$;
it is easily verified that $\theta$ is a skeletal map, 
$\theta : S \skm G$, so assertion~(1) is true.

Let $f: S \skm S'$ be as in assertion~(2).
By the above properties (i) and (ii), it follows that $f$ is a
bijection $\Vtx(S) \to \Vtx(S')$; consequently,
$\vec f : P(S) \to P(S')$ is bijective.
Since $S$ is a skeleton,
$P(S)$ is exactly the set of
ordered pairs $(u,v)$ such that $\{ u,v \}$ is an edge of $S$
(and similarly for $P(S')$).
So the bijectivity of $\vec f$ implies that $f$ is 
an isomorphism of graphs, which proves (2).
\end{proof}

\begin{lemma}\label{ExistsUniqueSkel}
Given a forest $G$, there exist a skeleton $S$ and a skeletal
map $\sigma : S \skm G$.
Moreover, the pair $(S,\sigma)$ is unique up to isomorphism of
graphs, i.e., if $S,S'$ are skeletons and 
$S \overset{\sigma}{\skm} G \overset{\sigma'}{\lskm} S'$ are skeletal
maps then there exists an isomorphism of graphs $\theta: S'\to S$
such that $\sigma' = \sigma\circ\theta$.
\end{lemma}

\begin{proof}
To prove uniqueness, consider skeletal maps
$S \overset{\sigma}{\skm} G \overset{\sigma'}{\lskm} S'$
where $S$ and $S'$ are skeletons.
Part~(1) of \ref{Basic:Skel} gives
$\sigma' = \sigma\circ\theta$ for some $\theta : S' \skm S$;
then $\theta$ is an isomorphism of graphs,
by assertion~(2) of \ref{Basic:Skel}.
We prove the existence of $(S,\sigma)$ by induction on
$| \Vtx_{ \in\{0,2\} }(G) |$;
it suffices to show that if $| \Vtx_{ \in\{0,2\} }(G) | > 0$
then there exists a pair $(G',\sigma)$ where $G'$ is a forest
satisfying $| \Vtx_{ \in\{0,2\} }(G') | < | \Vtx_{ \in\{0,2\} }(G) |$
and $\sigma : G' \skm G$ is a skeletal map.

If $v$ is a vertex of degree two in $G$ with neighbors
$v_1$ and $v_2$, then remove $v$ and the edges
$\{ v_1, v \}$ and $\{ v, v_2 \}$ and add the edge $\{ v_1, v_2\}$;
let $G'$ be the resulting graph and note that
$\Vtx_{\neq2}(G) = \Vtx_{\neq2}(G')$ and that
the identity map of $\Vtx_{\neq2}(G)$ is a skeletal map $G' \skm G$.
Clearly, $| \Vtx_{ \in\{0,2\} }(G') | < | \Vtx_{ \in\{0,2\} }(G) |$.

If $v$ is a vertex of degree zero in $G$, then add a new vertex
$v^*$ and an edge $\{ v, v^* \}$; 
let $G'$ be the resulting graph and note that
$\Vtx_{\neq2}(G') = \Vtx_{\neq2}(G) \cup \{ v^* \}$.
Define a map $\sigma : \Vtx_{\neq2}(G') \to \Vtx_{\neq2}(G)$
by $\sigma(v^*)=v$ and, for all $u\in \Vtx_{\neq2}(G)$, $\sigma(u)=u$.
Then $\sigma : G' \skm G$ is a skeletal map
and, again, $| \Vtx_{ \in\{0,2\} }(G') | < | \Vtx_{ \in\{0,2\} }(G) |$.
\end{proof}


So far, we defined skeletal maps for forests which are not weighted.
We need that notion in the weighted case as well:

\begin{definition}\label{SkOfWeighted}
By a \textit{skeletal map from $\GG$ to $\GG'$},
where $\GG$ and $\GG'$ are forests which may or may not be weighted,
we mean a skeletal map from 
the underlying graph of $\GG$ to the underlying graph of $\GG'$.
The symbol $f : \GG \skm \GG'$ means that 
$f$ is a skeletal map from $\GG$ to $\GG'$.
If $\GG$ is a weighted forest then, by \ref{ExistsUniqueSkel},
there exists a pair $(S,\sigma)$ where $S$ is a skeleton
(so $S$ is not weighted)
and $\sigma : S \skm \GG$ is a skeletal map, and moreover $(S,\sigma)$
is unique up to isomorphism of graphs.
By \textit{the skeleton of a weighted forest\/} $\GG$,
we mean a skeleton $S$ such that there exists a skeletal map $S \skm \GG$.
So the skeleton of $\GG$ is not weighted.
\end{definition}

\section{Strict equivalence of weighted graphs}
\label{Sec:StrictEqu}

\begin{definition}\label{DefStrictEquiv}
Let $\GG$ and $\GG'$ be weighted graphs, where $\GG'$ is a blowing-up
of $\GG$. If $\GG'$ is either a subdivisional blowing-up of $\GG$,
or a blowing-up of $\GG$ at a vertex $v$ such
that $\deg(v,\GG)\le1$,
we say that $\GG'$ is a \textit{strict blowing-up\/} of $\GG$
and that $\GG$ is a \textit{strict blowing-down\/} of $\GG'$.
Two weighted graphs are \textit{strictly equivalent\/} if
one can be transformed into the other by a sequence of
strict blowings-up and strict blowings-down.
\end{definition}

\begin{theorem}\label{EquStrict}
For minimal weighted graphs, equivalence implies strict equivalence.
\end{theorem}

Although this result is stated and proved for general weighted graphs, 
we will only use it on forests.
The result is of fundamental importance for the remainder of this paper
because it allows us to restrict our attention to strict equivalence,
which has the effect (as we will see) of fixing the skeleton.
See also \ref{RegStrict}, \ref{AllReg}.
For the proof of \ref{EquStrict} we need:

\begin{sublemma}\label{SeqNotMin}
Consider a sequence $( \Gn0, \dots, \Gn n)$ of weighted graphs
satisfying:
\begin{enumerate}

\item There exists an integer $m$ satisfying $0<m<n$ and such that
$\Gn i$ is a blowing-up (resp.\ blowing-down) of $\Gn{i-1}$
for all $i$ such that $0<i\le m$ (resp.\ $m<i\le n$);

\item $\Gn1$ is a blowing-up of $\Gn0$ at a vertex $v$;

\item $v\in \Vtx(\Gn n)$ and $\deg(v,\Gn{n-1}) > \deg(v,\Gn{n})$.

\end{enumerate}
Then there exists a sequence of blowings-up and blowings-down
which transforms $\Gn0$ into $\Gn n$ in fewer than $n$ operations.
\end{sublemma}

Regarding the statement of \ref{SeqNotMin}, two remarks are in order.
First, $( \Gn0, \dots, \Gn n)$ is a sequence of blowings-up
and blowings-down which transforms $\Gn0$ into $\Gn n$ in exactly $n$
operations; the lemma claims that there is a shorter sequence
achieving the same thing, but not necessarely satisfying (1--3).
Secondly, the blowings-up come with injective maps
$$
\Vtx(\Gn{i-1}) \hookrightarrow \Vtx(\Gn{i})\ \ ( 0<i\le m )
\qquad \text{and} \qquad
\Vtx(\Gn{i-1}) \hookleftarrow \Vtx(\Gn{i})\ \ ( m<i\le n )
$$
which allow us to identify every set $\Vtx(\Gn j)$ with a subset
of $\Vtx(\Gn m)$; so the statement of condition (3) makes sense.

\begin{proof}[Proof of \ref{SeqNotMin}]
The graph $\Gn n$ is the blowing-down of $\Gn{n-1}$ at some vertex
$u$; since this decreases the degree of $v$,
$\{ u \}$ is in fact a branch of $\Gn{n-1}$ at $v$.
We proceed by induction on $n$. The case $n=2$ is either as in
\eqref{n=2.1} or as in \eqref{n=2.2}:
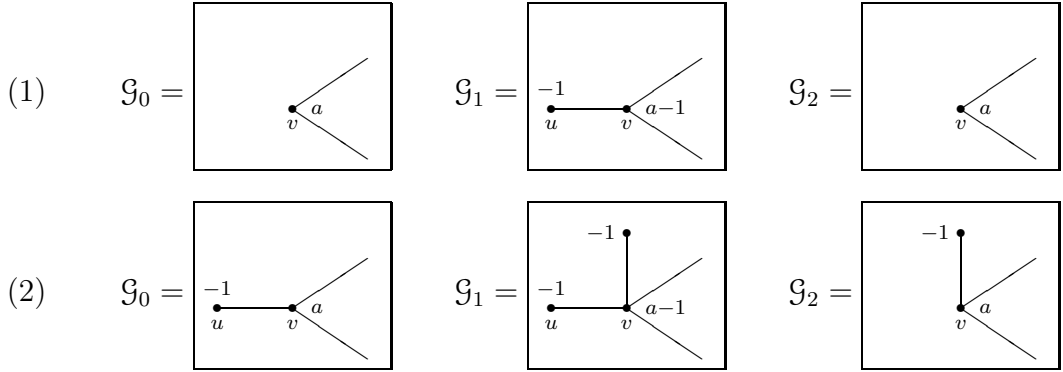
\begin{equation}\label{n=2.1}
\Gn0 =\raisebox{-8\unitlength}{\fbox{%
\begin{picture}(24,20)(-2,-7)
\put(10,0){\circle*{1}}
\put(10,0){\line(3,2){10}}
\put(10,0){\line(3,-2){10}}
\put(10,-1.5){\makebox(0,0)[t]{$\scriptstyle v$}}
\put(12.5,0){\makebox(0,0)[l]{$\scriptstyle a$}}
\end{picture}%
}}
\qquad
\Gn1 =\raisebox{-8\unitlength}{\fbox{%
\begin{picture}(24,20)(-2,-7)
\put(0,0){\circle*{1}}
\put(10,0){\circle*{1}}
\put(0,0){\line(1,0){10}}
\put(10,0){\line(3,2){10}}
\put(10,0){\line(3,-2){10}}
\put(0,1.5){\makebox(0,0)[b]{$\scriptstyle -1$}}
\put(0, -1.5){\makebox(0,0)[t]{$\scriptstyle u$}}
\put(10,-1.5){\makebox(0,0)[t]{$\scriptstyle v$}}
\put(12.5,0){\makebox(0,0)[l]{$\scriptstyle a-1$}}
\end{picture}%
}}
\qquad
\Gn2 =\raisebox{-8\unitlength}{\fbox{%
\begin{picture}(24,20)(-2,-7)
\put(10,0){\circle*{1}}
\put(10,0){\line(3,2){10}}
\put(10,0){\line(3,-2){10}}
\put(10,-1.5){\makebox(0,0)[t]{$\scriptstyle v$}}
\put(12.5,0){\makebox(0,0)[l]{$\scriptstyle a$}}
\end{picture}%
}}
\end{equation}
\begin{equation}\label{n=2.2}
\Gn0 =\raisebox{-8\unitlength}{\fbox{%
\begin{picture}(24,20)(-2,-7)
\put(0,0){\circle*{1}}
\put(10,0){\circle*{1}}
\put(0,0){\line(1,0){10}}
\put(10,0){\line(3,2){10}}
\put(10,0){\line(3,-2){10}}
\put(0,1.5){\makebox(0,0)[b]{$\scriptstyle -1$}}
\put(0, -1.5){\makebox(0,0)[t]{$\scriptstyle u$}}
\put(10,-1.5){\makebox(0,0)[t]{$\scriptstyle v$}}
\put(12.5,0){\makebox(0,0)[l]{$\scriptstyle a$}}
\end{picture}%
}}
\qquad
\Gn1 =\raisebox{-8\unitlength}{\fbox{%
\begin{picture}(24,20)(-2,-7)
\put(0,0){\circle*{1}}
\put(10,0){\circle*{1}}
\put(10,10){\circle*{1}}
\put(0,0){\line(1,0){10}}
\put(10,0){\line(3,2){10}}
\put(10,0){\line(3,-2){10}}
\put(10,0){\line(0,1){10}}
\put(0,1.5){\makebox(0,0)[b]{$\scriptstyle -1$}}
\put(0, -1.5){\makebox(0,0)[t]{$\scriptstyle u$}}
\put(10,-1.5){\makebox(0,0)[t]{$\scriptstyle v$}}
\put(12.5,0){\makebox(0,0)[l]{$\scriptstyle a-1$}}
\put(8.5,10){\makebox(0,0)[r]{$\scriptstyle -1$}}
\end{picture}%
}}
\qquad
\Gn2 =\raisebox{-8\unitlength}{\fbox{%
\begin{picture}(24,20)(-2,-7)
\put(10,0){\circle*{1}}
\put(10,10){\circle*{1}}
\put(10,0){\line(3,2){10}}
\put(10,0){\line(3,-2){10}}
\put(10,0){\line(0,1){10}}
\put(10,-1.5){\makebox(0,0)[t]{$\scriptstyle v$}}
\put(12.5,0){\makebox(0,0)[l]{$\scriptstyle a$}}
\put(8.5,10){\makebox(0,0)[r]{$\scriptstyle -1$}}
\end{picture}%
}}
\end{equation}
In \eqref{n=2.1}, $\Gn2=\Gn0$;
in \eqref{n=2.2}, $\Gn2$ is isomorphic to $\Gn0$;
so in all cases $\Gn0$ can be transformed into $\Gn2$ in zero steps.
Hence, the case $n=2$ is true.

Let $n>2$ and assume that the result is true for smaller values of $n$.
Let $\Beul$ be the branch of $\Gn m$ at $v$ such that $u\in\Vtx(\Beul)$.
Define vertices $e_1,\dots, e_m$ by writing
$\Vtx(\Gn i) = \Vtx(\Gn{i-1}) \cup \{ e_i \}$ for $1\le i \le m$.
Consider the set $E = \{ e_1, \dots, e_m \} \cap \Vtx(\Beul)$.

If $E\neq\emptyset$ then there exists $e_j \in E$ which also satisfies
$w(e_j,\Gn m)=-1$ (namely, let $j = \max\setspec{i}{e_i \in E}$).
Since every
vertex of $\Beul$ disappears in the blowing-down process from $\Gn m$
to $\Gn n$, we may consider $k$ such that $m<k\le n$ and such that
$\Gn k$ is the blowing-down of $\Gn{k-1}$ at $e_j$.
Because $w(e_j,\Gn m)=-1$, it follows that
\begin{equation}\label{Condej}
\textit{$w(e_j, \Gn i)=-1$ for all $i$ such that $j\le i<k$.}
\end{equation}
Since $e_j$ is created by a blowing-up and later deleted by a
blowing-down, and since \eqref{Condej} holds, it follows that
these two steps, the creation and deletion of $e_j$, can be omitted.
So $\Gn0$ can be transformed into $\Gn{n}$ in $n-2$ steps and
we are done.

So we may assume that $E=\emptyset$.
Since $\Vtx(\Gn m) = \Vtx(\Gn0) \cup \{ e_1, \dots, e_m \}$,
this implies that
$ \Vtx(\Beul) \subset \Vtx(\Gn0)$,
i.e.,  that $\Beul$ is present in $\Gn0$ and is ``ready to be shrunk''.
So we may reorder the blowings-up and blowings-down in
$(\Gn0, \dots, \Gn n)$ in such a way that (i) we first perform a
sequence of $p = | \Vtx(\Beul) |$ blowings-down, at the end of which
$\Beul$ has disappeared; (ii) then we perform a sequence of $m$
blowings-up, corresponding exactly to the operations performed in
$(\Gn0, \dots, \Gn m)$; (iii) then we perform $n-m-p$ blowings-down
so as to obtain $\Gn n$ at the end.
In other words, there exists a sequence $(\Gn0', \dots, \Gn n')$ such that
$\Gn 0'=\Gn0$, $\Gn n'=\Gn n$ and the following hold
(where $p = | \Vtx(\Beul) |$):
\begin{itemize}

\item 
$\Gn i'$ is a blowing-down of $\Gn{i-1}'$ at a vertex of $\Beul$
for all $i$ such that $0<i\le p$;

\item $\Gn{p+1}'$ is the blowing-up of $\Gn p'$ at $v$;

\item 
$\Gn i'$ is a blowing-up or a blowing-down of $\Gn{i-1}'$  
for all $i$ such that $p+1 <i\le n$.

\end{itemize}
Since the last vertex of $\Beul$ disappears in the blowing-down which
transforms $\Gn{p-1}'$ into $\Gn p'$,
it follows that $\Gn{p-1}'$ is the blowing-up of $\Gn{p}'$ at $v$;
so $\Gn{p-1}'\isom \Gn{p+1}'$ and
$\Gn0$ can be transformed into $\Gn{n}$ in $n-2$ steps.
So we are done.
\end{proof}

\begin{proof}[Proof of \ref{EquStrict}]
Let $\GG \sim \GG'$ be minimal weighted graphs.
Then there exist sequences of blowings-up and blowings-down which
transform $\GG$ into $\GG'$. Among all such sequences, choose
one
$$
s = ( \Gn0, \dots, \Gn n)  \qquad \text{(where $\Gn0=\GG$ and $\Gn n=\GG'$)}
$$
of minimal length, i.e., it is impossible to transform 
$\GG$ into $\GG'$ in fewer than $n$ steps.
We may assume that $n>0$, otherwise the assertion holds trivially;
since $\GG$ and $\GG'$ are minimal, it follows that $s$ contains
both blowings-up and blowings-down.
As is well-known, one may arrange (without changing the number of
steps) that all blowings-up are performed before the blowings-down,
i.e., for some $m$ such that $0<m<n$ we have:
{\it $\Gn i$ is a blowing-up (resp.\ blowing-down) of $\Gn{i-1}$
for all $i$ such that $0<i\le m$ (resp.\ $m<i\le n$).}
As explained before the proof of \ref{SeqNotMin}, we may regard
the sets $\Vtx(\Gn i)$ as subsets of $\Vtx(\Gn m)$.  Define vertices
$e_1,\dots, e_m$ by writing
$\Vtx(\Gn i) = \Vtx(\Gn{i-1}) \cup \{ e_i \}$ for $1\le i \le m$.

We claim that every blowing-up and blowing-down in $s$ is strict.
Suppose the contrary. Then we may assume that $s$ contains a
blowing-up which is not strict (if not, interchange $\GG$ and $\GG'$
and work with $s = ( \Gn n, \dots, \Gn 0)$ instead). Let
$\Gn j$ be a non strict blowing-up of $\Gn{j-1}$ (where $j\le m$).

If $\Gn j$ is a free blowing-up of $\Gn{j-1}$ then consider the
connected component $\Beul$ of $\Gn m$ containing $e_j$.
Then $\Beul \sim \Empty$ and
$\Vtx(\Beul) \subseteq \{ e_j, e_{j+1}, \dots, e_m \}$.
Clearly, there exists $e_k \in \Vtx(\Beul)$ satisfying $w(e_k,\Gn m)=-1$.
By minimality of $n$, one sees that $e_k$
must still be present in $\Gn n$, 
otherwise the blowing-up which creates $e_k$ and the blowing-down at
$e_k$ are two operations which could be omitted from $s$;
so $\Vtx(\Beul)\cap\Vtx(\Gn n) \neq \emptyset$.
On the other hand, $\Beul \sim \Empty$ and the fact that
$\Beul$ is a connected component of $\Gn m$ imply that
all vertices of $\Beul$ must disappear in the course of the
blowing-down process from $\Gn m$ to $\Gn n$
(for $\Gn n$ is a minimal weighted graph);
so $\Vtx(\Beul)$ is disjoint from $\Vtx(\Gn n)$,
a contradiction.

So $\Gn j$ must be the blowing-up of $\Gn{j-1}$ at some vertex
$v$ satisfying $\deg(v,\Gn{j-1}) \ge2$.
Consider the branch $\Beul$ of $\Gn m$ at $v$ such that $e_j \in
\Vtx(\Beul)$.  Again, we have $\Beul \sim \Empty$ and
$\Vtx(\Beul) \subseteq \{ e_j, e_{j+1}, \dots, e_m \}$;
also, the set $E = \setspec{e \in \Vtx(\Beul)}{ w(e,\Gn m)=-1}$
is not empty and each $e\in E$ must still be present in $\Gn n$,
by minimality of $n$ (see the previous paragraph);
so $\emptyset \neq E \subseteq \Vtx(\Gn n)$.
Since $\Gn n$ is minimal, $w(e,\Gn n) \neq -1$ for all $e\in E$;
in particular, the weight of some vertex of $\Beul$ is increased
by the blowing-down process. It follows that some vertex of
$\{v\}\cup\Beul$ disappears in the blowing-down; consequently,
it makes sense to consider the least integer $k$ such that $m<k\le n$ and:
{\it $\Gn k$ is the blowing-down of $\Gn {k-1}$ at some
vertex $u$ of $\{v\}\cup\Beul$}. From the minimality of $k$ it follows
that $w(x,\Gn m) = w(x,\Gn{k-1})$ for all $x\in\Vtx(\Beul)$;
so, if $u\in\Vtx(\Beul)$, we must have $-1 = w(u,\Gn{k-1}) =  w(u,\Gn{m})$,
so $u\in E$, which is absurd (because $u\not\in\Vtx(\Gn n)$ and 
$E \subseteq \Vtx(\Gn n)$). This shows that 
$\Gn k$ is the blowing-down of $\Gn {k-1}$ at $v$.
In particular $\deg(v,\Gn{k-1}) \le 2 < \deg(v,\Gn m)$ so there exists
an integer $\ell$ such that $m<\ell\le k$ and
$\deg(v,\Gn{\ell-1}) > \deg(v,\Gn{\ell})$.
Then $(\Gn{j-1}, \Gn j, \dots, \Gn{\ell})$ satisfies the hypothesis
of \ref{SeqNotMin} and consequently there exists
a sequence of blowings-up and blowings-down
which transforms $\Gn {j-1}$ into $\Gn{\ell}$ in fewer than $\ell-j+1$
operations.
It follows that $\Gn 0$ can be transformed into $\Gn n$ in fewer than
$n$ operations, which is a contradiction.
\end{proof}

%

Result \ref{EquStrict} has several consequences; we begin with:

\begin{corollary}\label{AllLinChains}
Let $\LL$ and $\LL'$ be equivalent linear chains.
Then there exists a sequence of blowings-up and blowings-down 
which transforms $\LL$ into $\LL'$ and which has the additional
property that every graph which occurs in the sequence is itself
a linear chain.
\end{corollary}

\begin{proof}
There exists a sequence of blowings-down 
which transforms $\LL$ into a minimal weighted graph $\Meul$;
then $\Meul$ and every graph which occurs in this sequence
is a linear chain.
By this remark (also applied to $\LL'$), we may assume that
both $\LL$ and $\LL'$ are minimal. Then \ref{EquStrict} implies that
there exists a sequence of strict blowings-up and strict blowings-down
which transforms $\LL$ into $\LL'$; this sequence has the desired
additional property, because any weighted graph
strictly equivalent to a linear chain is itself a linear chain.
\end{proof}

See \ref{DefSkelMap} and \ref{SkOfWeighted} for the notion
of skeletal map.

\begin{definition}\label{DefBlowDownMap}
Given a weighted forest $\GG$ and a strict blowing-down $\GG'$ of
$\GG$, we shall now define a skeletal map $\pi: \GG \skm \GG'$,
which we call the \textit{blowing-down map}.

Say that $\GG'$ is the blowing-down of $\GG$ at $e$,
so $\Vtx(\GG) = \{ e \} \cup \Vtx(\GG')$.
We define the set map
$\pi : \Vtx_{\neq2}(\GG) \to \Vtx_{\neq2}(\GG')$ as follows.
\begin{itemize}

\item If $\deg(e,\GG)=2$ then $\Vtx_{\neq2}(\GG) = \Vtx_{\neq2}(\GG')$
and we let $\pi$ be the identity map.

\item If $\deg(e,\GG)=1$ then let $v\in\Vtx(\GG)$ be the unique
neighbor of $e$ in $\GG$. Since the blowing-down is strict,
we have $\deg(v,\GG)<3$; so $v\in\Vtx_{<2}(\GG')$, which allows us to
define $\pi(e)=v$.
In order to define $\pi$ on $\Vtx_{\neq2}(\GG) \setminus \{ e \}$,
we note:
$$
\Vtx_{\neq2}(\GG) \setminus \{ e \} =
\begin{cases}
\Vtx_{\neq2}(\GG'), & \text{if } \deg(v,\GG)=1 \\
\Vtx_{\neq2}(\GG') \setminus \{ v\}, & \text{if } \deg(v,\GG) \neq 1
\end{cases}
$$
so it makes sense to define $\pi(x)=x$
for all $x\in \Vtx_{\neq2}(\GG) \setminus \{ e \}$.
\end{itemize}
One can verify that $\pi$ is a skeletal map from $\GG$ to $\GG'$.
\end{definition}

\begin{lemma}\label{StrictEqSameSk}
Strictly equivalent weighted forests have isomorphic skeletons.
\end{lemma}

\begin{proof}
It suffices to verify that if $\GG$ is a weighted forest
and $\GG'$ is a strict blowing-down of $\GG$,
then $\GG$ and $\GG'$ have the same skeleton.  Consider the skeleton
$S$ of $\GG$; then (by definition) there exists a skeletal map
$\sigma : S \skm \GG$. Composing $\sigma$ with the blowing-down map
$\pi: \GG \skm \GG'$ gives a skeletal map $S \skm \GG'$, so
$S$ is the skeleton of $\GG'$.
\end{proof}

Together with \ref{EquStrict}, this gives:

\begin{corollary}\label{SameSkel}
Let $\Ceul$ be the equivalence class of some weighted forest.
Then any two minimal elements of $\Ceul$ have isomorphic skeletons.
\end{corollary}

The special case of \ref{SameSkel}
where $\Ceul$ is the equivalence class of a linear
chain is the following well-known fact:

\begin{center}
\it 
Any minimal weighted graph equivalent to a linear chain is a linear chain.
\end{center}

%

\section*{Pseudo-minimal forests}

\begin{definition}\label{DefWgamma}
Given a weighted forest $\GG$,  we define $P(\GG) = P(G)$ where $G$ is
the underlying graph of $\GG$ (see \ref{DefPT}).
Given $\gamma = (v_0,\dots,v_n) \in P(\GG)$,
let $i_1 < \cdots < i_p$ be the elements of the set
$\setspec{i}{\deg(v_i,\GG)<3}$. Then define
$$
W_\GG (\gamma) = \big( w(v_{i_1},\GG), \dots, w(v_{i_p},\GG) \big) \in \Integ^*
$$
(see \ref{NotaSequences} for the notation $\Integ^*$).
This defines a set map $W_\GG  : P(\GG) \to \Integ^*$.
\end{definition}

\begin{lemmadef}\label{LemDefRegForest}
For a weighted forest $\GG$, the following conditions are equivalent:
\begin{enumerate}

\item Every $\gamma\in P(\GG)$ satisfying $W_\GG(\gamma)\sim\emptyseq$
is of type $(+,+)$;

\item $\GG$ is strictly equivalent to a minimal weighted forest.

\end{enumerate}
{\rm We call $\GG$ a {\it pseudo-minimal forest} if conditions
(1), (2) hold.}
\end{lemmadef}

\begin{proof}
Let us say (only in this proof) that a weighted forest $\GG$ is
1-regular (resp.\ 2-regular) if it satisfies condition (1) 
(resp.\ condition~(2)). We first show that if weighted forests
$\GG$ and $\GG'$ are strictly equivalent, and if one of them
is 1-regular, then both are 1-regular.  We may assume that
$\GG'$ is a strict blowing-down of $\GG$;
let $\pi : \GG \skm \GG'$ be the blowing-down map
defined in \ref{DefBlowDownMap}. Recall
that $\vec{\pi} : P(\GG) \to P(\GG')$ is surjective and
(see \ref{PreservesType}) preserves type.  Also, it is
clear from the definition of $\pi$ that, for every $\gamma\in P(\GG)$,
the sequences $W_\GG(\gamma)$ and $W_{\GG'} \big( \vec{\pi}( \gamma ) \big)$
are equivalent. It follows that $\GG$ is 1-regular if and only if $\GG'$ is
1-regular.

Since every minimal weighted forest is 1-regular, the above paragraph
implies that every 2-regular forest is 1-regular.

Conversely, consider a 1-regular forest $\GG$. Note that 1-regularity
implies that if $\GG'$ is any blowing-down of $\GG$, then in fact
$\GG'$ is a strict blowing-down of $\GG$. Then, by the first
paragraph, $\GG'$ is 1-regular. Reiterating this argument shows
that $\GG$ is strictly equivalent to a minimal forest, i.e., $\GG$ is
2-regular.
\end{proof}

\begin{corollary}\label{RegStrict}
For pseudo-minimal forests, equivalence implies strict equivalence.
Consequently, equivalent pseudo-minimal forests have isomorphic skeletons.
\end{corollary}

\begin{proof}
If $\GG$ and $\GG'$ are pseudo-minimal forests then they are
strictly equivalent to minimal forests $\Meul$ and $\Meul'$ respectively.
If $\GG \sim \GG'$, then $\Meul \sim \Meul'$;
so $\Meul$ is strictly equivalent to $\Meul'$ by \ref{EquStrict}, 
and consequently $\GG$ is strictly equivalent to $\GG'$.
The last assertion follows from \ref{StrictEqSameSk}.
\end{proof}

The next statement is simply a reformulation of \ref{RegStrict}:

\begin{corollary}\label{AllReg}
Let $\GG$ and $\GG'$ be equivalent weighted forests.
If $\GG$ and $\GG'$ are minimal,
or more generally if they are pseudo-minimal,
then there exists a sequence of operations which transforms one into
the other and such that:
\begin{enumerate}
\item Every graph which occurs in the sequence is a pseudo-minimal forest
\item every operation in the sequence is a strict blowing-up or a
strict blowing-down.
\end{enumerate}
\end{corollary}

\section{Finite sequences of integers}
\label{SectionSequences}

We consider $\Integ^*$ and $\Neul^*$, 
where $\Neul = \setspec{ x\in\Integ }{ x<-1 }$
(see \ref{NotaSequences} for the notation $E^*$, where $E$ is a set).
As indicated in \ref{DefLinChain}, given $A\in\Integ^*$
we write $[A]$ for the corresponding linear chain.

This section classifies elements of $\Integ^*$ up to equivalence. 
From this, a classification of linear chains will be obtained,
in the next section,
by recalling that $[A] = [A^-]$ for all $A \in \Integ^*$.

The material up to \ref{SameSUB}
is well known when stated for linear chains.
The main results of the section are
\ref{SeqThm} (complemented by \ref{uniqueness}),
\ref{EquivSeq} and
\ref{BigPicture} (complemented by \ref{ListPrimes}).

\begin{notation}\label{NotaSeq}
For each $i\in\{1,\dots,r\}$, let $A_i$ be either an integer or an
element of $\Integ^*$.
We write $(A_1,\dots,A_r)$ for the concatenation of $A_1, \dots, A_r$;
that is, $(A_1,\dots,A_r)\in\Integ^*$ is a single sequence.
Also, we will use superscripts to indicate repetitions.
For instance, if $A= (0^3,-5,-1) \in \Integ^*$
and $B= (-2^3,3,-2) \in \Integ^*$ then 
$$
(A,-2,B) = ( 0^3,-5,-1,-2, -2^3,3,-2)
= (0,0,0,-5,-1,-2,-2,-2,-2,3,-2).
$$
Superscripts occurring in sequences (or linear chains) should always
be interpreted in this way, never as exponents.
\end{notation}

\begin{definition}\label{DefBupSeq}
If $X = (x_1,\dots,x_n) \in \Integ^*$ and $X\neq\emptyseq$,
then any of the following sequences $X'\in\Integ^*$ is called
a \textit{blowing-up\/} of $X$:
\begin{itemize}

\item $X' = (-1,x_1-1,x_2,\dots,x_n)$;
\item $X' = (x_1,\dots,x_{i-1}, x_i-1, -1, x_{i+1}-1, x_{i+2}, \dots,x_n)$
(where $1\le i < n$);
\item $X' = (x_1,\dots,x_{n-1},x_n-1,-1)$.
\end{itemize}
Moreover, we regard the one-term sequence $(-1)$ as a blowing-up
of the empty sequence $\emptyseq$.
If $X'$ is a blowing-up of $X$, we also say that $X$ is a blowing-down
of $X'$.
Two elements of $\Integ^*$ are said to be \textit{equivalent\/} if
one can be obtained from the other by a finite sequence of blowings-up
and blowings-down. This defines an equivalence relation ``$\sim$''
on the set $\Integ^*$.
We also consider the partial order relation ``$\le$'' on the set
$\Integ^*$ which is generated by the condition:
\begin{center}
$X \le X'$ whenever $X'$ is a blowing-up of $X$.
\end{center}
Thus a minimal element of $\Integ^*$ is a sequence which cannot
be blown-down, i.e., an element of $( \Integ \setminus \{ -1 \})^*$.
\end{definition}

\begin{lemma}\label{BasicBup}
Given $X,Y\in\Integ^*$,
\begin{enumerate}

\item $X\sim Y\ \iff\  X^-\sim Y^-$

\item $[X]\sim[Y] \iff X\sim Y \text{\ or\ } X\sim Y^-$. 

\end{enumerate}
\end{lemma}

\begin{proof}
The only nontrivial claim is implication ``$\Rightarrow$'' of assertion~(2),
and this easily follows from \ref{AllLinChains}.
\end{proof}

Refer to \ref{DefBilForm} and \ref{DefDet} for the following:

\begin{definition}
Given $X\in\Integ^*$, we define
$\det(X) = \det\big( [X] \big)$ and $\hodge{X} = \hodge{[X]}$.
\end{definition}

\begin{lemma}\label{TwoInvars}
If $X,Y\in\Integ^*$ satisfy $X\sim Y$, then
$\det(X) = \det(Y)$ and $\hodge{X} = \hodge{Y}$.
\end{lemma}

\begin{proof}
Follows from \ref{BasicBup}, \ref{HodgeInvar} and \ref{DetInvar}.
\end{proof}

Fact \ref{TwoInvars} allows us to define
$\det(\Ceul)$ and $\hodge{\Ceul}$ for any
equivalence class $\Ceul \subset \Integ^*$
(the definitions are the obvious ones).

\begin{notation}\label{deti}
Given $X=( x_1, \dots, x_n ) \in \Integ^*$, define:
\begin{align*}
{\textstyle\det_i}(X) &=
\begin{cases}
\det \big( x_{i+1},\dots,x_n  \big), &\text{ if } 0\le i<n,\\
1, &\text{ if }i=n,\\
0, &\text{ if }i>n;
\end{cases}\\
{\textstyle\det_*(X)} &= 
\begin{cases}
\det(x_2,\dots,x_{n-1}), & \text{if $n>2$,} \\
1,  & \text{if $n=2$,} \\
0,  & \text{if $n<2$}.
\end{cases}
\end{align*}

In particular, note that $\det_0(X)=\det(X)$.
The sequence $X$ determines the ordered pair
$$
\textstyle
\Sub(X) = \big( \det_1(X), \det_1(X^-) \big)
$$
which is an element of the $\Integ$-module $\Integ\times\Integ$.
This gives in particular $\Sub(\emptyseq)= (0,0)$ and
if $a\in\Integ$, $\Sub\big( (a) \big) = (1,1)$.
Finally, let $d = \det(X)$ and define the pair
$$
\textstyle
\SUB(X) = \big( \pi(\det_1(X)), \pi(\det_1(X^-)) \big) 
\in \Integ/ d\Integ \times \Integ/ d\Integ,
$$
where $\pi : \Integ \to \Integ/d\Integ$
is the canonical epimorphism and where we regard 
$\Integ/ d\Integ \times \Integ/ d\Integ$ as a $\Integ$-module.
\end{notation}

Facts \ref{remark*}--\ref{SameSUB} are, in one form or another,
contained in \cite{Hirzebruch53}. We omit some proofs.

\begin{lemma}\label{remark*}
If $X=( x_1, \dots, x_n ) \in \Integ^*$ then:
$$
\textstyle\det_i(X)= (-x_{i+1}) \det_{i+1}(X)-\det_{i+2}(X)
\qquad (0\le i<n).
$$
In particular, $\det X = (-x_1)\det_1(X) - \det_2(X)$.
\end{lemma}

\rien{
Proof is omitted

\begin{proof}
Let $v_1,\dots,v_n$ be the vertices of $\GG=[X]=[x_1,\dots,x_n]$,
where the labelling is such that $w(v_i,\GG)=x_i$ and
$\{v_i,v_{i+1}\}$ is an edge for every $i$.
Let $M$ be as in \ref{DefDet}.
The Laplace expansion of $\det(-M)$ along the first row gives
$\det_0(X) = (-x_1)\det_1(X) - \det_2(X)$.
Applying this formula to $(x_{i+1},\dots,x_n)$ gives the desired
result.
\end{proof} } 

\begin{lemma}\label{DetAdmiss}\label{EquivNeul}
The assignment $X\mapsto ( \det(X), \det_1(X) )$ is a well-defined bijection:
$$
\Neul^* \longrightarrow
\setspec{ (r_0,r_1) \in \Nat^2 }
{ \text{$0 \le r_1 < r_0$ and $\gcd(r_0,r_1)=1$} }.
$$
\end{lemma}


\begin{lemma}\label{Claim1Disc}
If $X\in\Integ^*$, $d=\det(X)$ and $(x,y)=\Sub(X)$, then
$xy\equiv 1\pmod{d}$.
\end{lemma}

\rien{
Proof is omitted

\begin{proof}
The result holds trivially when $X=\emptyseq$.
For $X \neq \emptyseq$, we prove 
\begin{equation}\label{xydd}
\textstyle xy=1+d \det_*(X).
\end{equation}
Write $X =(x_1,\dots,x_n)$.
We leave the cases $n=1,2$ to the reader.
Assume that $n>2$ and
that \eqref{xydd} holds for the shorter sequence $(x_2,\dots,x_n)$;
i.e., we are assuming that
\begin{equation}\label{IndHyp}
\textstyle
\det_2(X) \det_*(X)  = 1 + x \delta,
\end{equation}
where $\delta = \det_* (x_2,\dots,x_{n})$.
We obtain $d=-x_1x-\det_2(X)$ by \ref{remark*}, so
\begin{equation}\label{FirstEq}
\textstyle
\det_2(X)  = -x_1 x - d.
\end{equation}
Applying \ref{remark*} to $(x_1,\dots,x_{n-1})$ gives
$y = -x_1 \det_*(X) - \delta$ and hence
\begin{equation}\label{SecondEq}
\textstyle
\delta  = -x_1 \det_*(X) - y.
\end{equation}
Substituting \eqref{FirstEq} and \eqref{SecondEq} in \eqref{IndHyp}
yields the desired conclusion \eqref{xydd}.
\end{proof} } 

\begin{lemma}\label{SubSubxy}
Suppose that $A,B\in \Integ^*$ satisfy $A\sim B$ and let
$d = \det(A)=\det(B)$. Then there exists $(x,y)\in \Integ^2$ such that
\begin{equation}
\label{SubSubxyEQ}
\Sub(A) = \Sub(B) + d(x,y).
\end{equation}
\end{lemma}

\begin{proof}
Note that $\det(A)=\det(B)$ by \ref{TwoInvars}.
Since $A\sim B$,
performing a certain sequence of blowings-up and
blowings-down on $A$ produces $B$; if the same sequence of operations
is performed on $(0,A)$ then (obviously)
we obtain $(x,B)$ for some $x\in\Integ$, which shows that 
$(0,A)\sim(x,B)$.
By the same argument, $(A,0) \sim (B,y)$ for some $y\in\Integ$.
By \ref{remark*} we have
$ \det(0,A) = -\det_1(A)$ and $\det(x,B) = -x d -\det_1(B)$;
since $(0,A)\sim(x,B)$ implies
$\det(0,A) = \det(x,B)$, we obtain 
$$
\textstyle \det_1(A) = \det_1(B) + dx.
$$
Similarly, we have 
$(0,A^-) = (A,0)^- \sim (B,y)^- = (y,B^-)$, so
$\det(0,A^-) = \det(y,B^-)$ and consequently
$\det_1(A^-) = \det_1(B^-) + dy$.
So $(x,y)$ satisfies \eqref{SubSubxyEQ}.
\end{proof}

\begin{corollary} \label{SameSUB}
If $A,B\in\Integ^*$ and $A \sim B$,  then $\SUB(A) = \SUB(B)$.
\end{corollary}

\begin{proof}
Obvious consequence of \ref{SubSubxy}.
\end{proof}

%
%
%
%
%

%
%
%
%

\section*{Classification of sequences up to equivalence}

Sequences of the form $(0^{2i}, A)$ (see \ref{NotaSeq} for notations)
play an important role in the classification.
We need the following facts.

\begin{lemma}\label{0iA}
Let $i\in\Nat$ and $A\in\Integ^*$.
\begin{enumerate}

\item $\det( 0^{2i},A ) = (-1)^i \det( A )$

\item $\Sub( 0^{2i},A ) = (-1)^i \Sub( A )$.

\end{enumerate}
\end{lemma}

\begin{proof}
We may assume that $i>0$, then \ref{remark*} gives
$$
\textstyle
\det( 0^{2i},A )
= 0 \det_1( 0^{2i},A ) - \det_2( 0^{2i},A ) 
= - \det( 0^{2i-2},A )
$$
and assertion (1) follows by induction.
We also have:
\begin{multline}\label{1stHalf}
\textstyle
\det_1(0^{2i},A)
= \det(0^{2i-2},0,A) 
\overset{\mbox{\tiny (1)}}{=}
(-1)^{i-1} \det(0,A)
\\
\textstyle
\overset{\mbox{\tiny\ref{remark*}}}{=}
(-1)^{i-1} (0\det(A) - \det_1(A))
= (-1)^{i} \det_1(A)
\end{multline}
so, to prove (2), there remains only to show that
\begin{equation}\label{2ndHalf}
\textstyle
\det_1 \big( (0^{2i}, A)^- \big) = (-1)^i \det_1( A^- ).
\end{equation}
If $A=\emptyseq$ then \eqref{2ndHalf} reads $\det( 0^{2i-1} ) = 0$,
which is true by assertion~(1). 
So we may assume that $A= (a_1,\dots,a_n)$ with $n\ge1$,
in which case
\begin{multline*}
\textstyle
\det_1 \big( (0^{2i}, A)^- \big)
= \det ( a_{n-1}, \dots, a_1, 0^{2i} )
= \det (0^{2i}, a_1, \dots, a_{n-1} ) \\
\textstyle
\overset{\mbox{\tiny (1)}}{=}   
(-1)^i \det (a_1, \dots, a_{n-1} )
 = (-1)^i \det (a_{n-1}, \dots, a_1 )
 = (-1)^i \det_1( A^- ).
\end{multline*}
So \eqref{2ndHalf} holds and assertion~(2) follows from 
\eqref{1stHalf} and \eqref{2ndHalf}.
\end{proof}

\begin{lemma}\label{BasicHodge}
If $i\in\Nat$ and $A\in\Integ^*$,
then $\hodge{(0^{2i},A)}= i + \hodge{A}$.
\end{lemma}

\begin{proof}
This is an exercise in diagonalization.
It suffices to prove that $\hodge{(0,0,A)}= 1 + \hodge{A}$
for every $A \in \Integ^*$.
This is obvious if $A=\emptyseq$, so assume that $A\neq\emptyseq$
and write $A = (a_1, \dots, a_n)$. Consider the linear chain
$$
\LL=[0,0,A] =  \raisebox{-2\unitlength}{%
\begin{picture}(44,6)(-2,-3)
\put(0,0){\circle*{1}}
\put(10,0){\circle*{1}}
\put(20,0){\circle*{1}}
\put(30,0){\makebox(0,0){\mbox{\,\dots}}}
\put(40,0){\circle*{1}}
\put(0,0){\line(1,0){26}}
\put(40,0){\line(-1,0){6}}
\put(0,1.5){\makebox(0,0)[b]{$\scriptstyle 0$}}
\put(10,1.5){\makebox(0,0)[b]{$\scriptstyle 0$}}
\put(20,1.5){\makebox(0,0)[b]{$\scriptstyle a_1$}}
\put(40,1.5){\makebox(0,0)[b]{$\scriptstyle a_n$}}
\put(0, -1.5){\makebox(0,0)[t]{$\scriptstyle u_1$}}
\put(10,-1.5){\makebox(0,0)[t]{$\scriptstyle u_2$}}
\put(20,-1.5){\makebox(0,0)[t]{$\scriptstyle v_1$}}
\put(40,-1.5){\makebox(0,0)[t]{$\scriptstyle v_n$}}
\end{picture}}
$$
and let $V$ be the real vector space with basis $\Vtx(\LL)$.
Then the matrix representing $B_\LL$ with respect to the basis
$(u_1, u_2, v_1-u_1, v_2, \dots, v_n)$ of $V$ is:
\begin{equation}\label{DiagBlock}
\left(
\begin{array}{cc|ccc}
0 & 1 & 0 & \cdots & 0 \\
1 & 0 & 0 & \cdots & 0 \\
\hline
0      & 0       & &   & \\
\vdots & \vdots  & & M & \\
0      & 0       & &   & 
\end{array}
\right)
\end{equation}
where $M$ is the $n\times n$ matrix given by $M_{ii}=a_i$,
$M_{ij}=1$ if $|i-j|=1$ and $M_{ij}=0$ if $|i-j|>1$,
that is, $M$ is the matrix
representing the intersection form of the linear chain $[A]$.
Now $\left(\begin{smallmatrix} 0&1\\1&0 \end{smallmatrix}\right)$
can be diagonalized to 
$\left(\begin{smallmatrix} 1&0\\0&-1 \end{smallmatrix}\right)$
and we conclude that a diagonal matrix congruent to \eqref{DiagBlock}
has $1 + \hodge{A}$ nonnegative entries on its main diagonal,
i.e., $\hodge{\LL}=1 + \hodge A$.
\end{proof}

\begin{lemma}\label{BasicZero}
Let $a,b,x\in\Integ$ and $A,B\in\Integ^*$. Then
$$
(A,a,0,b,B) \sim (A,a-x,0,b+x,B) \sim (A,0,0,a+b,B).
$$
\end{lemma}

\begin{proof}
$(A,a,0,b,B) \sim (A,a-1,-1, -1 ,b,B) \sim (A,a-1,0,b+1,B)$,
from which the result follows.
\end{proof}

\begin{lemma}\label{MoveZeros}
Let $n\in\Nat$ and $A,B,C\in\Integ^*$. Then
$ (A,B,0^{2n},C) \sim (A,0^{2n},B,C) $.
\end{lemma}

\begin{proof}
If $A,C\in\Integ^*$ and $b\in\Integ$ then by \ref{BasicZero}
$$
(A,b,0,0,C) \sim (A,b-b,0,0+b,C) = (A,0,0,b,C),
$$
from which the result follows.
\end{proof}

%

\begin{lemma}\label{ChangeOddValue}
Let $n\in\Nat$, $x,y\in\Integ$ and $A\in\Integ^*$. Then
$ (0^{2n+1},x,A) \sim (0^{2n+1},y,A)$.
\end{lemma}

\begin{proof}
We first consider the case $n=0$:
$(0,x,A) \sim (-1,-1,x-1,A) \sim (0,x-1,A)$, from which we deduce
$(0,x,A) \sim (0,y,A)$.
Now the general case:
$$
(0^{2n+1},x,A)
\simref{\ref{MoveZeros}} (0,x,0^{2n},A)
\simref{$(n=0)$} (0,y,0^{2n},A)
\simref{\ref{MoveZeros}} (0^{2n+1},y,A).
$$
\end{proof}

\begin{lemma}\label{SubSeq}
Let $n\in\Nat$ and $A,B\in\Integ^*$. Then:
$$
A \sim B\ \ \implies\ \ (0^{2n},A) \sim (0^{2n},B).
$$
\end{lemma}

\begin{proof}
We may assume that $n\ge1$.
If $A\sim B$ then performing a certain sequence of blowings-up and
blowings-down on $A$ produces $B$; if the same sequence of operations
is performed on $(0^{2n},A) = (0^{2n-1},0,A)$,
then we obtain
$(0^{2n-1},x,B)$ for some $x\in\Integ$, i.e.,
only the rightmost zero in $0^{2n}$ is affected.
So 
$$
(0^{2n},A) \sim (0^{2n-1},x,B)
\simref{\ref{ChangeOddValue}} (0^{2n-1},0,B) = (0^{2n},B).
$$
\end{proof}

\begin{definition}
Let $B = (b_1, \dots, b_n ) \in \Integ^*$.
\begin{enumerate}

\item Given $x\in\Integ$, define
$_xB = (0,x,B)=(0,x,b_1,\dots,b_n)\in\Integ^*$ and
$B_x = (B,x,0)=(b_1,\dots,b_n, x, 0) \in \Integ^*$.

\item Suppose that $B \neq \emptyseq$.
Given $i \in \{ 1,\dots, n \}$ and $x,y\in\Integ$ such that $x+y=b_i$,
define
$B_{(i;x,y)} = (b_1, \dots, b_{i-1}, x, 0, y, b_{i+1}, \dots, b_n )
\in \Integ^*$.

\end{enumerate}
\end{definition}

\begin{definition}\label{M+}
Given a minimal element $M=(m_1,\dots,m_k)$ of $\Integ^*$,
let $M^{\oplus} $ be the set of sequences $Z\in\Integ^*$ which can be
constructed in one of the following ways.
\begin{enumerate}

\item Pick $x\in\Integ$ and let $Z$ be the unique minimal sequence
such that $Z\le {}_xM$.

\item Pick $x\in\Integ$ and let $Z$ be the unique minimal sequence
such that $Z\le M_x$.

\item Assuming that $M\neq\emptyseq$,
pick $j\in\{ 1, \dots, k\}$ and $x,y\in\Integ$ such that $x+y=m_j$
and let $Z$ be the unique minimal sequence such that $Z \le M_{(j;x,y)}$.

\item Pick $M' = (\mu_1, \dots, \mu_\ell)$ such that $M'\ge M$
and exactly one term $\mu_j$ is equal to $-1$;
pick $x,y\in\Integ\setminus\{ -1 \}$ such that $x+y=-1$ and let
$Z = M'_{(j;x,y)}$.
\end{enumerate}
Note that each element of $M^{\oplus} $ is a minimal element of $\Integ^*$.
\end{definition}

\begin{lemma}\label{LemZM}
If $M$ is a minimal element of $\Integ^*$ and $Z \in M^{\oplus} $ then
$Z \sim (0,0,M)$.
Moreover, $\det Z = -\det M$ and $\hodge{Z} = \hodge{M}+1$.
\end{lemma}

\begin{proof}
By definition \ref{M+} of $M^{\oplus} $, one of the following holds:
$$
Z \le {}_xM,
\quad Z \le M_x,
\quad Z \le  M_{(j;x,y)}
\quad \text{or} \quad
Z =  M'_{(j;x,y)}
\text{ where }
M' \sim M.
$$
Consequently, one of the following holds:
$$
Z \sim {}_xM,
\quad Z \sim M_x
\quad \text{or} \quad
Z \sim  M'_{(j;x,y)}
\text{ where }
M' \sim M.
$$
By \ref{ChangeOddValue}, ${}_xM=(0,x,M) \sim (0,0,M)$.
Since $X\sim Y$ implies $X^- \sim Y^-$, we also have
$M_x= \big({}_x(M^-)\big)^- \sim (0,0,M^-)^- = (M,0,0) \sim (0,0,M)$
by \ref{MoveZeros}.

Let $M'=(b_1,\dots,b_m)$ be any nonempty sequence equivalent to $M$
and let $j \in \{ 1, \dots, m \}$ and $x,y\in\Integ$ be such that $x+y=b_j$;
then
\begin{multline*}
M'_{(j;x,y)} = (b_1, \dots, b_{j-1}, x, 0, y, b_{j+1}, \dots, b_m )
\simref{\ref{BasicZero}}
(b_1, \dots, b_{j-1},0,0, x+y, b_{j+1}, \dots, b_m ) \\
= (b_1, \dots, b_{j-1},0,0, b_j, b_{j+1}, \dots, b_m )
\simref{\ref{MoveZeros}}
(0,0,M')
\simref{\ref{SubSeq}}
(0,0,M).
\end{multline*}
Thus $Z\sim(0,0,M)$ whenever $Z\in M^{\oplus} $.
By \ref{0iA} and \ref{BasicHodge} we get
$\det Z = -\det M$ and $\hodge Z = \hodge M + 1$.
\end{proof}

\begin{proposition}\label{ZM}
Let $Z$ be a minimal element of $\Integ^*$ such that $\hodge Z >0$
and $Z \neq (0)$. Then $Z \in M^{\oplus} $ for some minimal element $M$ of
$\Integ^*$.
\end{proposition}

\begin{proof}
Assume that $Z = (z_1,\dots,z_n)$ is minimal, $\hodge Z>0$ and $Z \neq (0)$.
In particular, $\hodge Z>0$ implies that $z_i\ge-1$ for some $i$;
so by minimality of $Z$ there exists $i$ such that $z_i\ge0$.
If $z_i=0$ for some $i$, we distinguish three cases:
\begin{enumerate}

\item[(i)]
If $z_1=0$ then, since $Z \neq (0)$, we have $Z = (0,x,M)= {}_xM$
for some $M \in \Integ^*$ and $x\in\Integ$;
then $M$ is minimal and $Z\in M^{\oplus} $.

\item[(ii)]
If $z_n=0$ then, similarly, $Z=M_x$ for some $M\in\Integ^*$
and $x\in\Integ$; 
then $M$ is minimal and $Z\in M^{\oplus} $.

\item[(iii)]
If $z_i=0$ for some $i$ such that $1<i<n$,
then $Z = (z_1,\dots, z_{i-1}, 0, z_{i+1}, \dots, z_n)
= B_{(i-1; z_{i-1}, z_{i+1} )}$ where
$B = (z_1,\dots, z_{i-2}, z_{i-1} + z_{i+1},  z_{i+2},\dots, z_n)$.
If $B$ is minimal then $Z\in M^{\oplus} $ where $M=B$.
If $B$ is not minimal then its $(i-1)$-th term ($z_{i-1} + z_{i+1}$)
is the only one which is equal to $-1$;
we have $B\ge M$ for some minimal $M$, then $Z\in M^{\oplus} $.

\end{enumerate}
From now-on, assume that $z_j\neq0$ for all $j\in \{1, \dots, n \}$.
Then $z_i>0$ for some $i$ and we have four cases:
\begin{enumerate}

\item[(iv)]
If $Z = (p)$ where $p>0$, then 
$Z \le (0,-1, -2^{p-1}) = {}_{-1}M$ where $M = ( -2^{p-1})$ is minimal;
then $Z \in M^{\oplus} $.

\item[(v)]
If $z_1>0$ and $n>1$ then 
$Z \le (0,-1, -2^{z_1-1}, z_2-1, z_3, \dots, z_n) = {}_{-1}M$
where $M = (-2^{z_1-1}, z_2-1, z_3, \dots, z_n)$ is minimal; 
then $Z \in M^{\oplus} $.

\item[(vi)]
If $z_n>0$ and $n>1$ then 
$Z \le (z_1, \dots, z_{n-2}, z_{n-1}-1, -2^{z_n-1}, -1,  0) = M_{-1}$
where $M = (z_1, \dots, z_{n-2}, z_{n-1}-1, -2^{z_n-1})$ is minimal;
then $Z \in M^{\oplus} $.

\item[(vii)]
If $z_i>0$ and $1<i<n$ then $Z \le
(z_1, \dots, z_{i-1}, 0, -1,  -2^{z_i-1}, z_{i+1}-1, z_{i+2}, \dots, z_n  )
= M_{(i-1; z_{i-1}, -1)}$, where $M = 
(z_1, \dots, z_{i-2}, z_{i-1}-1, -2^{z_i-1}, z_{i+1}-1,
z_{i+2}, \dots, z_n  )$ is minimal;
then $Z \in M^{\oplus} $.

\end{enumerate}
\end{proof}

\bigskip

\begin{definition}\label{DefCanSeq}
An element $C$ of $\Integ^*$ is
a \textit{canonical sequence} if it has the form
$$
\textit{$C = (0^r, A)$, where $r\in\Nat$,  $A \in\Neul^*$ and
if $A\neq\emptyseq$ then $r$ is even.}
$$
\end{definition}

Then one has the fundamental result:

\smallskip
\begin{theorem}\label{SeqThm}
Each element of $\Integ^*$ is equivalent to a unique canonical sequence.
\end{theorem}

\smallskip
The proof of \ref{SeqThm} consists of
\ref{EveryCan} and \ref{uniqueness}, below.

\begin{subnothing}\label{EveryCan}
Every element of $\Integ^*$ is equivalent to a canonical sequence.
\end{subnothing}

\begin{proof}
It suffices to show that every minimal element $Z$ of $\Integ^*$ is
equivalent to a canonical sequence.  We proceed by induction on
$\hodge Z$.  If $\hodge Z =0$ then $Z \in \Neul^*$, so $Z$ itself is
canonical.
If $\hodge Z>0$ then, by \ref{ZM}, either $Z=(0)$
or $Z\in M^{\oplus} $ for some minimal element $M$ of $\Integ^*$.
In the first case, $Z$ is canonical and we are done.
In the second case, \ref{LemZM} gives $\hodge M < \hodge Z$ so
we may assume by induction that $M$ is equivalent to a canonical sequence $C$;
then $ Z \sim (0,0,M) \sim (0,0,C) $
by \ref{LemZM} and \ref{SubSeq}, and clearly $(0,0,C)$ is canonical.
\end{proof}

\medskip
\begin{subnothing}\label{uniqueness}
Let $L\in\Integ^*$, let $n=\hodge{L}$ and
let $d$ be the absolute value of $\det(L)$.

If $(0^r, A)$ (where $r\in\Nat$ and $A \in \Neul^*$)
is a canonical sequence equivalent to $L$, then:
\begin{enumerate}

\item[(a)] If $d=0$ then 
$r=2n-1$ and $A=\emptyseq$.

\item[(b)] If $d\neq0$ then $r=2n$ and
$A$ is the unique element of $\Neul^*$ which satisfies:
$$
\text{$\det(A)=d$\ \ and\ \ $\SUB(A) = (-1)^n \SUB(L)$.}
$$

\end{enumerate}
In particular, $r$ and $A$ are uniquely determined by $L$.
\end{subnothing}

\begin{proof}
The claim that $r$ and $A$ are uniquely determined by $L$
is obvious in case (a), and follows from  \ref{EquivNeul}
in case (b).
Consider any canonical sequence $(0^r, A)$ equivalent to $L$;
we have $r\in \Nat$, $A\in\Neul^*$, and if $A\neq\emptyseq$ then $r$ is even.
To prove (a) and (b), it suffices to show:
\begin{enumerate}

\item[(a$'$)] \textit{If $r$ is odd then $d=0$ and $r=2n-1$.}

\item[(b$'$)] \textit{If $r$ is even then $\det(A)=d\neq0$, $r=2n$ and
$\SUB(A) = (-1)^n \SUB(L)$.}

\end{enumerate}

If $r$ is odd then $A=\emptyseq$; writing $r=2i+1$, we get
$
\pm d = \det(L)
= \det( 0^{2i+1} )
= \det( 0^{2i},0 )
= (-1)^i \det( 0 ) = 0
$
by \ref{0iA} and
$n=\hodge{L} = \hodge{ (0^{2i+1}) } = i+\hodge{(0)}=i+1$
by \ref{BasicHodge}.  This proves (a$'$).

If $r$ is even then \ref{BasicHodge} gives
$n=\hodge{L} = \hodge{ (0^r,A) } = \frac r2+ \hodge A=\frac r2$,
so $r=2n$. Then \ref{0iA} gives
$\pm d = \det(L) = \det( 0^{2n},A ) = (-1)^{n} \det(A)$;
since $\det(A)>0$ by \ref{DetAdmiss},
we obtain $\det(A)=d \neq0$.
Since $(0^{2n},A) \sim L$, \ref{SubSubxy} implies that there exist
$(u,v) \in \Integ^2$ such that 
$ \Sub(0^{2n},A) = \Sub(L) + d (u,v)$.
On the other hand,
\ref{0iA} gives $\Sub(A) = (-1)^n \Sub(0^{2n},A)$, so 
$$
\Sub(A) = (-1)^n \big( \Sub(L) + d(u,v) \big).
$$
It follows that $\SUB(A) = (-1)^n \SUB(L)$
and that (b$'$) is true.
\end{proof}

\begin{corollary}\label{EquivSeq}
For $L,L'\in\Integ^*$, the following are equivalent:
\begin{enumerate}

\item $L \sim L'$

\item $\hodge{L} = \hodge{L'}$, $\det (L) = \det (L')$ and
$\SUB(L) = \SUB(L')$.

\end{enumerate}
\end{corollary}

\begin{proof}
Immediate consequence of \ref{uniqueness}.
\end{proof}

\begin{remark}
One can state some variants of \ref{EquivSeq}, for instance:
\begin{itemize}

\item Suppose that $L, L' \in \Integ^*$ satisfy
$\det(L) = 0 = \det(L')$. Then 
$$
L \sim L' \iff \hodge{L} = \hodge{L'}.
$$

\item Suppose that $L, L' \in \Integ^*$ satisfy
$\det(L) = d = \det(L')$ and $\hodge{L} = \hodge{L'}$.
Then 
$$
\textstyle
L \sim L' \iff \det_1(L) \equiv \det_1(L') \pmod d.
$$

\end{itemize}
\end{remark}

\medskip
\begin{definition}\label{SeqTransp}
Let $C = (0^r,A) \in \Integ^*$ be a canonical sequence
(where $r\in\Nat$ and $A\in\Neul^*$).
The \textit{transpose} $C^t$ of $C$ is defined by $C^t = (0^r,A^-)$.
Note that $C^t$ is a canonical sequence.
\end{definition}

\begin{lemma}\label{RevTransp}
Let $X \in \Integ^*$.
If $C$ is the unique canonical sequence equivalent to $X$,
then $C^t$ is the unique canonical sequence equivalent to $X^-$.
\end{lemma}

\begin{proof}
Since $X^- \sim C^-$, it suffices to show that $C^- \sim C^t$.
Write $C = (0^r,A)$ with $A\in\Neul^*$.
If $r$ is odd then $A=\emptyseq$ and the result holds trivially.
Assume that $r$ is even, then:
$$
C^- = (0^r,A)^- = (A^-, 0^r) \simref{\ref{MoveZeros}} (0^r, A^-)= C^t.
$$
\end{proof}

\section*{Further results on the classification of sequences}

We write $\Ceul\in \SEC$ to indicate that $\Ceul \subset \Integ^*$ is an
equivalence class of sequences.
Result \ref{BigPicture}, below,
gives a surprisingly simple description of the set
$\SEC$.

\smallskip
Given $\Ceul\in \SEC$, 
let
$
\min\Ceul = \setspec{ M \in \Ceul }
{ \text{$M$ is a minimal element of $\Integ^*$} }
$
denote the set of minimal elements of $\Ceul$
(see \ref{DefBupSeq} for the notion of minimal sequence).

\begin{lemma}\label{BasicSuccLemma}
Suppose that $M_1, M_2$ are minimal elements of $\Integ^*$
and $X_i \in M_i^{\oplus} $ $(i=1, 2)$. Then
$$
X_1 \sim X_2  \iff  M_1 \sim M_2.
$$
\end{lemma}

\begin{proof}
For each $i$ we have $X_i \sim (0,0,M_i)$ by \ref{LemZM},
so it suffices to prove:
$$
(0,0,A_1) \sim (0,0,A_2)\ \iff\ A_1 \sim A_2, \quad
\textit{for all $A_1, A_2 \in \Integ^*$.}
$$
If $A_1\sim A_2$ then $(0,0,A_1) \sim (0,0,A_2)$ by \ref{SubSeq}.
For the converse, note that \ref{0iA} and \ref{BasicHodge}  give
$\det(0,0,A_i) = - \det A_i$,
$\Sub(0,0,A_i) = - \Sub A_i$ and
$\hodge{(0,0,A_i)} = 1+ \hodge{ A_i}$ (for $i=1,2$), so
if $(0,0,A_1) \sim (0,0,A_2)$ we obtain $A_1\sim A_2$ by \ref{EquivSeq}.
\end{proof}

\begin{definition}\label{DefCoplus}
For each  $\Ceul \in \SEC$  we define an element
$\Ceul^{\oplus} $ of $\SEC$ as follows: Pick any minimal element $M$
of $\Ceul$, pick any $X\in M^{\oplus} $ and let $\Ceul^{\oplus} $
be the class of $X$.
By \ref{BasicSuccLemma}, $\Ceul^{\oplus}$ is well-defined and:
\begin{equation}\label{UniquePred}
\Ceul \longmapsto \Ceul^{\oplus} 
\ \ \text{is an injective map from $\SEC$ to itself.}
\end{equation}
We call $\Ceul^{\oplus} $ the \textit{successor\/} of $\Ceul$.
If $\Ceul = \Ceul_1^{\oplus} $ for some $\Ceul_1$ then $\Ceul_1$ is unique
by \eqref{UniquePred}; in this case we
say that ``$\Ceul$ has a predecessor'' and we call $\Ceul_1$
the \textit{predecessor\/} of $\Ceul$.
If $\Ceul \in \SEC$ then \ref{LemZM} gives
\begin{equation}\label{HodgeC+}
\det\big( \Ceul^\oplus \big) = - \det \Ceul
\quad \text{and} \quad
\hodge{\Ceul^\oplus} = 1+ \hodge{\Ceul}.
\end{equation}
\end{definition}

\begin{lemmadef}\label{LDprime}
For an element $\Ceul$ of $\SEC$, the following are equivalent:
\begin{enumerate}

\item $\min \Ceul$ is a singleton
\item $\min \Ceul$ is a finite set
\item $r\le1$, where $r$ is defined by the condition:
The unique canonical sequence belonging to $\Ceul$ is $(0^r,A)$,
where $r\in\Nat$ and $A\in\Neul^*$.
\item $\Ceul$ does not have a predecessor.

\end{enumerate}
{\rm If $\Ceul$ satisfies these conditions
then we call it a {\em prime class}.}
\end{lemmadef}

\begin{proof}
Note that $\neg\mbox{(2)} \Rightarrow \neg\mbox{(1)}$ is trivial; we
prove 
$\neg\mbox{(1)} \Rightarrow \neg\mbox{(4)}
\Rightarrow \neg\mbox{(3)} \Rightarrow \neg\mbox{(2)}$.

If $\hodge{\Ceul}=0$ then the canonical
element of $\Ceul$ is a sequence $X\in\Neul^*$; clearly, $X$ is then
the unique minimal element of $\Ceul$, so the condition $\hodge{\Ceul}=0$ 
implies (1). 

Hence, if (1) is false then $\hodge{\Ceul}>0$; since $\min\Ceul$
has more than one element, we may pick a minimal $X\in\Ceul$ such that
$X \neq (0)$; then \ref{ZM} gives $X\in M^{\oplus} $ for some
minimal element $M$ of $\Integ^*$. Thus (4) is false.

If (4) is false then we may
consider the unique canonical sequence $C$ which belongs to the
predecessor of $\Ceul$.
Then $(0,0,C)$ is the unique canonical sequence belonging to $\Ceul$,
so (3) is false.

If (3) is false then the canonical element  $(0^r,A)$ of $\Ceul$ satisfies
$r\ge2$.  By \ref{ChangeOddValue}, $(0,x,0^{r-2},A) \in \min\Ceul$ for every
$x\in \Integ\setminus\{-1\}$, so (2) is false.
\end{proof}

\begin{corollary}\label{ListPrimes}
The set of prime classes is 
$\{ \Ceul_0 \} \cup \setspec{ \Ceul_X }{ X \in \Neul^* }$,
where $\Ceul_0$ denotes the equivalence class of the sequence $(0)$
and, for each $X\in\Neul^*$, $\Ceul_X$ is the class of $X$.
\end{corollary}

\begin{proof}
Follows immediately from condition (3) of \ref{LDprime}.
\end{proof}

\medskip
\begin{nothing*}\label{Poset}
Given $\Ceul, \Ceul' \in \SEC$, write $\Ceul \preceq \Ceul'$ 
to indicate that there exists a sequence $\Ceul_0, \dots, \Ceul_n$
in $\SEC$ satisfying $n\in\Nat$, $\Ceul_0=\Ceul$,
$\Ceul_n=\Ceul'$ and $\Ceul_{i+1} = \Ceul_i^\oplus$
for all $i$ such that $0\le i<n$.
Then $\preceq$ is a partial order on the set $\SEC$ such
that the minimal elements are precisely the prime classes.
Given $\Ceul \in \SEC$, define the interval
$$
[ \Ceul, \infty) = \setspec {\Ceul' \in \SEC\ } {\ \Ceul \preceq \Ceul' }.
$$
\end{nothing*}

The following is now clear:

\begin{subproposition}\label{BigPicture}
Let $P = \setspec{ [ \Ceul, \infty) }{ \text{$\Ceul$ is a prime class}}$.
\begin{enumerate}

\item The set $P$ is a partition of $\SEC$.

\item If $I \in P$ then $( I, \preceq )$ is isomorphic to $(\Nat, \leq)$
as a partially ordered set.

\item If $I,I'$ are distinct elements of $P$,
and if $\Ceul \in I$ and $\Ceul' \in I'$,
then $\Ceul$ and $\Ceul'$ are not comparable w.r.t.\ $\preceq$. 

\end{enumerate}
\end{subproposition}

Note that \ref{ListPrimes} and \ref{BigPicture}
completely describe the partially ordered set
$\big( \SEC,\,  \preceq \big)$.

%

\section{Classification of linear chains}
\label{Sec:ClassLinChains}

It is clear that every result about sequences gives rise to a result
about linear chains.
This section reformulates \ref{SeqThm} and \ref{EquivSeq}
in terms of linear chains but leaves it
to the reader to translate the other results of
Section~\ref{SectionSequences}
(in particular \ref{uniqueness}, the remark after \ref{EquivSeq},
\ref{LDprime}--\ref{BigPicture}).

\begin{definition}\label{DefCanChain}
By a \textit{canonical chain}, we mean a linear chain of the form
$[L]$ where $L\in\Integ^*$ is a canonical sequence. 
The \textit{transpose\/} $\LL^t$ of a canonical chain $\LL$ is defined by:
$$
\LL^t = [L^t]
$$
where $L \in \Integ^*$ is a canonical sequence satisfying $\LL = [L]$
and where $L^t$ was defined in \ref{SeqTransp}.
Note that $\LL^t$ is a canonical chain.
\end{definition}


\begin{remark}
The linear chain $\LL^t$ is well-defined even when
$L$ is not uniquely determined by $\LL$
(i.e., when $L$ and $L^-$ are canonical and distinct).
Indeed, this only happens when $L\in\Neul^*$ and in that
case we have $[L^t]= [(L^-)^t]=\LL$.
\end{remark}

Concretely, a linear chain is canonical if it is $[0^r]$ with $r$ odd
or if it has the form:
\begin{equation}\label{CanChainPict1}
\raisebox{-6\unitlength}{%
\begin{picture}(54,10)(-2,-7)
\put(0,0){\circle*{1}}
\put(10,0){\makebox(0,0){\mbox{\,\dots}}}
\put(20,0){\circle*{1}}
\put(30,0){\circle*{1}}
\put(40,0){\makebox(0,0){\mbox{\,\dots}}}
\put(50,0){\circle*{1}}
\put(0,0){\line(1,0){6}}
\put(20,0){\line(-1,0){6}}
\put(20,0){\line(1,0){16}}
\put(50,0){\line(-1,0){6}}
\put(0,1.5){\makebox(0,0)[b]{$\scriptstyle 0$}}
\put(20,1.5){\makebox(0,0)[b]{$\scriptstyle 0$}}
\put(30,1.5){\makebox(0,0)[b]{$\scriptstyle a_1$}}
\put(50,1.5){\makebox(0,0)[b]{$\scriptstyle a_n$}}
\put(10,-1.5){\makebox(0,0)[t]%
{$\displaystyle\underbrace{\rule{22\unitlength}{0mm}}_{r \text{ vertices}}$}}
\end{picture}}
\qquad \text{($r\ge0$ is even, $n\ge0$ and  $\forall_i\ a_i\le-2$).}
\end{equation}
The transpose of $[0^r]$ is the same graph $[0^r]$ and the transpose
of \eqref{CanChainPict1} is:
\begin{equation*}
\raisebox{-6\unitlength}{%
\begin{picture}(54,10)(-2,-7)
\put(0,0){\circle*{1}}
\put(10,0){\makebox(0,0){\mbox{\,\dots}}}
\put(20,0){\circle*{1}}
\put(30,0){\circle*{1}}
\put(40,0){\makebox(0,0){\mbox{\,\dots}}}
\put(50,0){\circle*{1}}
\put(0,0){\line(1,0){6}}
\put(20,0){\line(-1,0){6}}
\put(20,0){\line(1,0){16}}
\put(50,0){\line(-1,0){6}}
\put(0,1.5){\makebox(0,0)[b]{$\scriptstyle 0$}}
\put(20,1.5){\makebox(0,0)[b]{$\scriptstyle 0$}}
\put(30,1.5){\makebox(0,0)[b]{$\scriptstyle a_n$}}
\put(50,1.5){\makebox(0,0)[b]{$\scriptstyle a_1$}}
\put(10,-1.5){\makebox(0,0)[t]%
{$\displaystyle\underbrace{\rule{22\unitlength}{0mm}}_{r \text{ vertices}}$}}
\end{picture}}.
\end{equation*}

As a corollary to the classification of sequences, we obtain:

\begin{theorem}\label{ChainsThm}
Every linear chain is equivalent to a canonical chain.
Moreover, if $\LL$ and $\LL'$ are canonical chains then
$$
\LL \sim \LL'  \iff  \LL' \in \big\{ \LL, \LL^t \big\}.
$$
\end{theorem}

\begin{proof}
In view of \ref{BasicBup}, 
this is a corollary to \ref{SeqThm} and \ref{RevTransp}.
\end{proof}

\rien{ 
\begin{proof}
Consider an arbitrary linear chain $[X]$, where $X\in\Integ^*$.
By \ref{SeqThm}, there is a unique canonical sequence $L\in\Integ^*$
such that $X\sim L$; then $[L]$ is a canonical chain and
$[X] \sim [L]$.

Let $\LL=[L]$ and $\LL'=[L']$ be canonical chains, where
$L,L'\in\Integ^*$ are canonical sequences.
As noted in \ref{BasicBup}, the condition $\LL\sim\LL'$ is equivalent to:
\begin{equation}\label{FirstForm}
\text{$L'\sim L$ or $L' \sim L^-$}.
\end{equation}
Since \ref{RevTransp} gives $L^- \sim L^t$,
condition~\eqref{FirstForm} is equivalent to:
\begin{equation}\label{SecondForm}
\text{$L'\sim L$ or $L' \sim L^t$}.
\end{equation}
Since $L', L, L^t$ are canonical sequences,
the uniqueness part of \ref{SeqThm} implies the equivalence of 
condition~\eqref{SecondForm} with:
\begin{equation*}
\text{$L'= L$ or $L' = L^t$}.
\end{equation*}
This proves that 
$\LL \sim \LL'  \iff  \LL' \in \big\{ \LL, \LL^t \big\}$.
\end{proof}
} 

For the next result, we need:

\begin{definition}
Let $\LL$ be a linear chain.
Define a subset $\Sub(\LL)$ of $\Integ$ as follows:
Choose $L\in\Integ^*$ such that $\LL = [L]$,
let $(x,y) = \Sub(L) \in \Integ\times\Integ$ and set
$$
\Sub(\LL) = \{x, y \}.
$$
We also define the subset $\SUB(\LL)$ of $\Integ/d\Integ$,
where $d=\det(\LL)$, 
by taking the image of $\Sub(\LL)$ via the canonical epimorphism
$\Integ \to \Integ/d\Integ$.
\end{definition}

\begin{corollary}\label{Isoms}
For linear chains $\LL$ and $\LL'$, the following are equivalent:
\begin{enumerate}

\item $\LL \sim \LL'$

\item $\hodge{\LL} = \hodge{\LL'}$,
$\det\LL = \det\LL'$ and
$\SUB(\LL) \cap \SUB(\LL') \neq \emptyset$.

\end{enumerate}
\end{corollary}

\begin{proof}
Follows immediately from \ref{EquivSeq}.
Note that the condition
$\SUB(\LL) \cap \SUB(\LL') \neq \emptyset$
is equivalent to $\SUB(\LL) = \SUB(\LL')$
by \ref{Claim1Disc}.
\end{proof}

\section{$\tau$-equivalence of sequences}
\label{Section:TauEquiv}

\begin{definition}\label{DefBup}
The sentence \textit{``$\tau$ is a type''} means that $\tau$ is one of
the four symbols $(-,-)$, $(+,-)$, $(-,+)$, $(+,+)$.
For each type $\tau$, we define a subset $\Integ^*_\tau$ of $\Integ^*$
and an equivalence relation $\tequiv{\tau}$ on $\Integ^*_\tau$.
The sets $\Integ^*_\tau$ are defined by:
\begin{gather*}
\Integ^*_{(-,-)} = \Integ^* \\
\Integ^*_{(+,-)} = \Integ^*_{(-,+)}
= \setspec{ (x_1, \dots, x_n) \in \Integ^* }{ n\ge1 }
=\Integ^* \setminus \{ \emptyseq \} \\
\Integ^*_{(+,+)} 
= \setspec{ (x_1, \dots, x_n) \in \Integ^* }{ n\ge2 }.
\end{gather*}
Let $X = (x_1, \dots, x_n) \in \Integ^*$.
\begin{enumerate}

\item By a \textit{$(-,-)$-blowing-up\/} of $X$, we mean a blowing-up
of $X$ in the sense of \ref{DefBupSeq}.

\item If $X \in \Integ^*_{(+,-)}$, then any of the following
is called a \textit{$(+,-)$-blowing-up\/} of $X$:
	\begin{itemize}

	\item $X' =
	(x_1, \dots, x_{i-1}, x_i - 1, -1, x_{i+1}-1, x_{i+2}, \dots, x_n)$,
	for some $i$ such that $1\le i < n$;

	\item $X' =
	(x_1, \dots, x_{n-1}, x_n - 1, -1)$.

	\end{itemize}

\item If $X \in \Integ^*_{(-,+)}$, then any of the following
is called a \textit{$(-,+)$-blowing-up\/} of $X$:
	\begin{itemize}

	\item $X' =
	(x_1, \dots, x_{i-1}, x_i - 1, -1, x_{i+1}-1, x_{i+2}, \dots, x_n)$,
	for some $i$ such that $1 \le i < n$;

	\item $X' = (-1, x_1-1, x_2, \dots, x_{n})$.

	\end{itemize}

\item If $X \in \Integ^*_{(+,+)}$, then any of the following
is called a \textit{$(+,+)$-blowing-up\/} of $X$:
	\begin{itemize}

	\item $X' =
	(x_1, \dots, x_{i-1}, x_i - 1, -1, x_{i+1}-1, x_{i+2}, \dots, x_n)$,
	for some $i$ such that $1 \le i < n$.

	\end{itemize}

\end{enumerate}
Note that if $X \in \Integ^*_\tau$ and $X'$ is a $\tau$-blowing-up of
$X$ then $X' \in \Integ^*_\tau$; in this situation, we also say that
$X$ is a \textit{$\tau$-blowing-down\/} of $X'$.
Two elements of $ \Integ^*_\tau$ are \textit{$\tau$-equivalent\/} if
one can be obtained from the other by a finite sequence of
$\tau$-blowings-up and $\tau$-blowings-down.
We write $X \tequiv{\tau} X'$ for $\tau$-equivalence.
\end{definition}

\begin{remark}
The theory of $(-,-)$-equivalence is exactly the content of
section~\ref{SectionSequences}.
In fact we will rarely use the notation $X \tequiv{(-,-)} X'$,
since it simply means $X \sim X'$.
\end{remark}

\bigskip
Before developing the theory of $\tau$-equivalence, let us explain
how it will be used.

\begin{example}
Consider the weighted tree
$$
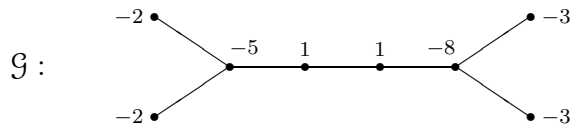
\GG : \qquad
\raisebox{-8\unitlength}{\begin{picture}(60,18)(-5,-9)
\put(0,6.66666){\circle*{1}}
\put(0,-6.66666){\circle*{1}}
\put(10,0){\circle*{1}}
\put(20,0){\circle*{1}}
\put(30,0){\circle*{1}}
\put(40,0){\circle*{1}}
\put(50,6.66666){\circle*{1}}
\put(50,-6.66666){\circle*{1}}
\put(10,0){\line(1,0){30}}
\put(10,0){\line(-3,2){10}}
\put(10,0){\line(-3,-2){10}}
\put(40,0){\line(3,2){10}}
\put(40,0){\line(3,-2){10}}
\put(-1.5,6.66666){\makebox(0,0)[r]{$\scriptstyle -2$}}
\put(-1.5,-6.66666){\makebox(0,0)[r]{$\scriptstyle -2$}}
\put(51.5,6.66666){\makebox(0,0)[l]{$\scriptstyle -3$}}
\put(51.5,-6.66666){\makebox(0,0)[l]{$\scriptstyle -3$}}
\put(10,1.5){\makebox(0,0)[bl]{$\scriptstyle -5$}}
\put(20,1.5){\makebox(0,0)[b]{$\scriptstyle 1$}}
\put(30,1.5){\makebox(0,0)[b]{$\scriptstyle 1$}}
\put(40,1.5){\makebox(0,0)[br]{$\scriptstyle -8$}}
%
\end{picture}}
$$
Consider the path $\gamma$ in $\GG$ which goes from the vertex of
weight $-5$ to that of weight $-8$; note that $\gamma \in P(\GG)$ is
of type $(+,+)$.
One can see (e.g.\ by \ref{ZeroProp}, below) that 
$(-5,1,1,-8) \tequiv{(+,+)} (-6+x, 0,0,0, -8-x)$
for any $x\in\Integ$.
Then it follows that $\GG \sim \GG'$,
where
$$
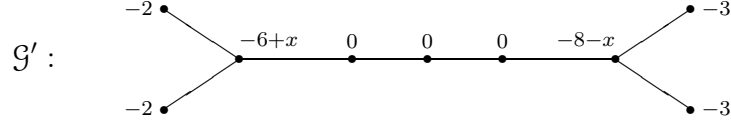
\GG' : \qquad
\raisebox{-8\unitlength}{\begin{picture}(80,18)(-10,-9)
\put(-5,6.66666){\circle*{1}}
\put(-5,-6.66666){\circle*{1}}
\put(5,0){\circle*{1}}
\put(20,0){\circle*{1}}
\put(30,0){\circle*{1}}
\put(40,0){\circle*{1}}
\put(55,0){\circle*{1}}
\put(65,6.66666){\circle*{1}}
\put(65,-6.66666){\circle*{1}}
\put(5,0){\line(1,0){50}}
\put(5,0){\line(-3,2){10}}
\put(5,0){\line(-3,-2){10}}
\put(55,0){\line(3,2){10}}
\put(55,0){\line(3,-2){10}}
\put(-6.5,6.66666){\makebox(0,0)[r]{$\scriptstyle -2$}}
\put(-6.5,-6.66666){\makebox(0,0)[r]{$\scriptstyle -2$}}
\put(66.5,6.66666){\makebox(0,0)[l]{$\scriptstyle -3$}}
\put(66.5,-6.66666){\makebox(0,0)[l]{$\scriptstyle -3$}}
\put(5,1.5){\makebox(0,0)[bl]{$\scriptstyle -6+x$}}
\put(20,1.5){\makebox(0,0)[b]{$\scriptstyle 0$}}
\put(30,1.5){\makebox(0,0)[b]{$\scriptstyle 0$}}
\put(40,1.5){\makebox(0,0)[b]{$\scriptstyle 0$}}
\put(55,1.5){\makebox(0,0)[br]{$\scriptstyle -8-x$}}
%
\end{picture}}
$$
Indeed, the definition of $(+,+)$-equivalence takes into account
that a graph cannot be blown-down at a vertex of degree greater than
two.
\end{example}

More generally, we state the following (trivial) fact:

\begin{nothing}\label{ApplTau}
Let $\GG$ be a weighted forest,
let $\gamma = (u_1, \dots, u_m) \in P(\GG)$
and let $\tau$ be the type of $\gamma$.
Write $x_i = w(u_i,\GG)$  for $i = 1, \dots, m$ and
suppose that $(y_1, \dots, y_n) \in \Integ_\tau^*$ satisfies
$$
(x_1, \dots, x_m) \tequiv{\tau} (y_1, \dots, y_n).
$$
Consider the weighted forest $\GG'$ obtained from $\GG$ by replacing
\begin{equation*}
\raisebox{-1\unitlength}{
\begin{picture}(34,6)(-2,-3)
\put(0,0){\circle*{1}}
\put(10,0){\circle*{1}}
\put(20,0){\makebox(0,0){\mbox{\,\dots}}}
\put(30,0){\circle*{1}}
\put(0,0){\line(1,0){16}}
\put(30,0){\line(-1,0){6}}
\put(0,1.5){\makebox(0,0)[b]{$\scriptstyle x_1$}}
\put(10,1.5){\makebox(0,0)[b]{$\scriptstyle x_2$}}
\put(30,1.5){\makebox(0,0)[b]{$\scriptstyle x_m$}}
\put(0, -1.5){\makebox(0,0)[t]{$\scriptstyle u_1$}}
\put(10,-1.5){\makebox(0,0)[t]{$\scriptstyle u_2$}}
\put(30,-1.5){\makebox(0,0)[t]{$\scriptstyle u_m$}}
\end{picture}}
\qquad
\text{by}
\qquad
\raisebox{-1\unitlength}{
\begin{picture}(34,6)(-2,-3)
\put(0,0){\circle*{1}}
\put(10,0){\circle*{1}}
\put(20,0){\makebox(0,0){\mbox{\,\dots}}}
\put(30,0){\circle*{1}}
\put(0,0){\line(1,0){16}}
\put(30,0){\line(-1,0){6}}
\put(0,1.5){\makebox(0,0)[b]{$\scriptstyle y_1$}}
\put(10,1.5){\makebox(0,0)[b]{$\scriptstyle y_2$}}
\put(30,1.5){\makebox(0,0)[b]{$\scriptstyle y_n$}}
%
%
\end{picture}}
\end{equation*}
and leaving the rest of the graph unchanged.  Then $\GG \sim \GG'$.
\end{nothing}

\bigskip
We shall now prove some properties of $\tau$-equivalence.
\begin{remark}
If $X, X' \in \Integ^*_{(+,-)} = \Integ^*_{(-,+)}$
then $X \tequiv{(+,-)} X' \iff X^- \tequiv{(-,+)} (X')^-$.
So any result regarding $(+,-)$-equivalence gives rise to
a result on $(-,+)$-equivalence, and vice-versa.
\end{remark}

Given $X,Y \in \Integ^*$, it is obvious that 
$X \tequiv{(+,+)} Y$ implies both 
$X \tequiv{(+,-)} Y$ and $X \tequiv{(-,+)} Y$,
and that ``$X \tequiv{(+,-)} Y$ or $X \tequiv{(-,+)} Y$''
implies $X \sim Y$.  Here is a slightly less obvious fact:

\begin{lemma}\label{LessObvious}
Let $X,Y\in \Integ^*$ and $a,b,\alpha,\beta \in \Integ$.
\begin{enumerate}

\item If $(a,X,b) \tequiv{(+,+)} (\alpha,Y,\beta)$
then $(a,X) \tequiv{(+,-)} (\alpha,Y)$
and $(X,b) \tequiv{(-,+)} (Y,\beta)$.

\item If $(a,X) \tequiv{(+,-)} (\alpha,Y)$ then $X \sim Y$.

\item If $(X,b) \tequiv{(-,+)} (Y,\beta)$ then $X \sim Y$.

\end{enumerate}
\end{lemma}

\begin{proof}
To prove (1), it suffices to consider the case where $(\alpha,Y,\beta)$
is a $(+,+)$-blowing-up of $(a,X,b)$; but the assertion is easily
verified in this case, so (1) is true.  Similar arguments prove (2) and (3).
\end{proof}

\begin{lemma}\label{ForSuitable}
Let $X,Y\in\Integ^*$ be such that $X\sim Y$. Then:
\begin{gather*}
\forall_{a,b\in\Integ}\ \exists_{\alpha,\beta\in\Integ}\ 
(a,X,b) \tequiv{(+,+)} ( {\alpha}, Y, {\beta}), \\
\forall_{a\in\Integ}\ \exists_{\alpha\in\Integ}\ 
(a,X) \tequiv{(+,-)} (\alpha, Y)
\qquad
\text{and}
\qquad
\forall_{b\in\Integ}\ \exists_{\beta\in\Integ}\ 
(X,b) \tequiv{(-,+)} (Y, \beta).
\end{gather*}
%
%
%
%
%
\end{lemma}

\begin{proof}
Let $a,b\in\Integ$.
Since $X\sim Y$, there exists a sequence $s$ of blowings-up and
blowings-down of sequences which transforms $X$ into $Y$.
If we now replace $X$ by $(a,X,b)$ and perform
the ``same'' sequence $s$ of operations, we get 
a sequence $s'$ of $(+,+)$-blowings-up and $(+,+)$-blowings-down
which transforms
$(a,X,b)$ into $(\alpha, Y, \beta)$,
for suitable $\alpha,\beta\in\Integ$ (this claim is obvious in the
case where $s$ consists of a single blowing-up or a single
blowing-down, and we may reduce to that case). 
Thus $(a,X,b) \tequiv{(+,+)} (\alpha,Y,\beta)$.
The same argument shows that
$(a,X) \tequiv{(+,-)} (\alpha, Y)$ for some $\alpha$
and $(X,b) \tequiv{(-,+)} (Y,\beta)$ for some $\beta$.
\end{proof}

\begin{lemma}\label{AllIJ}
Let $X,Y\in\Integ^*$ and $a,b,\alpha,\beta \in \Integ$.
Then:
\begin{enumerate}

\item $(a,X,b) \tequiv{(+,+)} (\alpha, Y, \beta)
\ \implies\ \forall_{i,j\in\Integ}\ 
(a+i, X, b+j) \tequiv{(+,+)} (\alpha+i, Y, \beta+j)$

\item $(a,X) \tequiv{(+,-)} (\alpha, Y)
\ \implies\ \forall_{i\in\Integ}\ 
(a+i, X) \tequiv{(+,-)} (\alpha+i, Y)$

\item $(X,b) \tequiv{(-,+)} (Y,\beta)
\ \implies\ \forall_{j\in\Integ}\ 
(X, b+j) \tequiv{(-,+)} (Y,\beta+j)$

\end{enumerate}
\end{lemma}

\begin{proof}
Suppose that
$(a, X, b) \tequiv{(+,+)} (\alpha, Y, \beta)$.
Then there exists a sequence of $(+,+)$-blowings-up and
$(+,+)$-blowings-down which transforms 
$(a, X, b)$ into $(\alpha, Y, \beta)$;
clearly, the same sequence of operations applied to 
$(a+i, X, b+j)$ yields
$(\alpha+i, Y, \beta+j)$
(this claim is obvious in the case of a single $(+,+)$-blowing-up or a single
$(+,+)$-blowing-down, and we may reduce to that case). 
This proves (1) and the
other assertions are proved by the same argument.
\end{proof}

For the next result, we need:

\begin{definition}\label{DefdeltaXY}
Let $X,Y\in\Integ^*$ be equivalent sequences, let $d = \det(X) =  \det(Y)$
and $n= \hodge{X} = \hodge{Y}$.
Define the integer $\delta(X,Y)$ by:
$$
\delta(X,Y) =
\begin{cases}
\frac 1d \big( \det_1(X) - \det_1(Y) \big), & \text{if $d\neq0$;} \\
(-1)^{n-1} \big( \det_*(X) - \det_*(Y) \big), & \text{if $d = 0$.}
\end{cases}
$$
\end{definition}

\begin{proposition}\label{ZeroProp}
Let $X,Y\in\Integ^*$ be such that $X\sim Y$ and let
$a,b,\alpha,\beta \in \Integ$.
\begin{enumerate}

\item If $\det(X)\neq0$ then 
\begin{enumerate}

\item $(a, X, b) \tequiv{(+,+)} (\alpha, Y, \beta)
\iff  \alpha = a + \delta(X,Y) \text{ and }  \beta = b + \delta(X^-,Y^-)$

\item $(a,X) \tequiv{(+,-)} (\alpha, Y)
\iff  \alpha = a + \delta(X,Y)$

\item $(X,b) \tequiv{(-,+)} (Y, \beta)
\iff  \beta = b + \delta(X^-,Y^-)$.

\end{enumerate}

\item If $\det(X)=0$ then 
\begin{enumerate}

\item $(a,X,b) \tequiv{(+,+)} (\alpha, Y, \beta)
\iff \alpha + \beta = a+b + \delta(X,Y)$

\item For all $i,j \in\Integ$, 
$(i,X) \tequiv{(+,-)} (j,Y)$
and 
$(X,i) \tequiv{(-,+)} (Y,j)$.

\end{enumerate}

\end{enumerate}
\end{proposition}

\begin{proof}
Throughout, we fix $X,Y\in\Integ^*$ such that $X\sim Y$.
Let $d=\det(X)$.

Assume that $d\neq0$.

If $a,\alpha\in\Integ$ are such that
$(a,X) \tequiv{(+,-)} (\alpha,Y)$ then $(a,X) \sim (\alpha,Y)$
and consequently $\det(a,X) = \det(\alpha,Y)$;
then \ref{remark*} gives
$$
\textstyle
(-a)d -\det_1(X)= \det(a,X)= \det(\alpha,Y)= (-\alpha)d -\det_1(Y),
$$
from which $\alpha = a + \delta(X,Y)$ follows. This proves
implication ``$\Rightarrow$'' of assertion (1b).
Conversely, let $a\in\Integ$ and $\alpha = a + \delta(X,Y)$.
By \ref{ForSuitable}, 
$(a,X) \tequiv{(+,-)} (\alpha_1, Y)$ for some $\alpha_1$,
and by ``$\Rightarrow$'' of (1b) we have 
$\alpha_1 = a + \delta(X,Y)=\alpha$.
So assertion~(1b) is true.

If $b,\beta\in\Integ$ then
$(X,b) \tequiv{(-,+)} (Y,\beta)
\Leftrightarrow
(X,b)^- \tequiv{(+,-)} (Y,\beta)^-
\Leftrightarrow
(b,X^-) \tequiv{(+,-)} (\beta,Y^-)$; so (1b) implies (1c).

If $a,b,\alpha,\beta\in\Integ$ are such that
$(a,X,b) \tequiv{(+,+)} (\alpha, Y, \beta)$
then, by \ref{LessObvious},
both 
$(a,X) \tequiv{(+,-)} (\alpha, Y)$
and
$(X,b) \tequiv{(-,+)} (Y,\beta)$ hold;
then (1b) and (1c) imply that
$\alpha = a + \delta(X,Y)$ and $\beta = b + \delta(X^-,Y^-)$,
so implication ``$\Rightarrow$'' of (1a) is proved.
The converse follows from ``$\Rightarrow$'' and the
fact (\ref{ForSuitable}) that 
$(a,X,b) \tequiv{(+,+)} (\alpha_1, Y, \beta_1)$ for
some $\alpha_1,\beta_1$ (see the proof of (1b)).
This completes the proof of assertion (1).

\medskip
To prove assertion~(2), assume that $\det(X)=0$.
Let $n = \hodge{X} = \hodge{Y}$;
then \ref{uniqueness} implies that $n\ge1$ and that $X\sim Z \sim Y$,
where $Z=(0^{2n-1})\in\Integ^*$.
If $a,b\in\Integ$ then
\begin{multline*}
(a,Z,b)
= (a,0^{2n-1},b)
\tequiv{(+,+)} (a-1,-1,-1, 0^{2n-2},b)
\simref{\ref{MoveZeros}} (a-1,0^{2n-2},-1,-1,b) \\
\tequiv{(+,+)} (a-1,0^{2n-2},0,b+1)
=(a-1,Z,b+1)
\end{multline*}
and it is easily verified that the equivalence given by
\ref{MoveZeros} is actually a $(+,+)$-equivalence.
So $(a,Z,b) \tequiv{(+,+)}  (a-1, Z, b+1)$
and consequently:
\begin{equation}\label{abimplies}
\forall_{a,b,\alpha,\beta \in \Integ}\quad
a+b=\alpha+\beta \implies
(a,Z,b) \tequiv{(+,+)} (\alpha,Z,\beta).
\end{equation}
By \ref{LessObvious},
$(a,Z,b) \tequiv{(+,+)} (\alpha,Z,\beta)$
implies both
$(a,Z) \tequiv{(+,-)} (\alpha,Z)$
and
$(Z,b) \tequiv{(-,+)} (Z,\beta)$.
So \eqref{abimplies} implies:
\begin{equation}\label{aimplies}
\forall_{i,j \in \Integ}\quad
(i,Z) \tequiv{(+,-)} (j,Z)
\text{ and }
(Z,i) \tequiv{(-,+)} (Z,j).
\end{equation}

\medskip
By \ref{ForSuitable} we have
$(0,X) \tequiv{(+,-)} (r,Z)$
and
$(0,Y) \tequiv{(+,-)} (s,Z)$
for suitable $r,s\in\Integ$;
so given any $i,j \in\Integ$ we have
$$
(i,X) \tequiv{(+,-)} (r+i,Z) \tequiv{(+,-)} (s+j, Z)
\tequiv{(+,-)}
(j,Y)
$$
by \ref{AllIJ} and \eqref{aimplies},
and similarly $(X,i) \tequiv{(-,+)} (Y,j)$,
which proves assertion (2b).

There remains only to prove (2a). 
In view of \ref{SameSUB}, the condition
$X\sim Z$ implies that
$\det_1(X)$ is congruent to $\det_1(Z)$ modulo $\det(X)$;
so $\det_1(X) = \det_1(Z) = \det(0^{2n-2}) = (-1)^{n-1}$
(by \ref{0iA}) and for the same reason we have:
$$
\textstyle
\det_1(X) = 
\det_1(X^-) = 
\det_1(Y) = 
\det_1(Y^-) = (-1)^{n-1}.
$$

Let $a,b,\alpha,\beta \in \Integ$ be such that
$(a,X,b) \tequiv{(+,+)} (\alpha, Y, \beta)$.
Then $(a,X,b) \sim (\alpha,Y,\beta)$
and consequently
\begin{equation}\label{aXbaYb}
\det (a,X,b) = \det(\alpha,Y,\beta).
\end{equation}
Using \ref{remark*} three times gives:
\begin{align}
\label{aXb}
\det(a,X,b) & \textstyle = (-a) \det(X,b) - \det_1(X,b), \\
\label{Xb}
\det(X,b) & \textstyle =
\det(b,X^-) =  (-b) \det(X^-) - \det_1(X^-) =  - \det_1(X^-)
= (-1)^n, \\
\textstyle
\label{1Xb}
\det_1(X,b) & = {\textstyle \det(x_2,\dots,x_m,b)
= (-b) \det_1(X) - \det_*(X) = (-1)^n b - \det_*(X)}
\end{align}
where we wrote $X = (x_1,\dots,x_m)$.
Substituting \eqref{Xb} and \eqref{1Xb} into \eqref{aXb} gives
$ \textstyle \det(a,X,b) = (-1)^{n-1} (a+b) + \det_*(X) $
and similarly we obtain
$ \textstyle \det(\alpha,Y,\beta) = (-1)^{n-1} (\alpha+\beta) + \det_*(Y)$.
So \eqref{aXbaYb} gives 
$(-1)^{n-1} (a+b) + \det_*(X) = (-1)^{n-1} (\alpha+\beta) + \det_*(Y)$
and consequently
$$
\textstyle
\alpha + \beta = (a+b) + (-1)^{n-1} (\det_*(X) - \det_*(Y))
= a+b + \delta(X,Y).
$$

Conversely, let $a,b,\alpha,\beta \in \Integ$ be such that
$\alpha + \beta = a+b + \delta(X,Y)$; we show that
$(a,X,b) \tequiv{(+,+)} (\alpha,Y,\beta)$.
By \ref{ForSuitable}, we have
$(a,X,b) \tequiv{(+,+)} (\alpha_0, Y, \beta_0)$
for suitable $\alpha_0, \beta_0 \in\Integ$;
so the above paragraph implies that
$\alpha_0 + \beta_0 = a+b + \delta(X,Y)$, hence
$\alpha_0 + \beta_0 =\alpha+\beta$.
By \ref{ForSuitable} we have
$(0,Y,0) \tequiv{(+,+)} (r,Z,s)$ for suitable $r,s\in\Integ$, so
$$
(\alpha_0, Y, \beta_0) \tequiv{(+,+)} (r+\alpha_0, Z, s+\beta_0)
\tequiv{(+,+)} (r+\alpha, Z, s+\beta)
\tequiv{(+,+)}
(\alpha, Y, \beta)
$$
by \ref{AllIJ} and \eqref{abimplies}.
Since
$(a,X,b) \tequiv{(+,+)} (\alpha_0, Y, \beta_0)$,
we conclude that
$(a,X,b) \tequiv{(+,+)} (\alpha, Y, \beta)$
and the proof is complete.
\end{proof}

\section{Minimal sequences and linear chains}
\label{Sec:FurtherLin}

See \ref{DefMinWG} for the notion of minimal weighted graph
and \ref{DefBupSeq} for equivalence and minimality in $\Integ^*$.
Consider the following:

\begin{problem}\label{MinWG}
List all minimal elements of a given equivalence class of weighted
forests.
\end{problem}

In Section~\ref{Sec:MinimalReduc}, it is shown that 
the above problem reduces to the special case of linear chains,
or more precisely to:

\begin{problem}\label{MinSQ}
Given $X\in\Integ^*$, list all minimal sequences equivalent to $X$.
\end{problem}

Apparently,
very little is known about these problems.
One notable exception is \cite{Morrow:MinNorm}, which can be interpreted as
solving Problem~\ref{MinSQ} for $X = (1)$.

\medskip
This section begins by solving Problem~\ref{MinSQ} recursively (\ref{RecSol});
together with \ref{MoplusExpl},
and keeping in mind \ref{BigPicture},
this gives substantial information about Problem~\ref{MinSQ}. 
Then we make use of those results to describe explicitely all minimal
elements of certain classes $\Ceul \in \SEC$;
in fact we can do this when $\Ceul$ is either a prime
class (\ref{LDprime}) or the successor of a prime class.
Finally, we show (\ref{IdentGeom}) that the cases that
we can describe explicitely 
are precisely those which arise in the study of algebraic surfaces.

The notations $\SEC$ and $\min(\Ceul)$ are defined before
\ref{BasicSuccLemma}.

\begin{proposition}\label{RecSol}
If $\Ceul \in \SEC$ then
$\displaystyle
\min\big( \Ceul^{\oplus}  ) =  \bigcup_{M \in \min\Ceul}  M^{\oplus}$.
\end{proposition}

\begin{proof}
The inclusion ``$ \supseteq $'' is trivial by definition
\ref{DefCoplus} of $\Ceul^\oplus$.
Consider $Z \in \min(\Ceul^\oplus)$.
Since $\Ceul^\oplus$ has a predecessor (namely $\Ceul$) but $\Ceul_0$
doesn't by \ref{ListPrimes}, we have $\Ceul^\oplus\neq\Ceul_0$ and hence
$Z \neq (0)$; we also have $\hodge Z > 0$ by \eqref{HodgeC+};
so \ref{ZM} gives $Z\in M^{\oplus} $ for some
minimal element $M$ of $\Integ^*$. We have $M \in \Ceul$ by uniqueness
of the predecessor of $\Ceul^\oplus$,
so $Z \in \cup_{M \in \min\Ceul}  M^{\oplus}$.
\end{proof}

In order to derive explicit results from \ref{RecSol},
we need to describe the elements of $M^\oplus$
where $M$ is a minimal element of $\Integ^*$.
In other words, we have to describe the sequences
$M'$ which occur in part (4) of \ref{M+}.
Some preliminary work is needed.
If $x\in\Reals$, let $\lceil x \rceil$ denote the least integer $n$ such
that $x\le n$.

\begin{lemma}\label{variation}
Let $\Meul =
\setspec{(x_1,\dots,x_n) \in \Integ^*\setminus \{ \emptyseq \}}
{x_1 \neq -1 \text{ and } \forall_{i>1}\, x_i \le -2 }$
and $\Meul^- = \setspec{ X^- }{ X \in \Meul }$.
Then $X \mapsto ( \det X, \det_1X )$ is a bijection from $\Meul$ to
\begin{equation}\label{TheSetcp}
\setspec{(r_0,r_1) \in \Integ^2}{r_1>0,\ \gcd(r_0,r_1)=1 \text{ and }
\big\lceil {\textstyle \frac {r_0}{r_1} \big\rceil } \neq1 }
\end{equation}
and $X \mapsto ( \det(X^-), \det_1(X^-) )$ is a bijection
from $\Meul^-$ to \eqref{TheSetcp}.
\end{lemma}

\begin{proof}
It is well known that
$\gcd\big( \det(X) , \det_1(X) \big) = 1$
holds for every $X \in \Integ^*$.
Consider an element $X = (-q,N)$ of $\Integ^* \setminus \{ \emptyseq \}$,
where $q\in\Integ$ and $N\in\Integ^*$. By \ref{remark*},
$$
\textstyle
\det(X) = q \det(N) - \det_1(N).
$$
If $N\in\Neul^*$ then by \ref{DetAdmiss} we have
$0 \le \det_1(N) < \det(N)$,
so $q = \big\lceil \frac {\det X}{\det N} \big\rceil$;
thus $\Meul$ is mapped into the set \eqref{TheSetcp}.
If $(r_0,r_1)$ belongs to the set \eqref{TheSetcp}, there is a unique
pair $(q,r_2) \in \Integ^2$ such that $r_0 = qr_1-r_2$ and $0\le r_2 < r_1$;
by \ref{DetAdmiss}, a unique $N\in\Neul^*$ satisfies
$\det(N)=r_1$ and $\det_1(N)=r_2$; then $(-q,N)\in\Meul$ and this defines a map
from the set \eqref{TheSetcp} to $\Meul$.
It is clear that the two maps are inverse of each other, so
the first assertion is proved.
The second assertion follows from the first.
\end{proof}

\begin{definition}\label{NewDefEalpha}
Let $Z,Z' \in \Integ^*$.
We say that $Z$ can be \textit{$(+,-)$-contracted\/} to $Z'$
(resp.\ \textit{$(-,+)$-contracted, $(+,+)$-contracted})
if there exists a sequence of blowings-down which transforms $Z$ into
$Z'$ and such that no blowing-down is performed at the leftmost
(resp.\ rightmost, leftmost or rightmost) term of a sequence.
Observe that this condition is stronger than
$Z \tequiv{\tau} Z'$ (see \ref{DefBup}),
where $\tau = (+,-)$ (resp.\ $\tau = (-,+)$, $\tau = (+,+)$).
\end{definition}

\begin{definition}\label{DefEalpha}
Given $\alpha,\beta\in\Integ$,
define the following subsets of $\Nat^3$:
\begin{align*}
P &= \setspec{ (n,p,c) \in \Nat^3 } { 1\le p \le c \text{ and } \gcd(c,p)=1 }\\
{}^\alpha\!P &= P^\alpha = \setspec{ (n,p,c) \in \Nat^3 }
{ 1 \le p \le c,\ \gcd(p,c)=1 \text{ and }
\big\lceil{\textstyle \frac{c}{nc+p} } \big\rceil \neq \alpha+1} \\
{}^\alpha\!P^\beta &= \setspec{ (n,p,c) \in \Nat^3 }
{ 1 \le p \le c,\ \gcd(p,c)=1,\ 
\big\lceil{\textstyle \frac{c}{nc+p} } \big\rceil \neq \alpha+1
\text{ and } n \neq \beta }
\end{align*}
and the following subsets of $\Integ^* \times \Integ^*$:
\begin{align*}
E &= \setspec{ (X,Y)\in\Neul^* \times \Neul^* } {(X,-1,Y) \sim
\emptyseq }\\
{}^\alpha\!E &= \setspec{ (X,Y) \in \Meul \times \Neul^* }
	{ (X,-1,Y) \text{ can be $(+,-)$-contracted to } (\alpha) } \\
E^\alpha  &= \setspec{ (X,Y) \in \Neul^* \times \Meul^- }
	{ (X,-1,Y) \text{ can be $(-,+)$-contracted to } (\alpha) } \\
{}^\alpha\!E^\beta &= \setspec{ (X,Y) \in \Meul \times \Meul^- }
	{ (X,-1,Y) \text{ can be $(+,+)$-contracted to } (\alpha,\beta) }
\end{align*}
where $\Meul$ and $\Meul^-$ are defined in \ref{variation}.
Then define four maps
\begin{enumerate}

\item $f : P \to E$, $(n,p,c) \mapsto (X,Y)$,

\item ${}^\alpha\!f : {}^\alpha\!P \to {}^\alpha\!E$, 
$(n,p,c) \mapsto (X,Y)$

\item $f^\alpha : P^\alpha \to E^\alpha$,
$(n,p,c) \mapsto (X,Y)$

\item ${}^\alpha\!f^\beta  \ :\  {}^\alpha\!P^\beta \to {}^\alpha\!E^\beta$,
$(n,p,c) \mapsto (X,Y)$
\end{enumerate}
by declaring in each case that
$(X,Y)$ is the unique pair of sequences satisfying:
\begin{enumerate}

\item[$(1')$]
$(X,Y) \in \Neul^* \times \Neul^*$ and:
\begin{align*}
\det(X) &= nc+p & \det(Y^-) &= c \\
\textstyle \det_1(X)&\equiv -c \pmod{ nc+p } & \textstyle \det_1(Y^-) &= c-p
\end{align*}

\item[$(2')$]
$(X,Y) \in \Meul \times \Neul^*$ and:
\begin{align*}
\det(X) &= c - \alpha(nc+p) & \det(Y^-) &= c \\
\textstyle \det_1(X)&= nc+p & \textstyle \det_1(Y^-) &= c-p
\end{align*}

\item[$(3')$]
$(X,Y) \in \Neul^* \times \Meul^-$ and: 
\begin{align*}
\det(X) &= c & \det(Y^-) &= c - \alpha(nc+p) \\
\textstyle \det_1(X) &= c-p & \textstyle \det_1(Y^-)&= nc+p
\end{align*}

\item[$(4')$]
$(X,Y) \in \Meul \times \Meul^-$ and: 
\begin{align*}
\det(X) &= c-\alpha(nc+p)	&\det(Y^-) &= (n-\beta)c +p	 \\
\textstyle \det_1(X) &= nc+p	&	\textstyle \det_1(Y^-) &= c.
\end{align*}

\end{enumerate}
\end{definition}

\begin{sublemma}\label{3bij}
The maps $f$, ${}^\alpha\!f$, $f^\alpha$ and ${}^\alpha\!f^\beta$
are well-defined and bijective.
\end{sublemma}

\begin{proof}[Proof of \ref{3bij}]
Let us first argue that if $Z \in \Integ^*_{(+,+)}$ satisfies
$Z \tequiv{(+,+)} (\alpha, \beta)$ then $Z$ can be $(+,+)$-contracted
to $(\alpha,\beta)$. 
Write $Z = (a,Z_1,b)$ where $a,b\in\Integ$ and $Z_1 \in \Integ^*$.
Then $Z_1 \sim \emptyseq$ by \ref{LessObvious} so, as is well-known,
$Z_1$ ``contracts'' to $\emptyseq$;
it follows that $Z$ can be $(+,+)$-contracted to $(\alpha',\beta')$, 
for suitable $\alpha',\beta'\in\Integ$.
Then $Z \tequiv{(+,+)} (\alpha', \beta')$ and hence
$ (\alpha', \emptyseq, \beta') \tequiv{(+,+)}  (\alpha, \emptyseq, \beta) $,
so \ref{ZeroProp} gives $(\alpha', \beta') =  (\alpha, \beta)$,
which proves our claim.
Similar remarks apply to $(+,-)$- and $(-,+)$-contraction, so
the alternative definitions
\begin{align*}
{}^\alpha\!E &= \setspec{ (X,Y) \in \Meul \times \Neul^* }
	{ (X,-1,Y) \tequiv{(+,-)} (\alpha) } \\
E^\alpha  &= \setspec{ (X,Y) \in \Neul^* \times \Meul^- }
	{ (X,-1,Y) \tequiv{(-,+)} (\alpha) } \\
{}^\alpha\!E^\beta &= \setspec{ (X,Y) \in \Meul \times \Meul^- }
	{ (X,-1,Y) \tequiv{(+,+)} (\alpha,\beta) }
\end{align*}
can be used if convenient.

By \ref{variation} and \ref{DetAdmiss}, four injective maps
$$
P \to \Neul^* \times \Neul^*,
\qquad
{}^\alpha\!P \to  \Meul \times \Neul^*,
\qquad
P^\alpha \to  \Neul^* \times \Meul^-
\quad \text{and} \quad
{}^\alpha\!P^\beta \to  \Meul \times \Meul^-
$$
are defined by stipulations ($1'$--$4'$).
The fact that the image of the second map is ${}^\alpha\!E$
can be derived from 3.23 of \DR, or from the reader's favorite
technique for handling linear chains. Then
it immediately follows that the third map has image $E^\alpha$
(simply because $Z_1 \tequiv{(-,+)} Z_2 \iff Z_1^- \tequiv{(+,-)} Z_2^-$).
Let us deduce that the fourth map has image ${}^\alpha\!E^\beta$
(the case of the first map is easier and is left to the reader).
We begin with:
\begin{equation}\label{ClaimInProof}
\begin{minipage}[t]{.9\linewidth}
\it 
Let $(X,Y) \in \Meul\times\Meul^-$ and write
$X = (a,X_1)$ and $Y = (Y_1,b)$, where $a,b\in\Integ$ and
$X_1,Y_1\in\Neul^*$. Then
$(X,Y) \in {}^\alpha\!E^\beta  \iff
\text{$(X,Y_1) \in {}^\alpha\!E$ and $(X_1,Y) \in E^\beta $}$.
\end{minipage}
\end{equation}
Indeed, if 
$(X,Y) \in {}^\alpha\!E^\beta$ then
$(X,-1,Y) \tequiv{(+,+)} (\alpha,\beta)$,
so 
$(X,-1, Y_1) \tequiv{(+,-)} (\alpha)$ and $(X_1,-1, Y) \tequiv{(-,+)} (\beta)$
by \ref{LessObvious}, so 
$(X,Y_1) \in {}^\alpha\!E$ and $(X_1,Y) \in E^\beta $.

Conversely, suppose that 
$(X,Y_1) \in {}^\alpha\!E$ and $(X_1,Y) \in E^\beta $.
Then
\begin{equation}\label{LuckyOne}
(X,-1, Y_1) \tequiv{(+,-)} (\alpha)
\end{equation}
and $(X_1,-1, Y) \tequiv{(-,+)} (\beta)$.
Applying \ref{LessObvious} to \eqref{LuckyOne} gives
$(X_1, -1, Y_1) \sim \emptyseq$,
so \ref{ForSuitable} gives
\begin{equation}\label{NextLu}
(X,-1,Y) \tequiv{(+,+)} (\alpha',\beta')
\end{equation}
for suitable $\alpha',\beta'\in\Integ$.
Applying \ref{LessObvious} to \eqref{NextLu} gives
$(X,-1, Y_1) \tequiv{(+,-)} (\alpha')$, so by \eqref{LuckyOne} we have
$(\alpha) \tequiv{(+,-)} (\alpha')$ and hence $\alpha=\alpha'$ by
comparing determinants.
Similarly, $\beta=\beta'$. 
So $(X,Y) \in {}^\alpha\!E^\beta$ and \eqref{ClaimInProof} is proved.

Let $(n,p,c) \in {}^\alpha\!P^\beta$ and define
$(X,Y) \in \Meul\times\Meul^-$ by condition ($4'$). 
We show that $(X,Y) \in {}^\alpha\!E^\beta$.
Write $X = (a,X_1)$ and $Y = (Y_1,b)$,
with $a,b\in\Integ$ and $X_1,Y_1\in\Neul^*$.
Note that ($4'$) implies
\begin{equation}\label{Comp1}
\det(Y^-) = (n-\beta+1)c - (c-p),
\qquad \text{where } 0 \le c-p < c.
\end{equation}
From ($4'$) we also have $\det(Y_1^-) = c$,
so applying \ref{remark*} to $Y^- = (b,Y_1^-)$ gives
\begin{equation}\label{Comp2}
\textstyle
\det(Y^-) = -bc - \det_1(Y_1^-),
\qquad \text{where } 0 \le \det_1(Y_1^-) < c
\end{equation}
(note that $0\le \det_1(Y_1^-) < \det(Y_1^-)= c$ follows
from $Y_1 \in \Neul^*$ and \ref{DetAdmiss});
comparing \eqref{Comp1} and \eqref{Comp2}, we obtain $\det_1(Y_1^-) = c-p$.
Together with ($4'$), this implies
\begin{align*}
\det(X) &= c - \alpha(nc+p) & \det(Y_1^-) &= c \\
\textstyle \det_1(X)&= nc+p & \textstyle \det_1(Y_1^-) &= c-p.
\end{align*}
Since $(n,p,c) \in {}^\alpha\!P$, these equations imply that
$(X,Y_1) = {}^\alpha\!f(n,p,c)$, so $(X,Y_1) \in {}^\alpha\!E$.

Define $n' = \big\lceil \frac{c}{nc+p} \big\rceil -1$, $c'=nc+p$
and $p' = c - n'c'$. Then one can verify that
$(n', p', c') \in P^\beta$ and that $f^\beta(n', p', c') = (X_1,Y)$,
so $(X_1,Y) \in E^\beta$.  By \eqref{ClaimInProof}, it follows that
$(X,Y) \in {}^\alpha\!E^\beta$. 

This shows that 
${}^\alpha\!f^\beta : {}^\alpha\!P^\beta \to {}^\alpha\!E^\beta$ is a
well-defined injective map.  Surjectivity is also proved by
using \eqref{ClaimInProof} and the fact that ${}^\alpha\!f$ and
$f^\beta$ are well-defined and bijective;
that part is easier, and is left to the reader.
\end{proof}

The next result lists all elements of $M^\oplus$,
where $M$ is any minimal element of $\Integ^*$.
Since $E$, ${}^\alpha\!E$, $E^\alpha$ and ${}^\alpha\!E^\beta$
are explicitely described by \ref{3bij}, this description
of $M^\oplus$ is explicit.

\begin{proposition}\label{MoplusExpl}
The elements of $\emptyseq^\oplus$ are:
\begin{enumerate}

\item[(i)] $(1)$

\item[(ii)] $(0,x)$ where $x \in \Integ \setminus \{ -1 \}$
\item[(iii)] $(x,0)$ where $x \in \Integ \setminus \{ -1 \}$

\item[(iv)] $(X,x,0,y,Y)$, where $(X,Y) \in E$
and $x,y \in \Integ \setminus \{ -1 \}$ satisfy $x+y=-1$.

\end{enumerate}

If $M = (m_1, \dots, m_k) \neq \emptyseq$ is a minimal element of $\Integ^*$,
the elements of $M^\oplus$ are:
\begin{enumerate}

\item \begin{enumerate}

	\item $(0,x,m_1,\dots,m_k)$, for all $x\in\Integ\setminus\{-1\}$

	\item the unique minimal sequence obtained by blowing-down
	$(0,-1,m_1,\dots,m_k)$

	\end{enumerate}

\item \begin{enumerate}

	\item $(m_1,\dots,m_k,x,0)$, for all $x\in\Integ\setminus\{-1\}$

	\item the unique minimal sequence obtained by blowing-down
	$(m_1,\dots,m_k,-1,0)$

	\end{enumerate}

\item For each $j \in \{1,\dots, k\}$,
\begin{enumerate}

	\item $(m_1,\dots, m_{j-1},x,0,y,m_{j+1},\dots, m_k)$,
	for all $x,y\in\Integ\setminus\{-1\}$ satisfying $x+y=m_j$

	\item the unique minimal sequence obtained by blowing-down
	$$
	(m_1,\dots, m_{j-1},-1,0,m_j+1,m_{j+1},\dots, m_k)
	$$

	\item the unique minimal sequence obtained by blowing-down
	$$
	(m_1,\dots, m_{j-1}, m_{j}+1,0,-1,m_{j+1},\dots, m_k)
	$$

	\end{enumerate}

\item \begin{enumerate}

	\item $(X,x,0,y,Y,m_2,\dots, m_k)$,
	for all $x,y\in\Integ\setminus\{-1\}$ satisfying $x+y=-1$ and all
	$(X,Y) \in E^{m_1}$

	\item $(m_1,\dots, m_{i-1},X,x,0,y,Y,m_{i+2},\dots, m_k)$,
	for all $x,y\in\Integ\setminus\{-1\}$ satisfying $x+y=-1$,
	all $(X,Y) \in {}^{m_i}\!E^{m_{i+1}}$
	and all $i$ such that $1\le i < k$

	\item $(m_1,\dots, m_{k-1},X,x,0,y,Y)$, for all
	$x,y\in\Integ\setminus\{-1\}$ satisfying $x+y=-1$ and all
	$(X,Y) \in {}^{m_k}\!E$.
	
	\end{enumerate}
	
\end{enumerate}
\end{proposition}

\begin{proof}
Follows from definitions  \ref{M+} (of $M^\oplus$)
and \ref{DefEalpha} (of $E$, ${}^{\alpha}\!E$, $E^{\alpha}$
and ${}^{\alpha}\!E^{\beta}$).
\end{proof}

\medskip
It is now clear that we can list the minimal elements 
of any class $\Ceul \in \SEC$ which is either a prime class
or the successor of a prime class.
Indeed, the problem is trivial if $\Ceul$ is a prime class,
and if $\Ceul$ is the successor of a prime class $\Ceul_1$ then
\ref{RecSol} gives $\min( \Ceul ) = M^\oplus$
where $M$ denotes the unique minimal element of $\Ceul_1$;
since the set $M^\oplus$ is described in \ref{MoplusExpl},
we obtain the desired list.
We give two concrete examples of this process:

\begin{example}\label{Min(1)}
Let $\Ceul$ denote the equivalence class of the sequence $(1)$.
Then $\Ceul = \Ceul_\emptyseq^\oplus$, where 
$\Ceul_\emptyseq$ is the equivalence class of the empty sequence
$\emptyseq$.
We have $\min\Ceul = \emptyseq^\oplus$ by \ref{RecSol} so,
by \ref{MoplusExpl}, the minimal elements of $\Ceul$ are:
\begin{itemize}

\item $(1)$

\item $(0,x)$ where $x \in \Integ \setminus \{ -1 \}$
\item $(x,0)$ where $x \in \Integ \setminus \{ -1 \}$

\item $(X,x,0,y,Y)$, where $(X,Y) \in E$
and $x,y \in \Integ \setminus \{ -1 \}$ satisfy $x+y=-1$.

\end{itemize}
See \ref{3bij} for an explicit description of $E$.
\end{example}

\begin{remark}
The result contained in \ref{Min(1)} first appeared in
\cite{Morrow:MinNorm} and was later reproved by several authors.
\end{remark}

\begin{example}\label{FirstEx}
Let $\Ceul = \Ceul_0^\oplus$, where $\Ceul_0 \in \SEC$
is the equivalence class of the sequence $(0)$.
By \ref{ListPrimes}, $(0)$ is the unique minimal element of $\Ceul_0$;
so \ref{RecSol} gives $\min\Ceul = (0)^\oplus$ and, by \ref{MoplusExpl},
the complete list of minimal elements of $\Ceul$ is:
\begin{itemize}

\item $(1,1)$

\item $(0,x,0)$ where $x \in \Integ \setminus \{ -1 \}$
\item $(x,0,-x)$ where $x \in \Integ \setminus \{ 1, -1 \}$

\item $(X,x,0,y,Y)$,
where $(X,Y) \in E^{0} \cup {}^0\!E$
and $x,y \in \Integ \setminus \{ -1 \}$ satisfy $x+y=-1$.

\end{itemize}
See \ref{3bij} for an explicit description of $E^{0}$ and ${}^0\!E$.
\end{example}

We leave it to the reader to reformulate the above facts
(\ref{RecSol}--\ref{FirstEx}) in terms of linear chains.

\section*{Geometric weighted graphs}

If $S$ is a smooth projective algebraic surface over an algebraically
closed field, and if $D$ is an SNC-divisor of $S$, then the pair
$(D,S)$ determines a weighted graph $\GG(D,S)$
called \textit{the dual graph of $D$ in $S$}
(see for instance \DR, \cite{Rus:formal} or \cite{Shastri:FinFundGp}).
A weighted graph $\GG$ is said to be \textit{geometric\/} if it is
isomorphic to $\GG(D,S)$ for some pair $(D,S)$, where we require that
every irreducible component of $D$ is a rational curve.
The purpose of this subsection is to point out:

\begin{proposition}\label{IdentGeom}
For a linear chain $\LL$, the following conditions are equivalent:
\begin{enumerate}

\item $\LL$ is geometric

\item $\hodge{\LL} \le 1$ or $\LL \sim [0,0,0]$

\item $\LL$ is equivalent to one of the following:
$$
\text{$[0]$, $[0,0,0]$, $[A]$ or $[0,0,A]$ (for some $A\in\Neul^*$)}
$$

\item Let $X\in\Integ^*$ be such that $\LL = [X]$; then the
equivalence class of $X$ is either a prime class or the successor of a
prime class.

\end{enumerate}
\end{proposition}

This fact is interesting in connection with the paragraph before \ref{Min(1)}.
We don't know a reference for \ref{IdentGeom},
but at least part of it is known.
Compare with 3.2.4 of \cite{Rus:formal}.
The proof of \ref{IdentGeom} requires the following fact:

\begin{nothing}\label{geometric}
Let $\GG$ be a geometric weighted graph.
\begin{enumerate}

\item $\hodge{\GG} \le 1$ or $\det(\GG)=0$.

\item If $\GG' \sim \GG$ then $\GG'$ is geometric.

\item Every induced subgraph of $\GG$ is geometric.

\item Let $\GG'$ be a weighted graph with the same underlying
graph as $\GG$ and such that $w(v, \GG') \le w(v, \GG)$
holds for every vertex $v$.  Then $\GG'$ is geometric.

\end{enumerate}
\end{nothing}

Note that a subgraph $G'$ of a graph $G$ is ``induced''
if every edge of $G$ which has its two endpoints in $G'$ is an edge of
$G'$.

Result \ref{geometric} is well known (the first assertion
is a consequence of the Hodge Index Theorem,
see for instance \cite{Rus:formal}; (2) and (3) are trivial and (4)
follows from (2) and (3)).
%
%

\begin{proof}[Proof of \ref{IdentGeom}]
It is clear that (3) is equivalent to (4);
we prove $(1) \Rightarrow (2) \Rightarrow (3) \Rightarrow (1)$.
Suppose that $\LL$ is geometric and that $\det(\LL)=0$.
By \ref{uniqueness}, $\LL\sim [0^{2n+1}]$ for some $n\in\Nat$;
by parts (2) and (3) of \ref{geometric}, it follows that
$[0^{2n+1}]$ is geometric and then that $[0^{2n}]$ is geometric.
We have $\det [0^{2n}] = (-1)^n$ and $\hodge{[0^{2n}]}=n$
by \ref{0iA} and \ref{BasicHodge}, so $n\le1$ by part~(1)
of \ref{geometric}. Thus:
\begin{equation*}
\textit{If $\LL$ is geometric and $\det\LL=0$ then
$\LL$ is equivalent to $[0]$ or $[0,0,0]$.}
\end{equation*}
The fact that (1) implies (2) follows from this and part~(1)
of \ref{geometric}.

Consider a canonical linear chain $[0^{r},A]$ equivalent to $\LL$
(with $r\in\Nat$, $A\in\Neul^*$ and $r$ is even if $A\neq\emptyseq$).
If (2) holds then $r<4$,
so $[0^{r},A]$ is one of the chains displayed in assertion~(3).
So (2) implies (3).

To show that (3) implies (1), we have to check that each of
$[0]$, $[0,0,0]$, $[A]$, $[0,0,A]$ (where $A\in\Neul^*$) is geometric;
by part~(3) of \ref{geometric}, it suffices to prove that 
$[0,0,0]$ and $[0,0,A]$ are geometric, where we may assume that
$A\neq\emptyseq$.
Considering a pair of lines in $\proj^2$ shows that
$[1,1]$ is geometric; so $[0,0,0] \sim [1,1]$ is geometric.
Let $n\ge1$ be such that $A=(a_1,\dots,a_n)$.
If $n=1$ then $[0,0,A]$ is geometric by applying part~(4) of
\ref{geometric} to $[0,0,A]$ and $[0,0,0]$;
if $n>1$ then $[0,0,-1,-2^{n-2},-1] \sim [0,0,0]$ is geometric and,
by part~(4) of \ref{geometric} applied to $[0,0,A]$ and $[0,0,-1,-2^{n-2},-1]$,
$[0,0,A]$ is geometric.
\end{proof}

\section{Pseudo-minimal forests with a given skeleton}
\label{Section:RegForSkel}

Throughout this section we fix a skeleton $S$ (see \ref{DefSkel}).

\begin{notation}\label{DefRFSK}
Let $\RFSK(S)$ denote the set of pseudo-minimal forests $\GG$ satisfying:
$$
\textit{There exists at least one skeletal map $\sigma : S \skm \GG$.}
$$
\end{notation}

\begin{problem}\label{ProbS}
Classify the elements of $\RFSK(S)$
up to equivalence of weighted graphs.
\end{problem}

This is the fundamental problem that has to be solved since,
by \ref{RegStrict} and \ref{AllReg}, a solution to Problem~\ref{ProbS}
for every $S$ includes a classification of all minimal weighted forests
(and hence of all weighted forests).

\begin{definition}\label{DefFor}
Let $\For(S)$ be the set of ordered pairs $(\sigma,\GG)$ where
$\GG$ is a pseudo-minimal forest and $\sigma : S \skm \GG$ is a skeletal map.
Two elements $(\sigma,\GG),(\sigma',\GG')$ of $\For(S)$ are
\textit{isomorphic\/} if there exists an isomorphism $f : \GG \to \GG'$
of weighted graphs such that $f\circ\sigma = \sigma'$.
Let $p_2 : \For(S) \to \RFSK(S)$ be the surjection
given by $p_2(\sigma,\GG) = \GG$.
\end{definition}

\begin{definition}\label{EquFor}
Let $(\sigma,\GG), (\sigma',\GG') \in \For(S)$.
If $\GG'$ is a strict blowing-down of $\GG$,
and if the blowing-down map $\pi : \GG \skm \GG'$
(defined in \ref{DefBlowDownMap}) satisfies 
$\pi \circ \sigma = \sigma'$, then we say that 
$(\sigma',\GG')$ is a \textit{blowing-down\/} of $(\sigma,\GG)$
and that $(\sigma,\GG)$ is a \textit{blowing-up\/} of
$(\sigma',\GG')$.
Two elements of $\For(S)$ are \textit{equivalent\/}
(notation: ``$\sim$'') if one can be obtained from the other via
a finite sequence of blowings-up and blowings-down.
\end{definition}

\begin{problem}\label{ProbForS}
Classify the elements of $\For(S)$ up to equivalence $( \sim )$.
\end{problem}

\medskip
We shall first solve Problem~\ref{ProbForS} and then
derive a solution to Problem~\ref{ProbS}.
The precise relation between the two problems will be described after
the solution to Problem~\ref{ProbForS};
we will see that Problem~\ref{ProbS} is essentially 
Problem~\ref{ProbForS} modulo automorphisms of $S$.
For instance, if $S$ is 
\begin{picture}(14,3)(-2,-1)
\put(0,0){\circle*{1}}
\put(10,0){\circle*{1}}
\put(0,0){\line(1,0){10}}
\end{picture}
then Problem~\ref{ProbS} asks for the classification of all linear
chains not equivalent to $\Empty$
(solved in section~\ref{Sec:ClassLinChains})
and Problem~\ref{ProbForS} can be seen to be equivalent to the
classification of sequences not equivalent to $\emptyseq$ (solved in
section~\ref{SectionSequences}).

Paragraphs \ref{DefSWw}--\ref{MainForS} solve Problem~\ref{ProbForS}.
The machinery developed for solving the problem is, we think,
as meaningful as the final answer, stated in \ref{MainForS}.

\medskip

\begin{definition}\label{DefSWw}
An \textit{edge map\/} for $S$ is a set map
$W: P(S) \to \Integ^*$ satisfying the two conditions:
\begin{enumerate}

\item $\forall_{ \gamma \in P(S) }\ \ W( \gamma^- ) = W(\gamma)^-$

\item $\forall_{ \gamma \in P(S) }\ \ \ 
W( \gamma ) \sim \emptyseq \implies \text{$\gamma$ is of type $(+,+)$}$.

\end{enumerate}
Two edge maps $W$ and $W'$ for $S$ are \textit{equivalent\/} if
$ \forall_{\gamma \in P(S)}\ \ W(\gamma) \sim W'(\gamma)$
(equivalences in $\Integ^*$).
The symbol $\Omega(S)$ denotes the set of equivalence classes of edge
maps for $S$.
If $\omega\in\Omega(S)$ and $\gamma \in P(S)$, we define
$$
\omega(\gamma) = \setspec{ W(\gamma) }{ W\in\omega },
$$
so $\omega(\gamma) \subset \Integ^*$ is an equivalence class of sequences of
integers and it  makes sense to speak of the determinant of $\omega(\gamma)$
(see \ref{TwoInvars}).
\end{definition}

Clearly, if $(\sigma,\GG)\in\For(S)$ then the composite
$P(S) \xrightarrow{\ \vec{\sigma}\ } P(\GG) \xrightarrow{W_\GG} \Integ^*$
is an edge map for $S$
(see \ref{DefSkelMap} for $\vec{\sigma}$ and \ref{DefWgamma} for $W_\GG$).


\begin{definition}
Given $\omega \in \Omega(S)$, let $\For(S,\omega)$ be the set of
pairs $(\sigma,\GG)\in\For(S)$ satisfying:
\begin{equation*}
\textit{The composite
$P(S) \xrightarrow{\ \vec{\sigma}\ } P(\GG) \xrightarrow{W_\GG} \Integ^*$
is an element of $\omega$.}
\end{equation*}
Note that $\setspec{ \For(S,\omega) }{ \omega \in \Omega(S) }$
is a partition of $\For(S)$.
\end{definition}

\begin{lemma}\label{SameOmega}
Let $(\sigma,\GG), (\sigma',\GG')  \in \For(S)$ and $\omega \in \Omega(S)$.
If $(\sigma,\GG) \sim (\sigma',\GG')$ and 
$(\sigma,\GG) \in \For(S,\omega)$,
then $(\sigma',\GG') \in \For(S,\omega)$.
\end{lemma}

\begin{proof}
If $(\sigma',\GG')$ is a blowing-down of $(\sigma,\GG)$,
and if $\pi : \GG \skm \GG'$ is the blowing-down map, 
then for each $\gamma \in P(\GG)$ the sequences
$W_\GG(\gamma)$ and $W_{\GG'} \big( \vec\pi(\gamma) \big)$ 
are equivalent; thus
the composites
$P(S) \xrightarrow{\ \vec{\sigma}\ } P(\GG) \xrightarrow{W_\GG} \Integ^*$
and
$P(S) \xrightarrow{\ \vec{\sigma'}\ } P(\GG') \xrightarrow{W_{\GG'}} \Integ^*$
are equivalent edge maps for $S$ and, consequently,
$(\sigma,\GG) \in \For(S,\omega) \Leftrightarrow
(\sigma',\GG') \in \For(S,\omega)$.
The desired result follows.
\end{proof}

By \ref{SameOmega}, Problem~\ref{ProbForS} reduces to classifying
the elements of $\For(S,\omega)$ for each $\omega\in\Omega(S)$.

\medskip
\noindent {\bf Set-up.} Recall that $S$ was fixed at the beginning of
the section.
From here to \ref{Equ=action}, we fix $\omega \in \Omega(S)$ and classify
elements of $\For(S,\omega)$.

\medskip

\begin{definition}\label{DefMapT}
If $(\sigma,\GG) \in \For(S,\omega)$,
define $T(\sigma,\GG) = (W,w)$ where
$W: P(S) \to \Integ^*$ is the composite
$P(S) \xrightarrow{\ \vec{\sigma}\ } P(\GG) \xrightarrow{W_\GG} \Integ^*$
and $w : \Vtx_{>2}(S) \to \Integ$ is the map given by
$w(v) = w( \sigma(v), \GG)$.
It is clear that 
$$
T : \For(S,\omega) \to \omega \times \Integ^{\Vtx_{>2}(S)}
$$
is surjective and that the inverse image of any element of
$\omega \times \Integ^{\Vtx_{>2}(S)}$
is an isomorphism class of pairs $(\sigma,\GG)$ (isomorphism is
defined in \ref{DefFor}).
\end{definition}

\section*{The transplant operation}

\begin{nothing*}
Given $\gamma =(v_0,v_1) \in P(S)$ of type $\tau$,
define a map
$\Xeul_\gamma : \omega\times\Integ^{\Vtx_{>2}(S)} \to \Integ^*_\tau$
by
$$
\Xeul_\gamma (W,w) =
\begin{cases}
W(\gamma) & \text{if $\tau=(-,-)$} \\
\big( w(v_0), W(\gamma) \big)
& \text{if $\tau=(+,-)$} \\
\big( W(\gamma), w(v_1) \big)
& \text{if $\tau=(-,+)$} \\
\big( w(v_0), W(\gamma), w(v_1) \big)
& \text{if $\tau=(+,+)$}.
\end{cases}
$$
\end{nothing*}

\begin{definition}\label{DefTransplant}
Let $(W,w) \in \omega\times\Integ^{\Vtx_{>2}(S)}$,
let $\gamma = (v_0, v_1) \in P(S)$ be of type $\tau$
and let $\Yeul \in \Integ^*_\tau$
be such that $\Xeul_\gamma(W,w) \tequiv{\tau} \Yeul$.
Then a unique pair $(W',w') \in \omega\times\Integ^{\Vtx_{>2}(S)}$ is
determined by
\begin{itemize}

\item $\Xeul_\gamma(W',w') = \Yeul$

\item $W'$ agrees with $W$ on $P(S) \setminus \{ \gamma, \gamma^- \}$

\item $w'$ agrees with $w$ on $\Vtx_{>2}(S) \setminus  \{ v_0, v_1 \}$.

\end{itemize}
We say that $(W',w')$ is obtained by
\textit{transplanting $(\gamma, \Yeul)$ into $(W,w)$}
and write
$$
(W',w') = \transp(\gamma, \Yeul; W,w).
$$
If this is the case, and if 
$(\sigma,\GG), (\sigma',\GG') \in \For(S,\omega)$
satisfy $T(\sigma,\GG) = (W,w)$ and $T(\sigma',\GG')= (W',w')$,
we also say that $(\sigma',\GG')$ is obtained by
\textit{transplanting $(\gamma, \Yeul)$ into $(\sigma,\GG)$.}
\end{definition}

\begin{sublemma}\label{Transp>Equiv}
Let $(\sigma,\GG), (\sigma',\GG') \in \For(S,\omega)$
and suppose that $(\sigma',\GG')$ is obtained by
transplanting some pair $(\gamma, \Yeul)$ into $(\sigma,\GG)$.
Then $(\sigma,\GG) \sim (\sigma',\GG')$.
\end{sublemma}

\begin{proof}
The weighted graph $\GG'$ is obtained from $\GG$ by performing the
operation described in \ref{ApplTau}; more precisely, the operation
is performed on $(u_1,\dots,u_m) = \vec{\sigma}(\gamma) \in P(\GG)$.
So \ref{ApplTau} gives $\GG \sim \GG'$ and it is easy to see that
$(\sigma,\GG) \sim (\sigma',\GG')$.
\end{proof}

\section*{Congruence}

See \ref{DefSWw} for the definition of $\omega(\gamma)$.

\begin{definition}\label{DefN}
\label{DefSkelTriple}\label{Correspondence}
\begin{enumerate}

\item Let $(S,\omega)^\sharp$ denote the forest (not weighted) whose vertex
set is $\Vtx_{>2}(S)$ and whose edges are the pairs
$\{ u,v \}$ of vertices satisfying:
$$
\text{$\gamma=(u,v)$ belongs to $P(S)$ and $\omega (\gamma)$
has determinant zero.}
$$

\item A vertex $u$ of $(S,\omega)^\sharp$ is \textit{special\/} if it satisfies:
\begin{enumerate}

\item[] there exists $v\in\Vtx(S)$ such that $\gamma=(u,v)$ belongs
to $P(S)$, $\gamma$ is of type $(+,-)$ and 
$\omega (\gamma)$ has determinant zero.

\end{enumerate}

\item Let $Z(S,\omega)$ be the set of all maps $z: \Vtx_{>2}(S)\to\Integ$
satisfying:
\textit{For each connected component $C$ of $(S,\omega)^\sharp$
which contains no special vertex,}
$$
\sum_{v\in\Vtx(C)}  z(v) = 0.
$$
We give generators for the submodule $Z(S,\omega)$
of the $\Integ$-module $\Integ^{\Vtx_{>2}(S)}$.
For each $\epsilon=(u,u') \in P(S)$ of type $(+,+)$ and
such that $\omega(\epsilon)$ has determinant zero,
define $z_\epsilon \in Z(S,\omega)$ by 
$z_\epsilon (u)=1$,
$z_\epsilon (u')=-1$ and 
$z_\epsilon (v)=0$ for all $v\in\Vtx_{>2}(S) \setminus \{ u,u' \}$;
let $Z_e(S,\omega)$ be the set of these $z_\epsilon$.
For each special vertex $u$ of $(S,\omega)^\sharp$,
define $z_u \in Z(S,\omega)$
by $z_u(u)=1$ and $z_u (v)=0$ for all $v\in\Vtx_{>2}(S) \setminus \{ u \}$;
let $Z_s(S,\omega)$ be the set of these $z_u$.
Then the reader may verify that
\begin{equation*}
\textit{$Z(S,\omega)$ is generated as a $\Integ$-module
by $Z_e(S,\omega) \cup Z_s(S,\omega)$.}
\end{equation*}

\item Define the $\Integ$-module
$ N(S,\omega) = \Integ^{\Vtx_{>2}(S)} /  Z(S,\omega) $.

\item Let $\bar T : \For(S,\omega) \to \omega \times N(S,\omega)$
be the composite:
$$
\begin{array}{rcl}
\For(S,\omega) \xrightarrow{\ T\ } \omega \times \Integ^{\Vtx_{>2}(S)}
& \longrightarrow & \omega \times N(S,\omega)\\
\big( W,w \big) & \longmapsto & \big( W,\pi(w) \big)
\end{array}
$$
where $\pi :  \Integ^{\Vtx_{>2}(S)} \to  N(S,\omega)$ is the canonical
epimorphism and $T$ is defined in \ref{DefMapT}.
Note that $\bar T$ is surjective.

\item Let $(\sigma,\GG), (\sigma',\GG') \in \For(S,\omega)$.
If $\bar T(\sigma,\GG) = \bar T(\sigma',\GG')$ then
we write $(\sigma,\GG) \equiv (\sigma',\GG')$ and 
say that $(\sigma,\GG), (\sigma',\GG')$ are \textit{congruent}.

\end{enumerate}
\end{definition}

\begin{lemma}\label{Cong=>Equiv}
If $(\sigma,\GG), (\sigma',\GG') \in \For(S,\omega)$
and $(\sigma,\GG) \equiv (\sigma',\GG')$,
then $(\sigma,\GG) \sim (\sigma',\GG')$.
\end{lemma}

\begin{proof}
Write $(W,w) = T(\sigma,\GG)$ and $(W,w') = T(\sigma',\GG')$. 
Since $(\sigma,\GG) \equiv (\sigma',\GG')$,
we have $w' - w \in Z(S,\omega)$ and consequently there exists a sequence
$(w_0, \dots, w_n)$ in $\Integ^{\Vtx_{>2}(S)}$ satisfying $w_0=w$,
$w_n=w'$ and, for all $i>0$,
$\pm(w_{i+1}-w_i) \in Z_e(S,\omega) \cup Z_s(S,\omega)$
(notation as in part (3) of \ref{DefN}).
Define $(\sigma_0,\Gn0) = (\sigma,\GG)$,
$(\sigma_n,\Gn n) = (\sigma',\GG')$ and, for each $i$ such that
$0<i<n$, choose $(\sigma_i,\Gn i) \in \For(S,\omega)$ such that
$T( \sigma_i,\Gn i) = (W,w_i)$.
Note that $( \sigma_i,\Gn i) \equiv ( \sigma_{i+1},\Gn{i+1})$
for every $i$;
so it suffices to prove that 
$(\sigma,\GG) \sim (\sigma',\GG')$ under the assumption that
$$
w' - w \in Z_e(S,\omega) \cup Z_s(S,\omega).
$$
Assume that $w' - w = z_\epsilon \in Z_e(S,\omega)$
and write $\epsilon = (v_0, v_1)$;
recall that $\epsilon \in P(S)$ is of type $(+,+)$ and that
$\det W(\epsilon)=0$.
From $w' = w + z_\epsilon$ and part (2a) of \ref{ZeroProp} we get
$$
\big( w'(v_0), W(\epsilon), w'(v_1) \big)
=
\big( w(v_0)+1, W(\epsilon), w(v_1)-1 \big)
\tequiv{(+,+)}
\big( w(v_0), W(\epsilon), w(v_1) \big),
$$
which we rewrite as
$
\Xeul_\epsilon(W,w') \tequiv{(+,+)} \Xeul_\epsilon(W,w).
$
So it makes sense to
transplant $(\epsilon,\Yeul)$ into $(\sigma,\GG)$,
where $\Yeul =\Xeul_\epsilon(W,w')$, and clearly
$\transp(\epsilon, \Yeul; W,w) = (W,w')$.
So $(\sigma',\GG')$ can be obtained by transplanting
$(\epsilon,\Yeul)$ into $(\sigma,\GG)$ and, by \ref{Transp>Equiv},
$(\sigma,\GG) \sim (\sigma',\GG')$.

The case where $w' - w \in Z_s(S,\omega)$ is proved in a similar way,
except that part (2b) of \ref{ZeroProp} is used in place of (2a).
\end{proof}

\section*{Transplant as an action}

\begin{nothing*}
For each type $\tau$ (\ref{DefBup}),
define a restriction map $\Integ^*_\tau \to \Integ^*$, $X\mapsto X|_\tau$,
by declaring that if
$X = (x_1,\dots,x_n)\in\Integ^*_\tau$ then
$$
X|_\tau = \begin{cases}
X, & \text{if } \tau = (-,-) \\
(x_2, \dots, x_n), & \text{if } \tau = (+,-) \\
(x_1, \dots, x_{n-1}), & \text{if } \tau = (-,+) \\
(x_2, \dots, x_{n-1}), & \text{if } \tau = (+,+).
\end{cases}
$$
\end{nothing*}

\begin{definition}\label{DefR}
Let $\Reul(S,\omega)$ be the set of ordered pairs $\binom{Y}{\gamma}$
(written vertically) such that $\gamma \in P(S)$ and $Y \in \omega(\gamma)$.
\end{definition}

\begin{proposition}\label{PropTheta}
Let $y = \binom{Y}{\gamma} \in \Reul(S,\omega)$.
There is a unique set map
$$
\Theta_y \ :\ 
\omega \times N(S,\omega) \longrightarrow \omega \times N(S,\omega)
$$
which satisfies the following condition:

Let $(W,w) \in \omega \times \Integ^{\Vtx_{>2}(S)}$
and $\Yeul \in \Integ_\tau^*$ (where $\tau$ is the type of $\gamma$)
be such that $\Yeul |_\tau = Y$ and $\Yeul \tequiv{\tau} \Xeul_\gamma(W,w)$;
let $(W',w') = \transp(\gamma,\Yeul; W,w)$; then
$ \Theta_y(W,\pi(w)) = (W', \pi(w'))$, where
$\pi : \Integ^{\Vtx_{>2}(S)} \to N(S,\omega)$ is the canonical epimorphism.
\end{proposition}

\begin{proof}
Given $(W,w) \in \omega \times  \Integ^{\Vtx_{>2}(S)}$,
choose $\Yeul \in \Integ^*_\tau$ satisfying
\begin{equation}\label{CondYeul}
\Yeul|_\tau = Y \quad \text{and} \quad \Yeul \tequiv{\tau} \Xeul_\gamma(W,w)
\end{equation}
(by \ref{ForSuitable}, there exists at least one such $\Yeul$);
then let $(W',w') = \transp(\gamma, \Yeul; W,w)$ and define
$\theta_y(W,w) = (W', \pi(w'))$.
The proof of \ref{PropTheta} consists in showing that
$$
\theta_y : \omega \times  \Integ^{\Vtx_{>2}(S)}
\to \omega \times N(S,\omega) 
$$
is a well-defined map and satisfies
\begin{equation}\label{wIndep}
\text{$\theta_y(W,w) = \theta_y(W,w+z)$, \quad
for all $W \in \omega$, $w \in \Integ^{\Vtx_{>2}(S)}$ and $z\in Z(S,\omega)$.}
\end{equation}

Let $(W,w) \in \omega \times  \Integ^{\Vtx_{>2}(S)}$.
For $i=1,2$, consider $\Yeul_i \in \Integ^*_\tau$
such that $\Yeul_i|_\tau = Y$ and $\Yeul_i \tequiv{\tau} \Xeul_\gamma(W,w)$;
let $(W_i',w_i') = \transp(\gamma, \Yeul_i; W,w)$. 
From $\Yeul_1|_\tau = \Yeul_2|_\tau$, we obtain $W_1'=W_2'$.
To prove that $\pi(w_1')=\pi(w_2')$, we may assume that
$\Yeul_1 \neq \Yeul_2$;
then \ref{ZeroProp} implies that $\det(Y)=0$
and also describes how $\Yeul_1$ differs from $\Yeul_2$;
that description implies that
$w_1' - w_2'$ is a multiple of an element of 
$Z_e(S,\omega) \cup Z_s(S,\omega)$, so in particular
$w_1' - w_2' \in Z(S,\omega)$ and $\pi(w_1')=\pi(w_2')$.
Thus $(W_1',\pi(w_1')) = (W_2',\pi(w_2'))$, which shows that
$\theta_y$ is a well-defined map.

Let $W\in\omega$,  $w \in \Integ^{\Vtx_{>2}(S)}$ and
$z\in Z(S,\omega)$.
Choose $\Yeul \in \Integ^*_\tau$ 
such that $\Yeul |_\tau = Y$ and $\Yeul \tequiv{\tau} \Xeul_\gamma(W,w)$;
let $(W',w') = \transp(\gamma, \Yeul; W,w)$.
By \ref{AllIJ}, there exits $\Yeul^* \in \Integ^*_\tau$ 
such that 
$$
\Yeul^* |_\tau = Y, \quad
\Yeul^* \tequiv{\tau} \Xeul_\gamma(W,w+z) \quad \text{and} \quad
\transp(\gamma, \Yeul^*; W,w+z) = (W',w'+z)
$$
(for instance, if $\gamma = (v_0,v_1)$ is of type $(+,-)$ then
$\Xeul_\gamma(W,w+z)$ is obtained from $\Xeul_\gamma(W,w)$
by adding $z(v_0)$ to the leftmost term of the sequence;
then let $\Yeul^*$ be the sequence obtained by adding $z(v_0)$
to the leftmost term of $\Yeul$).
By definition of $\theta_y$, we have
$ \theta_y(W,w+z) = \big( W', \pi(w'+z) \big) = \big( W', \pi(w') \big)
= \theta_y(W,w)$, which proves 
that $\theta_y$ satisfies \eqref{wIndep}.
\end{proof}

\begin{definition}\label{DefLeftAction}
\begin{enumerate}

\item The symbol $\Reul^*(S,\omega)$ denotes the free monoid on the
set $\Reul(S,\omega)$.  That is, the elements of 
$\Reul^*(S,\omega)$ are the words $y_1\cdots y_n$ where
$n\in\Nat$ and $y_1,\dots,y_n\in \Reul(S,\omega)$, and 
the operation is concatenation.

\item
If $y \in \Reul(S,\omega)$ and
$(W,\eta) \in \omega \times N(S,\omega)$,
define
$y(W,\eta) = \Theta_y (W,\eta)$.
This extends uniquely to a left-action of the monoid
$\Reul^*(S,\omega)$ on the set $\omega \times N(S,\omega)$.

\end{enumerate}
\end{definition}

Because of \ref{Equ=action} below,
this left-action of $\Reul^*(S,\omega)$ on $\omega \times N(S,\omega)$
plays a central role in the classification.
We now investigate the properties of that action.
The first fact is easily verified:

\begin{nothing}\label{BasicPptiesAction}
Let $(W,\eta) \in \omega \times N(S,\omega)$.
\begin{enumerate}

\item If $\gamma \in P(S)$ then the element
$y=\binom{W(\gamma)}{\gamma}$ of $\Reul(S,\omega)$ satisfies
$y(W,\eta) = (W,\eta)$.

\item Let $y=\binom{Y}{\gamma} \in \Reul(S,\omega)$ and let
$y^-=\binom{Y^-}{\gamma^-}$.  Then $y^-$ belongs to $\Reul(S,\omega)$
and satisfies $y^- (W,\eta) = y (W,\eta)$.

\item If $y = \binom Y{\gamma} \in \Reul(S,\omega)$ and
$(W',\eta') = y(W,\eta)$, then $W'$ is determined by:
$$
W'(\gamma) = Y,
\quad W'(\gamma^-) = Y^-
\quad \text{and} \quad
W'(\gamma') = W(\gamma')
\text{ for all } \gamma' \in P(S) \setminus \{ \gamma, \gamma^- \}.
$$

\end{enumerate}
\end{nothing}

\begin{definition}
\begin{enumerate}


\item Two elements $\binom{Y_1}{\gamma_1}$ and $\binom{Y_2}{\gamma_2}$ of 
$\Reul(S,\omega)$ are said to be \textit{disjoint\/} if
$\{ \gamma_1,\gamma_1^- \} \cap \{ \gamma_2,\gamma_2^- \} 
= \emptyset$.

\item A word $y \in \Reul^*(S,\omega)$ is \textit{self-disjoint\/} if
the unique $y_1,\cdots, y_n \in \Reul(S,\omega)$ satisfying
$y=y_1\cdots y_n$ are pairwise disjoint
(i.e., $y_i,y_j$ are disjoint whenever $i\neq j$).
In particular, the empty word and all elements of
$\Reul(S,\omega)$  are self-disjoint.

\end{enumerate}
\end{definition}

\begin{lemma}\label{AdvPptiesAction}
Let $\xi \in \omega \times N(S,\omega)$.
\begin{enumerate}

\item If $y_1,y_2\in \Reul(S,\omega)$ are not disjoint,
then $y_2y_1 \xi = y_2 \xi$.

\item If $y_1,y_2 \in \Reul(S,\omega)$ are disjoint,
then $y_2y_1 \xi = y_1y_2 \xi$.

\item Given any $y \in \Reul^*(S,\omega)$, there exists a
self-disjoint word $y' \in \Reul^*(S,\omega)$ such that $y \xi = y' \xi$.

\end{enumerate}
\end{lemma}

\begin{proof}
The last assertion easily follows from the first two,
so we prove only (1) and (2).
Let $\pi : \Integ^{\Vtx_{>2}(S)} \to N(S,\omega)$ be the canonical epimorphism
and consider an arbitrary element $(W,\pi(w))$ of $\omega \times N(S,\omega)$
(where $w \in \Integ^{\Vtx_{>2}(S)}$).
We have to prove:
\begin{equation}\label{ThetaCondition}
(\Theta_{y_2} \circ \Theta_{y_1}) (W,\pi(w))
= \begin{cases}
\Theta_{y_2}(W,\pi(w)) & \text{if $y_1, y_2$ are not disjoint,} \\
(\Theta_{y_1} \circ \Theta_{y_2})(W,\pi(w))
&  \text{if $y_1, y_2$ are disjoint.}
\end{cases}
\end{equation}

For each $i \in \{ 1, 2 \}$,
write $y_i = \binom{Y_i}{\gamma_i}$ and $\gamma_i = (v_0^i, v_1^i )$,
and let $\tau_i$ be the type of $\gamma_i$.
Choose $\Yeul_i \in \Integ_{\tau_i}^*$
satisfying $\Yeul_i |_{\tau_i} = Y_i$
and $\Yeul_i \tequiv{\tau_i} \Xeul_{\gamma_i}(W,w)$;
let 
\begin{equation}\label{WiTransp}
(W_i,w_i) = \transp(\gamma_i,\Yeul_i; W,w);
\end{equation}
then $\Theta_{y_i} (W,\pi(w) ) = (W_i, \pi(w_i))$.
We record the following consequences of \eqref{WiTransp}
(cf.\ \ref{DefTransplant}):
\begin{gather}
\label{iTransp1}
\Xeul_{\gamma_i}(W_{i},w_{i}) = \Yeul_i \\
\label{iTransp2}
\text{$W_{i}(\gamma') = W(\gamma')$
for all $\gamma' \in P(S)\setminus\{\gamma_i, \gamma_i^- \}$}\\
\label{iTransp3}
\text{$w_{i}(v) = w(v)$
for all $v\in \Vtx_{>2}(S)\setminus \{ v_0^i, v_1^i \}$.}
\end{gather}

For each $(i,j) \in \{ (1,2), (2,1) \}$,
choose $\Yeul_{ij} \in \Integ_{\tau_j}^*$
satisfying $\Yeul_{ij} |_{\tau_j} = Y_j$
and $\Yeul_{ij} \tequiv{\tau_j} \Xeul_{\gamma_j}(W_i,w_i)$;
let $(W_{ij},w_{ij}) = \transp(\gamma_j,\Yeul_{ij}; W_i,w_i)$; then
$\Theta_{y_j} (W_i,\pi(w_i) ) = (W_{ij}, \pi(w_{ij}))$
or equivalently
$$
(\Theta_{y_j} \circ \Theta_{y_i}) (W,\pi(w)) =  (W_{ij}, \pi(w_{ij})).
$$

\smallskip
Suppose that $y_1,y_2$ are not disjoint. 
Then $\gamma_2 \in \{ \gamma_1, \gamma_1^- \}$ and,
by the second part of \ref{BasicPptiesAction},
we may assume that $\gamma_2 = \gamma_1$;
in fact we write $\gamma_1 = \gamma_2 = \gamma = (v_0,v_1)$
and $\tau_1 = \tau_2 = \tau$.
Note that 
$$
\Yeul_{2} |_\tau = Y_2
\quad \text{and} \quad
\Yeul_2
\tequiv{\tau}
\Xeul_{\gamma}(W,w)
\tequiv{\tau}
\Yeul_1
\overset{\text{\tiny \eqref{iTransp1}}}{=}
\Xeul_{\gamma}(W_1,w_1)
$$
so, when we choose $\Yeul_{12}$, we may set $\Yeul_{12}=\Yeul_2$.
If we do this, then
\begin{equation}\label{W12W2}
(W_{12}, w_{12})
= \transp(\gamma, \Yeul_2; W_1,w_1)
= \transp(\gamma, \Yeul_2; W,w)
= (W_{2}, w_{2})
\end{equation}
where the middle equality is a consequence of
\eqref{iTransp2} and \eqref{iTransp3}.
Now \eqref{W12W2} implies that
$$
(\Theta_{y_2} \circ \Theta_{y_1}) (W, \pi(w))
=  (W_{12}, \pi(w_{12}))
= (W_{2}, \pi(w_{2}))
= \Theta_{y_2} (W, \pi(w)),
$$
which proves the first part of \eqref{ThetaCondition}.

\smallskip
Suppose that $y_1,y_2$ are disjoint. 
Then $\{ v_0^1, v_1^1 \} \neq \{ v_0^2, v_1^2 \}$, but these sets may or
may not be disjoint.  We prove the second part of \eqref{ThetaCondition}
in the case
where these sets are not disjoint, as this is the most delicate of the
two cases.  In view of part~(2) of \ref{BasicPptiesAction}, we may arrange:
\begin{equation*}
\text{$\gamma_1 = (v_0^1, v_*)$ and $\gamma_2 = (v_0^2, v_*)$,
where $v_0^1,  v_0^2,  v_*$ are distinct.}
\end{equation*}
Note that $\deg(v_*,S)>1$,
so $\deg(v_*,S)>2$ and $\tau_1,\tau_2 \in \{ (-,+), (+,+) \}$.
We adopt the following notation: If $\tau \in \{ (-,+), (+,+) \}$,
$\Yeul \in \Integ_\tau^*$ and $n\in\Integ$, then
$\Yeul^n \in \Integ_\tau^*$ is the sequence obtained by adding $n$ to
the rightmost term of $\Yeul$.

Let $(i,j) \in \{ (1,2), (2,1) \}$.
With $n_i = w_i(v_*) - w(v_*)$, we have
\begin{equation}\label{AllowsSet}
\Yeul_j^{n_i} |_{\tau_j} = Y_j
\quad \text{and} \quad
\Yeul_j^{n_i}
\tequiv{\tau_j}
\Xeul_{\gamma_j}(W,w)^{n_i}
=
\Xeul_{\gamma_j}(W_i,w_i)
\end{equation}
where the $\tau_j$-equivalence follows from \ref{AllIJ}
and where the last equality is verified directly:
If $\tau_j = (-,+)$ (resp.~$(+,+)$) then
\begin{align*}
\Xeul_{\gamma_j}(W,w)^{n_i}
&= (W(\gamma_j), w(v_*) + n_i)
& \text{(resp.\ } & ( w(v_0^j), W(\gamma_j), w(v_*) + n_i )) \\
\Xeul_{\gamma_j}(W_i,w_i)
&= (W_i(\gamma_j), w_i(v_*))
& \text{(resp.\ } & ( w_i(v_0^j), W_i(\gamma_j), w_i(v_*) ))
\end{align*}
and these two sequences are equal by \eqref{iTransp2} and \eqref{iTransp3}.
Statement \eqref{AllowsSet} allows us to set $\Yeul_{ij} = \Yeul_j^{n_i}$. 
Then $(W_{ij},w_{ij}) = \transp(\gamma_j,\Yeul_j^{n_i}; W_i,w_i)$
and consequently
\begin{gather}
\label{Transp1}
\Xeul_{\gamma_j}(W_{ij},w_{ij}) = \Yeul_j^{n_i} \\
\label{Transp2}
\text{$W_{ij}(\gamma) = W_i(\gamma)$
for all $\gamma \in P(S)\setminus\{\gamma_j, \gamma_j^- \}$}\\
\label{Transp3}
\text{$w_{ij}(v) = w_i(v)$
for all $v\in \Vtx_{>2}(S)\setminus \{ v_0^j, v_* \}$.}
\end{gather}
We claim that $ w_{12} = w_{21}$.
If $v \in \Vtx_{>2}(S) \setminus \{ v_0^1, v_0^2, v_* \}$
then, by \eqref{Transp3} and \eqref{iTransp3},
$w_{ij}(v) = w_i(v) = w(v)$ and hence $w_{12}(v)=w_{21}(v)$. 
By \eqref{Transp1},
$w_{ij}(v_*) = w_j(v_*) + n_i = w_j(v_*) + w_i(v_*) - w(v_*)$,
so $w_{12}(v_*) = w_{21}(v_*)$.
Finally, let $k\in\{1,2\}$ and assume that $v_0^k \in \Vtx_{>2}(S)$.
If $k=j$ (resp.\ $k=i$)
then by \eqref{Transp1} (resp.  by \eqref{Transp3})
we obtain $w_{ij}(v_0^k) = w_k(v_0^k)$; thus
$w_{12}(v_0^k) = w_{21}(v_0^k)$ and
we proved that $ w_{12} = w_{21}$.

Again, let $k\in\{1,2\}$.
If $k=j$ (resp.\ $k=i$)
then \eqref{Transp1} (resp. \eqref{Transp2}) implies that
$W_{ij}(\gamma_k) = Y_k$; it follows that 
$W_{12}(\gamma_k) = W_{21}(\gamma_k)$ and that 
$W_{12}(\gamma_k^-) = W_{21}(\gamma_k^-)$.
If $\gamma \in P(S) \setminus
\{ \gamma_1, \gamma_1^-, \gamma_2, \gamma_2^- \}$
then, by \eqref{Transp2} and \eqref{iTransp2},
$W_{ij}(\gamma) = W_{i}(\gamma) =  W(\gamma)$,
so $W_{12}(\gamma) = W_{21}(\gamma)$.
Hence, $W_{12} = W_{21}$ and we conclude that
$$
(\Theta_{y_2} \circ \Theta_{y_1}) (W, \pi(w))
=  (W_{12}, \pi(w_{12}))
=  (W_{21}, \pi(w_{21}))
= (\Theta_{y_1} \circ \Theta_{y_2}) (W, \pi(w)),
$$
which proves the second part of \eqref{ThetaCondition}.
\end{proof}

Since $\Reul^*(S,\omega)$ is not a group, the following remark is
relevant:

\begin{corollary}\label{Orbits}
If $\xi \in \omega\times N(S,\omega)$, define
$\Oeul_\xi = \setspec{y \xi }{ y \in \Reul^*(S,\omega) }$.
\begin{enumerate}

\item $\setspec{\Oeul_\xi }{ \xi \in \omega\times N(S,\omega) }$
is a partition of the set $\omega\times N(S,\omega)$.
The elements of this partition are called the ``orbits'' of the
left-action \ref{DefLeftAction}.

\item If $\xi, \xi'$ belong to the same orbit, then there exists
$y \in\Reul^*(S,\omega)$ such that $y \xi = \xi'$.

\end{enumerate}
\end{corollary}

\begin{proof}
This reduces to showing that,
if $y \in \Reul^*(S,\omega)$ and $\xi \in \omega\times N(S,\omega)$,
then there exists $y'\in\Reul^*(S,\omega)$
such that $y' y \xi = \xi$.

It suffices to prove this in the case where
$y = \binom{Y}{\gamma} \in \Reul(S,\omega)$.
In this case we define $y' = \binom{W(\gamma)}{\gamma}$,
where $W$ is defined by $\xi = (W,\eta)$.
Then $y,y' \in \Reul(S,\omega)$ are not disjoint and,
by \ref{AdvPptiesAction}, $y'y\xi = y'\xi$.
We have $y'\xi = \xi$ by the first part of \ref{BasicPptiesAction},
so we are done.
\end{proof}

Note the following consequence of part~(3) of \ref{BasicPptiesAction}:

\begin{nothing}\label{ObsSelf}
Consider $y = y_1 \cdots y_n \in \Reul^*(S,\omega)$,
where for each $i$ we have $y_i = \binom{Y_i}{\gamma_i} \in \Reul(S,\omega)$.
Let $(W',\eta') = y \xi$, where $\xi$ is any element of
$\omega \times N(S,\omega)$.
If $y$ is self-disjoint, then $W'(\gamma_i) = Y_i$ for
all $i \in \{ 1, \dots, n \}$.
\end{nothing}

\begin{corollary}\label{OrbitBij}
Let $\Oeul$ be an orbit of the left-action defined in \ref{DefLeftAction}.
Then the map
$$
\begin{array}{rcl}
\Oeul & \overset{}{\longrightarrow} & \omega \\
(W,\eta) & \longmapsto & W
\end{array}
$$
is bijective.
\end{corollary}

\begin{proof}
Let the map $(W,\eta) \mapsto W$ be denoted $p : \Oeul \to \omega$.
Suppose that $\xi=(W,\eta)$ and $\xi'=(W',\eta')$ belong to $\Oeul$.
By \ref{Orbits}, $\xi' = y \xi$ for some $y \in \Reul^*(S,\omega)$;
by \ref{AdvPptiesAction}, we may choose $y$ to be a self-disjoint word.
Write $y = y_1\cdots y_n$ with 
$y_i = \binom{Y_i}{\gamma_i} \in \Reul(S,\omega)$.
Since $y$ is self-disjoint, \ref{ObsSelf} gives
$W'(\gamma_i)=Y_i$ for all $i\in\{ 1, \dots, n\}$.
If $p(\xi)=p(\xi')$, then $W'=W$ and consequently
\begin{equation}\label{WordY}
\textstyle
y = \binom{W(\gamma_1)}{\gamma_1}\cdots \binom{W(\gamma_n)}{\gamma_n};
\end{equation}
then $y \xi = \xi$ by part~(1) of \ref{BasicPptiesAction},
showing that $p$ is injective.

To prove surjectivity,
select a subset $\Sigma = \{ \gamma_1, \dots, \gamma_n \}$ of
$P(S)$ satisfying
\begin{equation}\label{chooseOrientation}
\forall_{\gamma\in P(S)}\ \ \{ \gamma, \gamma^- \} \cap \Sigma
\text{ is a singleton}.
\end{equation}
Given $W \in \omega$,
define a word $y\in\Reul^*(S,\omega)$ by \eqref{WordY} and
note that $y$ is self-disjoint by \eqref{chooseOrientation}.
This and \ref{ObsSelf} imply that
if $\xi$ is an arbitrary element of $\Oeul$
then $y \xi = (W,\eta)$ for some $\eta$, i.e., $p(y \xi) = W$.
So $p$ is surjective.
\end{proof}

\begin{nothing*}
The following remarks are useful for computing
elements of $N(S,\omega)$.
\end{nothing*}

\begin{subnothing*}\label{N-Interp}
Let $C_1$, \dots, $C_p$ be the distinct connected components of
$(S,\omega)^\sharp$ which contain no special vertex.
Then $N(S,\omega)$ may be identified with the free module
$\Integ^{\{ C_1, \dots, C_p \}} = \Integ^p$.
Indeed, given any map $w : \Vtx_{>2}(S) \to \Integ$ define 
$\bar w : \{ C_1, \dots, C_p \} \to \Integ$
by $\bar w(C_i) = \sum_{v\in C_i} w(v)$.
Then 
$\Integ^{\Vtx_{>2}(S)} \to \Integ^{\{ C_1, \dots, C_p \}}$,
$w\mapsto\bar w$, is surjective and has kernel $Z(S,\omega)$.
\end{subnothing*}

\begin{subdefinition}
Given a connected component $C$ of $(S,\omega)^\sharp$
which contains no special vertex, define
\begin{equation*}
P_C^1(S) =
\setspec{ (u,v) \in P(S) }{ u \in C \text{ and } v \not\in C },\quad
P_C^0(S) =
\setspec{ (u,v) \in P(S) }{ u, v \in C }.
\end{equation*}
See \ref{DefdeltaXY} for $\delta(X,Y)$;
given $(W,W') \in \omega \times \omega$, define
$$
\delta_C(W,W') =
\sum_{\gamma \in P_C^1(S)} \delta( W(\gamma), W'(\gamma) )
\ +\  \frac 12
\sum_{\gamma \in P_C^0(S)} \delta( W(\gamma), W'(\gamma) ).
$$
Note that $\delta_C(W,W')$ is an integer. 
Indeed, if $\gamma \in P_C^0(S)$ then
$\gamma^- \in P_C^0(S)$ and $\det W(\gamma)=0$;
so  $\delta( W(\gamma), W'(\gamma) ) = \delta( W(\gamma^-), W'(\gamma^-) )$
and $\sum_{\gamma \in P_C^0(S)} \delta( W(\gamma), W'(\gamma) )$ is even.
Also,  if $W'' \in \omega$ then
\begin{equation}\label{EquWW'W''}
\delta_C(W,W') + \delta_C(W',W'') = \delta_C(W,W'').
\end{equation}
\end{subdefinition}

\begin{sublemma}\label{EtaPrimeEta}
Let $\Oeul \subseteq \omega \times N(S,\omega)$ be an orbit of
the left-action \ref{DefLeftAction}.
If $(W,\eta), (W', \eta') \in \Oeul$, then
\begin{equation}\label{EqEtaPrimeEta}
\eta'(C_i) = \eta(C_i) + \delta_{C_i}( W, W')\qquad (1 \le i \le p)
\end{equation}
where we view $\eta$ and $\eta'$ as maps $\{ C_1, \dots, C_p \}\to\Integ$,
as in \ref{N-Interp}.
\end{sublemma}

\begin{proof}
Since the map $\Oeul\to\omega$ of \ref{OrbitBij} is injective,
it is a priori clear that $\eta'$ is uniquely determined by
$\eta$, $W$ and $W'$.
We have $(W',\eta') = y (W,\eta)$ for some
$y = y_n\cdots y_1 \in \Reul^*(S,\omega)$ with 
$y_j \in \Reul(S,\omega)$.
Formula \eqref{EqEtaPrimeEta} holds trivially if $n=0$;
the case $n=1$ follows from part~(3) of \ref{BasicPptiesAction}
and \ref{ZeroProp}. If $n>1$ write 
$(W_0,\eta_0) = (W,\eta)$ and, 
for $0<j\le n$, $(W_j,\eta_j) = y_j (W_{j-1},\eta_{j-1})$;
then $\eta_j(C_i) - \eta_{j-1}(C_i) = \delta_{C_i}(W_{j-1},W_j)$
holds for every $j$, so
$ \eta'(C_i) - \eta(C_i)
= \sum_{j=1}^n (\eta_j(C_i) - \eta_{j-1}(C_i) )
= \sum_{j=1}^n \delta_{C_i}(W_{j-1},W_j)$,
and this is equal to $\delta_{C_i}(W,W')$ by \eqref{EquWW'W''}.
\end{proof}

\begin{definition}
An edge map $W$ for $S$ is \textit{canonical\/} if 
$$
\forall_{\gamma \in P(S)}\ \ 
\text{at least one of $W(\gamma)$, $W(\gamma^-)$ is a canonical sequence}
$$
(see \ref{DefCanSeq} for the notion of canonical sequence).
%
\end{definition}

\begin{lemma}\label{DeltaCan}
Let $W,W'$ be canonical elements of $\omega$ and let $C$ be a 
connected component of $(S,\omega)^\sharp$
which contains no special vertex. Then $\delta_C(W,W')=0$.
\end{lemma}

\begin{proof}
It suffices to show that if $X,Y \in \Integ^*$ are equivalent
sequences satisfying:
$$
\text{at least one of $X$, $X^-$ and one of $Y$, $Y^-$
is a canonical sequence,}
$$
then $\delta(X,Y)=0$.
Indeed, if this is true then
$\delta( W(\gamma), W'(\gamma) ) = 0$ for all $\gamma\in P(S)$,
and consequently $\delta_C(W,W')=0$.
We may assume that $X \neq Y$, otherwise $\delta(X,Y)=0$ holds trivially.
This implies that $\det X \neq 0$ and 
$\{ X, Y \} = \big\{ (0^{2i},A), (A, 0^{2i}) \big\}$, where $i>0$
and $\emptyseq \neq A \in \Neul^*$.
Since $\delta(Y,X) = - \delta(X,Y)$, we may
assume that $X = (0^{2i},A)$ and $Y = (A, 0^{2i})$.
Then one can obtain $\delta(X,Y) = 0$ by a direct calculation.
Alternatively, observe that the equivalence
$(0,0^{2i},A,0) \sim (0, A, 0^{2i},0)$ given by
\ref{MoveZeros} is in fact a $(+,+)$-equivalence;
thus $(0,X,0) \tequiv{(+,+)} (0,Y,0)$ and we obtain
$\delta(X,Y) = 0$ by part~(1a) of \ref{ZeroProp}.
\end{proof}

Recall from \ref{SeqThm} that each equivalence class of sequences in
$\Integ^*$ contains exactly one canonical sequence; so
$\setspec{W \in \omega}{ \text{$W$ is canonical} }$
is nonempty and finite.
So \ref{OrbitBij} implies:

\begin{corollary}\label{FiniteNonEmpty}
Let $\Oeul$ be an orbit of the left-action defined in \ref{DefLeftAction}.
Then the set
$\setspec{(W,\eta) \in \Oeul}{ \text{$W$ is canonical} }$
is finite and nonempty.
\end{corollary}

\begin{lemma}\label{AbsolLemma}
Let $(W,\eta), (W',\eta') \in \omega \times N(S,\omega)$
where $W$ and $W'$ are canonical edge maps. Then the following
are equivalent:
\begin{enumerate}

\item $(W,\eta)$ and $(W',\eta')$ belong to the same orbit
$\Oeul \subseteq \omega \times N(S,\omega)$

\item $\eta=\eta'$.

\end{enumerate}
\end{lemma}

\begin{proof}
The fact that (1) implies (2) is an immediate consequence of
\ref{EtaPrimeEta} and \ref{DeltaCan}.
To prove the converse,
assume that $(W,\eta), (W',\eta) \in \omega \times N(S,\omega)$
are such that $W$ and $W'$ are canonical.
By \ref{OrbitBij}, the orbit $\Oeul$ of $(W,\eta)$ contains an element
$(W',\eta^*)$ for some $\eta^*\in N(S,\omega)$. By applying the fact
that (1) implies (2) to the pairs $(W,\eta)$ and $(W',\eta^*)$,
we obtain $\eta=\eta^*$; so $(W',\eta)\in \Oeul$, as desired.
\end{proof}


\begin{corollary}\label{BijPi2}
There exists a unique map
$$
\pi_2\ :\ 
\setspec{\Oeul_\xi }{ \xi \in \omega\times N(S,\omega) }
\longrightarrow N(S,\omega)
$$
satisfying the following condition:
\begin{equation}
\tag{$*$} \text{If $\xi = (W,\eta) \in \omega \times N(S,\omega)$
and $W$ is a canonical edge map, $\pi_2 (\Oeul_\xi)=\eta$.}
\end{equation}
Moreover, $\pi_2$ is bijective.
\end{corollary}

\begin{proof}
That $(*)$ determines a unique set map is clear from
\ref{FiniteNonEmpty} and \ref{AbsolLemma}.
In fact, \ref{AbsolLemma} also implies that
$\pi_2$ is injective. If $\eta$ is any element of $N(S,\omega)$ then
pick a canonical edge map $W\in\omega$ and set $\xi = (W,\eta)$;
then $\pi_2( \Oeul_\xi ) = \eta$,
so $\pi_2$ is also surjective.
\end{proof}

\begin{proposition}\label{Equ=action}
Consider the map $\bar T : \For(S,\omega) \to \omega \times N(S,\omega)$
defined in \ref{Correspondence}.
If $\Ceul \subseteq \For(S,\omega)$ is an equivalence class then
$\bar T$ maps $\Ceul$ onto an orbit
$\Oeul \subseteq \omega \times N(S,\omega)$ of the
action~\ref{DefLeftAction}.
Moreover, $\bar T^{-1}(\Oeul) = \Ceul$.
\end{proposition}


\begin{proof}
Let $(\sigma, \GG), (\sigma', \GG') \in \For(S,\omega)$;
it suffices to show that the following are equivalent:
\begin{gather}
\label{EquGraphs}
(\sigma, \GG) \sim (\sigma', \GG')\\
\label{SameOrb}
\text{$y \bar T(\sigma,\GG) = \bar T(\sigma',\GG')$, 
for some $y\in\Reul^*(S,\omega)$}.
\end{gather}
Let $(W,w) = T(\sigma,\GG)$ and $(W',w') = T(\sigma',\GG')$.
Let $\pi : \Integ^{\Vtx_{>2}(S)} \to N(S,\omega)$ be the canonical
epimorphism.

If $(\sigma', \GG')$ is either a blowing-up or a blowing-down
of $(\sigma, \GG)$ then it is quite clear that 
$(\sigma', \GG')$ can be obtained by transplanting a suitable
$(\gamma,\Yeul)$ into $(\sigma, \GG)$, so 
$(W',w') = \transp(\gamma, \Yeul; W,w)$;
let $Y = \Yeul |_\tau$ (where $\tau$ is the type of $\gamma$)
and $y = \binom Y{\gamma} \in \Reul(S,\omega)$,
then $\Theta_y (W,\pi(w)) = (W',\pi(w'))$ and hence
\eqref{SameOrb} holds.
It follows that \eqref{EquGraphs} implies \eqref{SameOrb}.

Conversely, suppose that \eqref{SameOrb} holds.

If $y$ is the empty word then $(\sigma,\GG) \equiv (\sigma',\GG')$,
so \eqref{EquGraphs} holds by \ref{Cong=>Equiv}.

If $y = \binom Y{\gamma} \in \Reul(S,\omega)$ then, since
\eqref{SameOrb} holds,
we have
$\Theta_y( W, \pi(w) ) = (W', \pi(w'))$.
Choose $\Yeul \in \Integ^*_\tau$ (where $\tau$ is the type of $\gamma$)
such that $\Yeul |_\tau = Y$ and $\Yeul \tequiv{\tau} \Xeul_\gamma(W,w)$;
let
\begin{equation}\label{WwWw1}
(W^*,w^*) = \transp(\gamma, \Yeul; W,w).
\end{equation}
Then by definition of $\Theta_y$ we have
$\Theta_y( W, \pi(w) ) = (W^*, \pi(w^*))$, so
\begin{equation}\label{WwWw2}
(W', \pi(w')) = (W^*, \pi(w^*)).
\end{equation}
Let $(\sigma^*, \GG^*) \in \For(S,\omega)$ be such that 
$T(\sigma^*, \GG^*) = (W^*,w^*)$.
Then \eqref{WwWw1} means that 
$(\sigma^*, \GG^*)$ is obtained by transplanting $(\gamma,\Yeul)$ into
$(\sigma,\GG)$, so 
$(\sigma, \GG) \sim (\sigma^*, \GG^*)$ by \ref{Transp>Equiv};
and \eqref{WwWw2} means that $(\sigma',\GG') \equiv (\sigma^*, \GG^*)$,
so $(\sigma', \GG') \sim (\sigma^*, \GG^*)$ by \ref{Cong=>Equiv}.
So \eqref{EquGraphs} holds.

By induction on the length of the word $y$, it follows that
\eqref{SameOrb} implies \eqref{EquGraphs}.
\end{proof}

\section*{Solution to Problem~\ref{ProbForS}}

Recall that $S$ is fixed (but, from here on, $\omega$ is no longer fixed).

\begin{definition}
A weighted forest $\GG$ is \textit{canonical\/} if 
$$
\forall_{\gamma \in P(\GG)}\ \ 
\text{one of $W_\GG (\gamma)$, $W_\GG (\gamma^-)$ is a canonical sequence}
$$
(note that this condition implies that $\GG$ is minimal).
\end{definition}

\begin{lemma}\label{ClassContainsCan}
If $(\sigma,\GG) \in \For(S)$ then there exists
$(\sigma',\GG')\in \For(S)$ such that $(\sigma,\GG) \sim (\sigma',\GG')$
and $\GG'$ is a canonical forest.
\end{lemma}

\begin{proof}
Let $\omega\in\Omega(S)$ be such that $(\sigma,\GG)\in\For(S,\omega)$;
by \ref{Equ=action}, $\bar T : \For(S,\omega) \to \omega\times N(S,\omega)$
maps the equivalence class $\Ceul \subseteq \For(S,\omega)$ of
$(\sigma,\GG)$ onto an orbit $\Oeul \subseteq \omega\times N(S,\omega)$.
So the composite
$ \Ceul \xrightarrow{\ \bar T\ } \Oeul \xrightarrow{\ \ p\ \ } \omega$
(where $p$ is the bijection of \ref{OrbitBij}) is surjective
and some $(\sigma',\GG')\in\Ceul$ gets mapped to a canonical element of
$\omega$. Then $\GG'$ is canonical.
\end{proof}

\begin{definition}
Let $\Xgoth(S)$ be the set of pairs $(\omega,\eta)$ such that
$\omega\in\Omega(S)$ and $\eta \in N(S,\omega)$.
\end{definition}

Here is the solution to Problem~\ref{ProbForS}:

\begin{theorem}\label{MainForS}
There exists a unique surjective map
$$
Q : \For(S) \longrightarrow \Xgoth(S)
$$
which satisfies the following conditions:
\begin{enumerate}

\item For $(\sigma,\GG),(\sigma',\GG') \in \For(S)$,
$(\sigma,\GG) \sim (\sigma',\GG') \iff Q(\sigma,\GG)=Q(\sigma',\GG')$.

\item For any $(\sigma,\GG)\in \For(S)$ such that
$\GG$ is a canonical forest, 
let $\omega$ be the element of $\Omega(S)$ such that 
$(\sigma,\GG)\in \For(S,\omega)$
and let $\eta\in N(S,\omega)$ be such that $\bar T(\sigma,\GG) = (W,\eta)$
for some $W$;
then $Q(\sigma,\GG) = (\omega,\eta)$.

\end{enumerate}
\end{theorem}

\begin{proof}
In view of \ref{ClassContainsCan}, it is clear that $Q$ is completely
determined by conditions (1) and (2).  So it suffices to prove the
existence of $Q$.

For each $\omega \in \Omega(S)$, let 
$Q_\omega : \For(S,\omega) \to \Xgoth(S)$
be the composite
$$
\begin{array}{ccccccc}
\For(S,\omega) & \xrightarrow{\ \alpha\ } &
\setspec{\Oeul_\xi }{ \xi \in \omega\times N(S,\omega) } &
\xrightarrow{\ \beta\ } & \{ \omega \} \times N(S,\omega) &
\hookrightarrow & \Xgoth(S) \\
(\sigma, \GG) & \longmapsto & \text{orbit of }\bar T(\sigma,\GG) \\
 & & \Oeul & \longmapsto & \big( \omega, \pi_2( \Oeul ) \big) 
\end{array}
$$
(see \ref{BijPi2} for $\pi_2$).
Since $\setspec{ \For(S,\omega) }{ \omega \in \Omega(S) }$ is a
partition of $\For(S)$, we obtain a map $Q : \For(S) \to \Xgoth(S)$
by taking the union of the $Q_\omega$.
Note that $\alpha$ is surjective, because 
$\bar T : \For(S,\omega) \to \omega\times N(S,\omega)$ is surjective;
since $\beta$ is bijective by \ref{BijPi2},
it follows that the image of $Q_\omega$ is 
$\{ \omega \} \times N(S,\omega)$.  Consequently, $Q$ is surjective.

If $(\omega, \eta) \in \{ \omega \} \times N(S,\omega)$ then, since
$\beta$ is bijective, $Q_\omega^{-1}(\omega, \eta)$ is equal to the
inverse image by $\bar T$ of some orbit;  by \ref{Equ=action},
this is an equivalence class in $\For(S,\omega)$ (hence also in $\For(S)$),
so $Q$ satisfies condition (1).
Because $\beta$ is defined in terms of $\pi_2$,
$Q$ satisfies condition (2)
(see condition $(*)$ in \ref{BijPi2}).
\end{proof}

\section*{Solution to Problem~\ref{ProbS}}

We consider the group $A=\Aut(S)$ of graph automorphisms of $S$.

\begin{definition}\label{DefAaction}
We define a right-action of the group $A$ on the set $\Xgoth(S)$.
\begin{enumerate}

\item If $W$ is an edge map for $S$ and $\alpha \in A$, let
$W^\alpha = W \circ \vec{\alpha} : P(S) \to \Integ^*$ and note that
$W^\alpha$ is an edge map for $S$; this defines a right-action of $A$
on the set of edge maps for $S$.

\item If $\omega\in\Omega(S)$ and $\alpha \in A$, let
$\omega^\alpha = \setspec{ W^\alpha }{ W\in\omega } \in \Omega(S)$.
This is a right-action of $A$ on the set $\Omega(S)$.

\item Again, let $\omega\in\Omega(S)$ and $\alpha \in A$.
Then the restriction $\alpha^\sharp : \Vtx_{>2}(S) \to \Vtx_{>2}(S)$
of $\alpha$ is an isomorphism of graphs,
$\alpha^\sharp : (S,\omega^\alpha)^\sharp \to (S,\omega)^\sharp$,
and maps the special vertices of 
$(S,\omega^\alpha)^\sharp$ to those of $(S,\omega)^\sharp$
(refer to \ref{DefN} for all this).
Consequently, the automorphism of $\Integ$-modules
$$
\Integ^{ \Vtx_{>2}(S) } \longrightarrow \Integ^{ \Vtx_{>2}(S) },\quad
w \longmapsto w \circ \alpha^\sharp
$$
maps $Z(S,\omega)$ onto $Z(S,\omega^\alpha)$ and hence determines
an isomorphism of $\Integ$-modules
$$
N(S,\omega) \to N(S,\omega^\alpha), \quad \eta \longmapsto \eta^\alpha.
$$

\item If $(\omega,\eta) \in\Xgoth(S)$ and $\alpha \in A$, let
$(\omega,\eta)^\alpha = (\omega^\alpha,\eta^\alpha) \in\Xgoth(S)$
where $\omega^\alpha$ and $\eta^\alpha$ are defined in parts (2) and (3)
of this definition. This defines a right-action of the group $A$
on the set $\Xgoth(S)$. The set of orbits is denoted $\Xgoth(S)/A$.

\end{enumerate}
\end{definition}

See \ref{DefRFSK} and \ref{DefFor} for the definitions of
$\RFSK(S)$, $\For(S)$ and $p_2: \For(S)\to\RFSK(S)$.
The following result solves Problem~\ref{ProbS}:

\begin{theorem}\label{MainRFSK}
There exists a unique map
$\bar Q : \RFSK(S) \longrightarrow \Xgoth(S) / A$
such that
\begin{equation}\label{CDQbar}
\begin{CD}
\For(S) @>{Q}>> \Xgoth(S) \\
@VV{p_2}V  @VVV \\
\RFSK(S) @>{\overline Q}>> \Xgoth(S)/A
\end{CD}
\end{equation}
commutes, where $Q$ is defined in \ref{MainForS} and where
$\Xgoth(S) \to \Xgoth(S)/A$ is the canonical quotient map.
Moreover, $\bar Q$ is surjective and given any
$\GG,\GG' \in \RFSK(S)$ we have
$$
\GG \sim \GG' \iff \bar Q(\GG)=\bar Q(\GG').
$$
\end{theorem}

Some facts are needed for the proof of \ref{MainRFSK}.

\begin{definition}\label{DefRightAction}
Given $(\sigma,\GG) \in \For(S)$ and $\alpha \in A$,
let $(\sigma,\GG)^\alpha = (\sigma \circ \alpha, \GG) \in \For(S)$.
This defines a right-action of the group $A$ on the set $\For(S)$.
The set of orbits is denoted $\For(S)/A$.
\end{definition}

\begin{lemma}\label{BijRFSK}
If $\GG \in \RFSK(S)$, then $p_2^{-1}( \GG )$ is an orbit of the
right-action \ref{DefRightAction}.
Moreover, $\GG \mapsto p_2^{-1}( \GG )$ defines 
a bijection $\RFSK(S) \to \For(S)/A$.
\end{lemma}

\begin{proof}
By \ref{ExistsUniqueSkel},
the skeletal map $S \skm \GG$ is unique up to automorphism of $S$.
\end{proof}

\begin{lemma}\label{Aaction-T}
Let $(\sigma,\GG) \in \For(S)$ and $\alpha \in A$.
If $\bar T( \sigma, \GG ) = (W,\eta) \in \omega \times N(S,\omega)$,
then 
$\bar T\big( ( \sigma, \GG )^\alpha \big) 
= (W^\alpha ,\eta^\alpha ) \in \omega^\alpha  \times N(S,\omega^\alpha)$.
\end{lemma}

\begin{proof}
Let $(W,w) = T(\sigma,\GG)$, then 
$T\big( (\sigma, \GG)^\alpha \big) = 
T(\sigma\circ\alpha, \GG) = (W^\alpha, w \circ \alpha^\sharp )$
is clear and the result follows.
\end{proof}

\begin{lemma}\label{Aaction-equiv}
Given $(\sigma,\GG), (\sigma',\GG') \in \For(S)$ and $\alpha \in A$,
$$
(\sigma,\GG) \sim (\sigma',\GG') \iff
(\sigma,\GG)^\alpha \sim (\sigma',\GG')^\alpha.
$$
\end{lemma}

\begin{proof}
It suffices to prove that if 
$(\sigma',\GG')$ is a blowing-down of $(\sigma,\GG)$
then
$(\sigma',\GG')^\alpha$ is a blowing-down of $(\sigma,\GG)^\alpha$.
Let $\pi : \GG \skm \GG'$ be the blowing-down map;
then $\sigma' = \pi \circ \sigma$ and consequently
$\sigma' \circ \alpha = \pi \circ (\sigma \circ \alpha)$,
so
$(\sigma' \circ \alpha ,\GG') = (\sigma',\GG')^\alpha$
is a blowing-down of
$(\sigma \circ \alpha ,\GG) = (\sigma,\GG)^\alpha$.
\end{proof}

\begin{lemma}\label{Qcommutes}
$Q\big( (\sigma,\GG)^\alpha \big) = Q(\sigma,\GG)^\alpha$,
for all $(\sigma,\GG) \in \For(S)$ and $\alpha\in A$.
\end{lemma}

\begin{proof}
By \ref{ClassContainsCan},
we may choose $(\sigma', \GG') \in \For(S)$ such that 
$(\sigma', \GG') \sim (\sigma, \GG)$ and 
$\GG'$ is a canonical forest.
Then $(\sigma, \GG) \sim (\sigma', \GG')$ implies
$Q(\sigma, \GG) = Q(\sigma', \GG')$ and hence
$ Q(\sigma, \GG)^\alpha = Q(\sigma', \GG')^\alpha $.
By \ref{Aaction-equiv} we have
$(\sigma,\GG)^\alpha \sim (\sigma',\GG')^\alpha$, so
$ Q\big( (\sigma, \GG)^\alpha \big) = Q\big( (\sigma', \GG')^\alpha \big)$.
So it suffices to prove that 
$Q\big( (\sigma',\GG')^\alpha \big) = Q(\sigma',\GG')^\alpha$,
i.e., we may assume that $\GG$ is canonical.

Let $\omega$ be the element of $\Omega(S)$ such that 
$(\sigma,\GG)\in \For(S,\omega)$ and 
let $(W,\eta) = \bar T( \sigma, \GG ) \in \omega \times N(S,\omega)$;
then \ref{Aaction-T} gives
$\bar T\big( ( \sigma, \GG )^\alpha \big) =
(W^\alpha ,\eta^\alpha ) \in \omega^\alpha \times N(S,\omega^\alpha)$.
Since $\GG$ is canonical,
$Q( \sigma, \GG)$ and 
$Q \big( (\sigma, \GG)^\alpha \big) = Q (\sigma\circ\alpha, \GG)$
may be computed by applying part~(2) of \ref{MainForS};
this gives $Q( \sigma, \GG) = ( \omega, \eta)$ and
$Q \big( (\sigma, \GG)^\alpha \big) =  ( \omega^\alpha, \eta^\alpha)$,
which is the desired statement.
\end{proof}

\begin{proposition}\label{EquEqu}
Let $(\sigma,\GG) \in \For(S)$.
For a pseudo-minimal forest $\GG'$, the following are equivalent:
\begin{enumerate}

\item $\GG \sim \GG'$
\item There exists $\sigma' : S \skm \GG'$ such that
$(\sigma,\GG) \sim (\sigma',\GG')$ in $\For(S)$.

\end{enumerate}
\end{proposition}

\begin{proof}
Only $(1) \Rightarrow (2)$ requires a proof.
Recall from \ref{AllReg} that
if (1) holds then $\GG$ is strictly equivalent to $\GG'$.
So we may assume that $\GG'$ is 
a strict blowing-down or a strict blowing-up of $\GG$.

Suppose that $\GG'$ is a strict blowing-down of $\GG$.
Let $\pi : \GG \skm \GG'$ be the blowing-down map
defined in \ref{DefBlowDownMap} and let $\sigma' = \pi \circ \sigma$;
then it is immediate that $(\sigma',\GG')$ belongs to  $\For(S)$
and is a blowing-down of $(\sigma,\GG)$.

If $\GG'$ is a strict blowing-up of $\GG$ then, by \ref{Basic:Skel},
$\sigma$ factors through the blowing-down map $\pi : \GG' \skm \GG$
and consequently
there exists a skeletal map $\sigma' : S \skm \GG'$ such that
$(\sigma',\GG')$ is a blowing-up of $(\sigma,\GG)$.
\end{proof}

\bigskip
\begin{proof}[Proof of \ref{MainRFSK}]
By \ref{BijRFSK} and \ref{Qcommutes}, if two elements of $\For(S)$ 
have the same image via $p_2$ then they have the same image via
the composite $\For(S) \xrightarrow{Q} \Xgoth(S)\to \Xgoth(S)/A$;
thus there exists a unique map $\bar Q$ such that \eqref{CDQbar} is
commutative. Since $\Xgoth(S)\to \Xgoth(S)/A$ and $Q$ are surjective,
so is $\bar Q$.
Let $\GG,\GG'\in\RFSK(S)$.
Pick any skeletal maps $\sigma:S \skm \GG$ and $\sigma':S\skm \GG'$
and consider $(\sigma,\GG), (\sigma',\GG')\in \For(S)$.
Then
\begin{align*}
\bar Q(\GG)=\bar Q(\GG')
&\overset{\eqref{CDQbar}}{\iff}
	\text{$Q(\sigma,\GG)$ and $Q(\sigma',\GG')$ belong
	to the same $A$-orbit}\\
&\iff
	\exists_{\alpha\in A}\ Q(\sigma,\GG)
= Q(\sigma',\GG')^\alpha
\overset{\ref{Qcommutes}}=  Q \big( (\sigma',\GG')^\alpha \big) \\
&\overset{\ref{MainForS}}{\iff}
	\exists_{\alpha\in A}\ (\sigma,\GG) \sim (\sigma',\GG')^\alpha \\
&\overset{\ref{BijRFSK}}{\iff}
	\exists_{\sigma^*:S\skm \GG'}\ (\sigma,\GG) \sim (\sigma^*,\GG')
\overset{\ref{EquEqu}}{\iff}
	\GG \sim \GG',
\end{align*}
which completes the proof of \ref{MainRFSK}.
\end{proof}

We conclude this section with:

\begin{problem}\label{DecisionProblem}
Given weighted forests $\GG$ and $\GG'$, decide whether they are equivalent.
\end{problem}
\noindent\textbf{Solution:}
\begin{enumerate}\it

\item Blow-down the two graphs until they are minimal; 
check that the two minimal weighted forests have isomorphic
skeletons.
\end{enumerate}
From now-on, assume that the skeletons are isomorphic
(by \ref{SameSkel}, this is a necessary condition for equivalence).
\begin{enumerate}\it
\addtocounter{enumi}{1}

\item Find canonical weighted forests $\GG_*$ and $\GG'_*$ such that
$\GG \sim \GG_*$ and $\GG' \sim \GG'_*$.

\end{enumerate}
One can design an algorithm which accomplishes step~(2);
we leave this to the reader.
\begin{enumerate}\it
\addtocounter{enumi}{2}

\item Let $S$ be the skeleton of $\GG_*$ (and of $\GG'_*$),
fix a skeletal map $\sigma : S \skm \GG_*$ and compute
$(W,\eta) = \bar T(\sigma,\GG_*)$.
For each 
$\sigma' : S \skm \GG'_*$, compute $(W',\eta') = \bar T(\sigma',\GG'_*)$.
The condition $\GG \sim \GG'$ is equivalent to:
$$
\text{For some $\sigma'$ such that $W'\sim W$, we have $\eta' = \eta$.}
$$

\end{enumerate}

The claim contained in step~(3) is a consequence of \ref{MainRFSK}.
Note that if $\sigma'$ satisfies $W' \sim W$
then $\eta$ and $\eta'$ belong to the same module $N(S,\omega)$
(where $\omega$ is the equivalence class of $W$) and hence can be compared.
One can use \ref{N-Interp} to compare $\eta$ and $\eta'$.

By \ref{BijRFSK}, the number of skeletal maps
$\sigma' : S \skm \GG'_*$ is equal to the order
of the group $A = \Aut(S)$,
which is often a small number (\ref{carterpillar}).
Moreover, we are only interested in those $\sigma'$ which satisfy
$W' \sim W$; if $\sigma_0'$ is such a map, then the set of $\sigma'$
satisfying $W' \sim W$ is
$\setspec{ \sigma_0' \circ \alpha }{ \alpha \in A \text{ and }
\omega^\alpha = \omega}$
(see \ref{DefAaction} for the right-action of $A$ on $\Omega(S)$).

\begin{example}\label{carterpillar}
In the study of algebraic surfaces, 
skeletons of the form
$$
S\ \ =\ \ \raisebox{-6\unitlength}
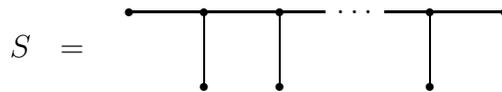
{%
\begin{picture}(54,14)(-2,-12)
\put(0,0){\circle*{1}}
\put(10,0){\circle*{1}}
\put(20,0){\circle*{1}}
\put(30,0){\makebox(0,0){\mbox{\,\dots}}}
\put(40,0){\circle*{1}}
\put(50,0){\circle*{1}}
\put(0,0){\line(1,0){26}}
\put(50,0){\line(-1,0){16}}
\put(10,-10){\circle*{1}}
\put(20,-10){\circle*{1}}
\put(40,-10){\circle*{1}}
\put(10,-10){\line(0,1){10}}
\put(20,-10){\line(0,1){10}}
\put(40,-10){\line(0,1){10}}
%
%
\end{picture}}
$$
are not uncommon. 
Such an $S$ satisfies $|\Aut(S)|=8$ if at least two vertices have degree $3$.
\end{example}

\section{Minimal weighted forests}
\label{Sec:MinimalReduc}

This section reduces Problem~\ref{MinWG} to Problem~\ref{MinSQ}.
Recall that Section~\ref{Sec:FurtherLin} partially solves
Problem~\ref{MinSQ}.

\begin{definition}
Let $S$ be a skeleton and $\omega \in \Omega(S)$.
By a \textit{minimal element of} $\omega$, we mean an element $W\in\omega$
such that
$$
\forall_{\gamma \in P(S)}\ \text{$W(\gamma)$ is a minimal element of
$\omega(\gamma)$}
$$
(see \ref{DefSWw} for the element $\omega(\gamma)$ of $\SEC$).
The symbol $\min(\omega)$ denotes the set of minimal elements of $\omega$.
\end{definition}

\begin{subnothing*}\label{FirstRed}
Pick a subset $\Sigma = \{ \gamma_1, \dots, \gamma_n \}$ of $P(S)$
such that $\Sigma \cap \{ \gamma, \gamma^-\}$ is a singleton 
for every $\gamma\in P(S)$.
Write $\Ceul_i = \omega(\gamma_i)$ for each $i$.
Then $W \mapsto (W(\gamma_1), \dots, W(\gamma_n))$ is a bijection
$\omega \to \prod_{i=1}^n \Ceul_i$ and maps $\min(\omega)$ onto
$\prod_{i=1}^n \min(\Ceul_i)$.
Consequently, a description of $\min(\omega)$ follows immediately from
a solution to Problem~\ref{MinSQ}.
\end{subnothing*}

\begin{nothing*}
Suppose that we want to describe the set $\Meul$ of minimal elements of an
equivalence class $\Ceul$ of weighted forests.

Pick $\GG \in \Meul$.
We may choose a skeleton $S$ and a skeletal map $\sigma : S \skm \GG$;
then $(\sigma,\GG) \in \For(S)$ and we may consider
the equivalence class $\Ceul^+ \subseteq \For(S)$ of $(\sigma,\GG)$.
By \ref{EquEqu},
$p_2(\Ceul^+)$ is the set of pseudo-minimal forests belonging to $\Ceul$;
so 
\begin{equation}\label{TwoInclusions}
\Meul \subseteq p_2(\Ceul^+) \subseteq \Ceul.
\end{equation}
By \ref{SameOmega},
we may choose $\omega \in \Omega(S)$ such that
$\Ceul^+ \subseteq \For(S,\omega)$.
Then, by \ref{Equ=action},
the map $\bar T :  \For(S,\omega) \to \omega \times N(S,\omega)$
maps $\Ceul^+$ onto an orbit $\Oeul \subseteq \omega \times N(S,\omega)$;
also consider the bijection $p : \Oeul \to \omega$ of \ref{OrbitBij},
so we have:
\begin{equation}\label{3maps}
\Ceul \xleftarrow{\ p_2\ } 
\Ceul^+ \xrightarrow{\ \bar T\ } \Oeul
\xrightarrow{\ p\ } \omega.
\end{equation}
\end{nothing*}

\begin{subproposition}\label{TechRed}
$\Meul = p_2 (\Meul^+)$, where
$\Meul^+ \subseteq \Ceul^+$ denotes the inverse image of $\min(\omega)$
via $p\circ\bar T : \Ceul^+ \to \omega$.
\end{subproposition}

\begin{proof}
It is clear that, for an element $(\sigma',\GG')$ of $\Ceul^+$,
$\GG'$ is a minimal weighted graph if and only if 
$p(\bar T( \sigma',\GG') )$ is a minimal element of $\omega$;
in other words, $\Meul^+ = p_2^{-1}(\Meul)$.
We get $p_2(\Meul^+) = \Meul$ by \eqref{TwoInclusions}.
\end{proof}

\begin{remark}
In applications of \ref{TechRed}, one needs to evaluate
$p^{-1}: \omega \to \Oeul$ at some elements of $\omega$;
\ref{EtaPrimeEta} can be used for this.
\end{remark}

\begin{example}
Fix $a,b,c \in \Integ$ and consider the weighted tree
$$
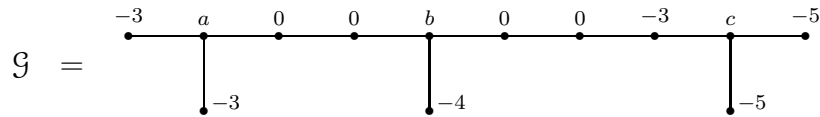
\GG\ \ =\ \ \raisebox{-6\unitlength}{%
\begin{picture}(94,14)(-2,-11)
\multiput(0,0)(10,0){10}{\circle*{1}}
\put(10,-10){\circle*{1}}
\put(40,-10){\circle*{1}}
\put(80,-10){\circle*{1}}
\put(0,0){\line(1,0){90}}
\put(10,-10){\line(0,1){10}}
\put(40,-10){\line(0,1){10}}
\put(80,-10){\line(0,1){10}}
\put(0,1.5){\makebox(0,0)[b]{$\scriptstyle -3$}}
\put(10,1.5){\makebox(0,0)[b]{$\scriptstyle a$}}
\put(20,1.5){\makebox(0,0)[b]{$\scriptstyle 0$}}
\put(30,1.5){\makebox(0,0)[b]{$\scriptstyle 0$}}
\put(40,1.5){\makebox(0,0)[b]{$\scriptstyle b$}}
\put(50,1.5){\makebox(0,0)[b]{$\scriptstyle 0$}}
\put(60,1.5){\makebox(0,0)[b]{$\scriptstyle 0$}}
\put(70,1.5){\makebox(0,0)[b]{$\scriptstyle -3$}}
\put(80,1.5){\makebox(0,0)[b]{$\scriptstyle c$}}
\put(90,1.5){\makebox(0,0)[b]{$\scriptstyle -5$}}
\put(11,-9){\makebox(0,0)[l]{$\scriptstyle -3$}}
\put(41,-9){\makebox(0,0)[l]{$\scriptstyle -4$}}
\put(81,-9){\makebox(0,0)[l]{$\scriptstyle -5$}}
\end{picture}}
$$
Write $X_1 = (0,0)$ and $X_2 = (0,0,-3)$ and,
for each $i\in\{1,2\}$,
let $\Ceul_i \in \SEC$ be the equivalence class of $X_i$.
Given $(Y_1, Y_2) \in \min(\Ceul_1) \times \min(\Ceul_2)$,
set
\begin{equation}\label{a'b'c'}
a' = a + \delta(X_1,Y_1),
\quad b' = b + \delta(X_1^-,Y_1^-) + \delta(X_2,Y_2),
\quad c' = c + \delta(X_2^-,Y_2^-)
\end{equation}
and define
$$
\GG'(Y_1,Y_2)\ \ =\ \ \raisebox{-6\unitlength}{%
\begin{picture}(64,17)(-2,-11)
\put(0,0){\circle*{1}}
\put(10,0){\circle*{1}}
\put(30,0){\circle*{1}}
\put(50,0){\circle*{1}}
\put(60,0){\circle*{1}}
\put(10,-10){\circle*{1}}
\put(30,-10){\circle*{1}}
\put(50,-10){\circle*{1}}
\put(0,0){\line(1,0){16}}
\put(30,0){\line(1,0){6}}
\put(30,0){\line(-1,0){6}}
\put(60,0){\line(-1,0){16}}
\put(10,-10){\line(0,1){10}}
\put(30,-10){\line(0,1){10}}
\put(50,-10){\line(0,1){10}}
\put(0,1.5){\makebox(0,0)[b]{$\scriptstyle -3$}}
\put(10,1.5){\makebox(0,0)[b]{$\scriptstyle a'$}}
\put(30,1.5){\makebox(0,0)[b]{$\scriptstyle b'$}}
\put(50,1.5){\makebox(0,0)[b]{$\scriptstyle c'$}}
\put(60,1.5){\makebox(0,0)[b]{$\scriptstyle -5$}}
\put(11,-9){\makebox(0,0)[l]{$\scriptstyle -3$}}
\put(31,-9){\makebox(0,0)[l]{$\scriptstyle -4$}}
\put(51,-9){\makebox(0,0)[l]{$\scriptstyle -5$}}
\put(20,0){\makebox(0,0){\mbox{\,\dots}}}
\put(40,0){\makebox(0,0){\mbox{\,\dots}}}
\put(20,1){\makebox(0,0)[b]{%
$\displaystyle\overbrace{ \rule{12\unitlength}{0mm} }^{Y_1}$ }}
\put(40,1){\makebox(0,0)[b]{%
$\displaystyle\overbrace{ \rule{12\unitlength}{0mm} }^{Y_2}$ }}
\end{picture}}
$$
By \ref{TechRed}
(and noting that \eqref{a'b'c'} is obtained from \ref{EtaPrimeEta}),
it follows that
$$
\setspec{ \GG'(Y_1,Y_2) } {(Y_1, Y_2) \in \min(\Ceul_1) \times \min(\Ceul_2)}
$$
is the set of 
minimal weighted graphs equivalent to $\GG$.
Since the sets $\min(\Ceul_1)$ and $\min(\Ceul_2)$
are explicitely described in section \ref{Sec:FurtherLin}
(because each $\Ceul_i$ is the successor of a prime class),
and since \eqref{a'b'c'} is an explicit formula for
$a'$, $b'$, $c'$, 
we know all minimal weighted graphs equivalent to $\GG$.
\end{example}

\medskip
We conclude that \ref{TechRed} reduces Problem~\ref{MinWG}
to the problem of describing $\min(\omega)$.
Hence, in view of \ref{FirstRed}, 
Problem~\ref{MinWG} reduces to Problem~\ref{MinSQ}.

\bibliographystyle{amsplain}
\bibliography{/home/ddaigle/articles/bib/dbase}

\end{document}